\documentclass[11pt,reqno]{amsart}
\usepackage[T1]{fontenc}
\usepackage[utf8]{inputenc}
\usepackage{lmodern}
\usepackage{amsmath,amsthm,amssymb,mathtools,mathrsfs}
\usepackage[margin=1.05in]{geometry}
\usepackage{microtype,enumitem,booktabs,xcolor,needspace}
\usepackage{palatino,mathpazo}
\definecolor{paperlink}{RGB}{0,0,100}
\usepackage[colorlinks=true,linkcolor=paperlink,citecolor=paperlink,urlcolor=paperlink]{hyperref}
\hypersetup{pdftitle={Entire Ginzburg-Landau solutions beyond the standard vortices},pdfsubject={Brezis Problem 2.4}}
\numberwithin{equation}{section}
\newtheorem{theorem}{Theorem}[section]
\newtheorem{mainresult}[theorem]{Theorem}
\newtheorem{proposition}[theorem]{Proposition}
\newtheorem{lemma}[theorem]{Lemma}
\theoremstyle{definition}

\theoremstyle{remark}

\newcommand{\C}{\mathbb C}
\newcommand{\R}{\mathbb R}
\DeclareMathOperator{\ind}{ind}
\title{On Brezis' open problem 2.4}

\author[H.G. Chen]{Hong-Ge Chen}
\address{Hong-Ge Chen, School of Mathematics and Statistics, Key Laboratory of Nonlinear Analysis \& Applications (Ministry of Education),
Central China Normal University, Wuhan, China}
\email{hongge\_chen@whu.edu.cn}

\author[Y. Liu]{Yong Liu}
\address{Yong Liu, School of Mathematics and Statistics, Beijing Technology and Business University, Beijing, China}
\email{yliumath@btbu.edu.cn}

\author[J.C. Wei]{Juncheng Wei}
\address{Juncheng Wei, Department of Mathematics, Chinese University of Hong Kong, Shatin, New Territories, Hong Kong}
\email{wei@math.cuhk.edu.hk}

\author[W. Yang]{Wen Yang}
\address{Wen Yang, Department of Mathematics, Faculty of Science, University of Macau, Taipa, Macau, China}
\email{wenyang@um.edu.mo}

\date{} 
\begin{document}
\begin{abstract}
We show the existence of a ``nonradial'' smooth entire solution of the Ginzburg--Landau equation
with degree $q\in\{8,9\}$ and modulus tending to one at infinity. This answers Brezis' Open Problem~2.4 negatively.
To construct the counterexample, we vary a cone parameter $q \in [8, 9]$ and keeping the degree on the cone equal to one. Several key ingredients are needed. First, we need to compute the Morse index of radial vortex solutions at $ q\in \{ 8, 9\}$. Second, we need to show compactness. Finally, Leray-Schauder degree theory and a homotopy argument  yield the desired existence result.
\end{abstract}
\maketitle

\section{Introduction}\label{intro}

We study smooth complex valued solutions of the Ginzburg--Landau equation
\begin{equation}\label{equation}
 -\Delta u=u(1-|u|^2)\quad\text{in }\R^2,
 \qquad \lim_{|x|\to\infty}|u(x)|=1.
\end{equation} We denote by $\deg(u,\infty)$  the degree of $u$.
The simplest solutions of positive degree $q$ are
\[
 V_q(re^{i\theta})=f_q(r)e^{iq\theta},
 \qquad f_q(0)=0,\qquad \lim_{r\to\infty}f_q(r)=1.
\]
The function $f_q$ is positive and strictly increasing on $(0,\infty)$.
It solves the ordinary differential equation
\[
 f_q''+\frac1r f_q'-\frac{q^2}{r^2}f_q+f_q(1-f_q^2)=0.
\]
The existence, uniqueness, and asymptotic properties of this solution are
classical. We refer to Herv\'e--Herv\'e \cite{HH} (see also the book by
Pacard--Rivi\`ere \cite{PR} for more references and related results). For negative degrees we set
$V_{-q}=\overline{V_q}$, and we set $V_0=1$. It is easy to see that translations and multiplication by a constant of modulus one preserve
\eqref{equation}. They therefore generate the family
$\zeta V_q(\,\cdot-a\,)$, where $a\in\R^2$ and $|\zeta|=1$. In Open Problem~2.4, Brezis \cite[pp.~310--311]{Brezis2023} asks whether
every solution of \eqref{equation} belongs to this family.
This will be the main objective in this paper.  

The potential energy  plays  important roles in our analysis. By definition, the potential energy is 
\[
 \mathcal P(u)=\int_{\R^2}(1-|u|^2)^2\,dx.
\]
Brezis' Open Problem~2.5  asks if all solutions of \eqref{equation} have finite potential energy $ \mathcal P(u)<+\infty.$ This has been solved recently  by Chen-Wei-Yan-Yang \cite{CWYY}, see also \cite{LYYYZ}.

Recall that for a solution of nonzero degree, the full Ginzburg-Landau energy on $B_R(0)$ in general grows 
logarithmically as $R\to\infty$. Indeed, for the standard vortex, the well known expansion
$f_q(r)=1-q^2/(2r^2)+O(r^{-4})$ yields
\[
 \int_{B_R(0)}
 \left(\frac12|\nabla V_q|^2+\frac14(1-|V_q|^2)^2\right)\,dx
 =\pi q^2\log R+O(1).
\]
Brezis--Merle--Rivi\`ere \cite[Theorem~1 and Remarks~1.1--1.2]{BMR}
proved that every entire solution with finite potential energy satisfies
$|u(x)|\to1$ and the quantized identity
\[
 \int_{\R^2}(1-|u|^2)^2\,dx=2\pi (\deg(u,\infty))^2.
\]
In particular, degree zero implies that the solution is constant with
modulus one. This identity connects a global integral to an integer and will be one of the main ingredients in our proof.

The  degree one solution has much stronger rigidity and has been studied for a long time.
Mironescu \cite[Theorem~1]{Mironescu96} classified all entire solutions
of degree $\pm1$ with finite potential energy.
The linearized equation at the standard vortex is treated in
\cite[Chapter~3]{PR}, and
del Pino--Felmer--Kowalczyk \cite{DPFK} proved minimality and
nondegeneracy through an inequality of Hardy type.
More recently, del Pino--Juneman--Musso \cite{DPJM} obtained estimates
for the inverse linearized operator by explicit formulas for its angular
Fourier components. We also mention that Gravejat--Pacherie--Smets \cite{GPS} constructed
a functional setting in which a renormalized energy has the degree one
vortex as its unique global minimizer, up to translation and phase.
They also proved nonlinear coercivity near this orbit.
These results explain the rigidity of an isolated vortex of degree one
and also provide analytic tools for studying several separated vortices.

The situation becomes much more complicated if the solution has multiple vortices (zero points). For several vortices, their individual degrees and mutual positions
both enter the energy. This relation is central in the theory on
bounded domains. Bethuel--Brezis--H\'elein
\cite{BBH93} studied the asymptotics of minimizers with boundary data of
degree zero. Their monograph \cite{BBH} develops the analysis of isolated
vortices for nonzero boundary degree. For the energy on a bounded domain $\Omega$,
\[
 \int_\Omega
 \left(\frac12|\nabla u|^2+
       \frac{(1-|u|^2)^2}{4\varepsilon^2}\right)\,dx,
\]
the logarithmic contribution accounts for the individual vortices,
while the renormalized energy describes their interaction and their
relation to the boundary. The energy lower bounds of Jerrard
\cite{Jerrard} and Sandier \cite{Sandier} relate the topology of a map
to the energy required near its zeros.
Jerrard--Soner \cite{JS} used the distributional Jacobian to express
compactness and quantization for sequences with bounded normalized
energy. Sandier--Serfaty \cite{SSbounds} refined the vortex ball estimates to
obtain lower bounds involving the renormalized interaction energy,
together with local compactness of the vorticity measures. For stationary solutions
as $\varepsilon\to0$, their study \cite{SSlimit} uses the stress tensor
to derive equilibrium conditions for limiting vorticity measures under
energy and gradient bounds. Both works include magnetic and nonmagnetic
models. We also refer to the book \cite{SandierSerfaty} which develops related techniques for
the magnetic model.  

Existence of several solutions on bounded domains can also be studied
through topological and perturbative methods.
For smooth bounded domains and fixed Dirichlet data in the unit circle,
Lin \cite{Lin95} constructed solutions whose vortices converge, as
$\varepsilon\to0$, to nondegenerate local minima
of the renormalized energy. His proof uses the Ginzburg--Landau heat
flow to control the vortex positions before taking the limit as time
tends to infinity. Later,
Lin--Lin \cite{LinLin} extended this construction to
nondegenerate critical points of the renormalized energy, including
saddle points, for smooth boundary data of positive degree. It is worth pointing out that
Almeida--Bethuel \cite{AB} related the topology of energy sublevel sets
to the renormalized energy and proved multiplicity results. On the other hand,
Del Pino--Kowalczyk--Musso \cite{DKM} reduced constructions with prescribed
degrees $\pm1$ to critical points of a perturbation of the renormalized
energy, for both Dirichlet and Neumann conditions. We emphasize that the distinction between multiplicity and uniqueness also depends on
the boundary data in a delicate way. For a disk with boundary value $x/|x|$,
Chen--Liu--Wei--Yang \cite{CLWY} proved that the global minimizer is
the radial solution. Their proof uses the entire vortex of degree one,
angular Fourier decomposition, and an identity of Picone type, and answers Brezis' Open Problem~2.2.
Deng--Li--Liu--Ma \cite[Theorem~1.1]{dllm} proved uniqueness of
classical and $H^1$ weak solutions on the unit disk with identity
boundary data for every $\varepsilon>0$, answering Brezis' Open Problem~2.1.

For the standard Ginzburg--Landau energy on the unit ball
$B^N\subset\R^N$, with maps into $\R^N$ and boundary value $u=x$,
Ignat--Nguyen \cite[Theorem~1.2]{IN24} proved local minimality
and positivity of the second variation at the radial vortex for
every $N\ge2$ and $\varepsilon>0$.
Their subsequent work \cite[Theorem~1.1]{IN26} proves uniqueness
of the global minimizer under the same assumptions.

It is also worth mentioning that vortex interactions also govern the slow evolution of the
Ginzburg--Landau heat flow. In smooth bounded domains with fixed unit boundary data, Lin
\cite[Theorems~3.7 and~5.1]{Lin96} studied initial maps satisfying
energy and convergence assumptions. He showed that the limiting
vortex positions stay fixed on time scales $o(|\log\varepsilon|)$
and established continuity of the limiting vortex sets on the
logarithmic time scale.
On $\R^2$, Lin--Xin \cite{LinXin}
considered initial data formed by products of distant vortices of
degree $+1$ and perturbations decaying at infinity. As the separations
tend to infinity, they derived a limiting motion law governed by the
negative gradient of the logarithmic interaction energy on a slow time scale.

On the whole plane, there has been earlier work towards construction nonradial solutions. Ovchinnikov--Sigal \cite{OSenergy} studied the
interaction energy of vortices with prescribed zeros and local degrees.
They derived its expansion as the mutual distances tend to infinity.
Their subsequent study of solutions with broken rotational
symmetry \cite{OS} involves a correlation coefficient. However,
Esposito \cite{Esposito} and Kurzke \cite{Kurzke} proved cancellations
in this coefficient, which indicates that such a construction seems to involve complicated analysis of the higher order terms and would be very difficult.

Related examples arise for the complex sine--Gordon II equation.
Chen--Gui--Liu--Yang--Yang \cite[Theorem~1]{CGLYY} constructed an
explicit family of vortices of degree two, depending on one real
parameter and containing nonradial solutions. Their family bifurcates
from the standard vortex along additional bounded solutions of the
linearized equation. These formulas describe how rotational symmetry
can break while the degree stays fixed.
The construction concerns an integrable equation and illustrates the
role of bounded Jacobi fields in bifurcation from a rotationally
symmetric solution.

For the Ginzburg--Landau equation, recent work also relates the angular
Fourier decomposition to the Morse index and bounded Jacobi fields of
the standard vortices. The Morse index is the largest dimension of a
real subspace on which the second variation is negative definite. 
Del Pino--Liu--Musso--Wei--Yang \cite[Theorems~1.1--1.2]{DPLMWY}
computed the indices of the standard vortices of degrees two and three,
and proved that their bounded Jacobi fields are generated by the phase
and translation symmetries. Their argument uses angular Fourier
decomposition, algebraic bounds for $f_q$, and positive functions for
the resulting scalar differential operators.

The geometry of the domain provides a further parameter in vortex
problems. For tangent vector fields on closed oriented smooth surfaces,
Ignat--Jerrard \cite{IJ21} derived the renormalized interaction
energy as a second order $\Gamma$ limit.
Their formula describes its dependence on Gaussian curvature.
Vortices for tangent vector fields on a circular cone were recently
studied by Cofoid--Golovaty--Sandier--Sternberg \cite{CGSS}.
Their asymptotic analysis determines vortex degrees and locations,
including a vortex at the vertex. Note that for complex valued functions on a
cone, the degree is the winding number of the image of a loop.
Thus the angle of the domain can vary continuously while this integer
degree stays fixed. This will also be used in our later analysis.

We can now state our main result.
\begin{mainresult}\label{main}
There exist $q\in\{8,9\}$ and $u\in C^\infty(\R^2;\C)$ satisfying
\eqref{equation}, with $\deg(u,\infty)=q$, such that
\[
 u\ne\zeta V_q(\,\cdot-a\,)
 \quad\text{for every }a\in\R^2,\ \zeta\in\C,\ |\zeta|=1.
\]
Moreover, for $z\in\C$, $r>0$, and $0<\theta<\pi/q$, we have
\begin{gather*}
 u(e^{2\pi i/q}z)=u(z),\qquad
 u(\bar z)=\overline{u(z)},\qquad
 \operatorname{Im}u(re^{i\theta})>0,\\
 \int_{\R^2}(1-|u|^2)^2\,dx=2\pi q^2.
\end{gather*}
\end{mainresult}
We remark that the theorem only guarantees the existence for at least one of the degrees
$8,9$. Due to the topological argument we are unable to identify the precise degree.

The proof begins with a family of equations on cone. We vary the cone angle through
$ds_p^2=dr^2+p^{-2}r^2d\phi^2$, with $p\in[8,9]$.
The equation on this cone is
\[
 -\Delta_pU=U(1-|U|^2),\qquad
 \Delta_p=\partial_r^2+r^{-1}\partial_r+p^2r^{-2}\partial_\phi^2.
\]
We consider solutions of degree one with conjugate reflection,
positive imaginary part for $0<\phi<\pi$, and normalization
$U(r,\cdot)\to e^{i\cdot}$ in $C^1(S^1)$. In this context, the standard solution is $V_p=f_p(r)e^{i\phi}$.
Note that at an integer $p=q$, the formula $u(re^{i\theta})=U(r,q\theta)$
produces an entire solution of degree $q$.

The first step is compactness of this positive class through deformation of $p$.
Quantization and reflection show that any sequence with zeros tending
to infinity separates into $N$ entire vortices of degrees $\pm1$,
where the limiting parameter is the odd integer $N$.
An expansion of the potential energy shows that $p_n>N$ for all
sufficiently large $n$. Since $9$ is the only possible odd integer in
$[8,9]$, this strict inequality excludes zeros tending to infinity
along the chosen interval.

 We can formulate the equation on a fixed real
Banach space using weighted estimates. The angular Fourier decomposition and positive supersolutions used in Del Pino-Liu-Musso-Wei-Yang
\cite{DPLMWY} yield the index calculation after the estimates for
$p=8,9$ are established. At these two parameters the linearizations
are invertible, with Morse indices $3$ and $4$ in the reflection space.
Their local degrees are
therefore $-1$ and $1$. Homotopy invariance of the degree of Leray and
Schauder implies the existence of an additional positive solution at
$p=8$ or $p=9$. Its lift proves the theorem. Thus the construction
combines a strict estimate on separating vortices with a change of
index, and produces an entire solution outside the orbit of the standard
vortex.

The argument presented in this paper should be able to be used for showing that additional solutions would exist when we deform the higher degree. It could involve more calculations but the main process should be almost the same. The choice of $8$ and $9$ reflects the spectral information established here. Proposition ~\ref{posiwcomp} already proves compactness on every fixed interval $[2m,2m+1]$ with integer $m\ge2$. Extending the argument to another such pair would require the corresponding formulation in weighted spaces and further spectral analysis in the same reflection class.   

This paper is organized as follows.
Section~\ref{cones} develops the cone equation and its
compactness properties, and Section~\ref{separation} estimates
separating vortices.
Sections~\ref{weighted} and~\ref{indices}
establish the weighted formulation, compute the indices, and complete
the degree argument, with elementary inequalities collected in
Appendix~\ref{radialineq}.
\raggedbottom

\section{The Ginzburg–Landau equation on cones}\label{cones}
 
In this section, we introduce a continuous parameter by studying
the Ginzburg--Landau equation on cones of angle $2\pi/p$,
where $p>2$. When $p=q$ is an integer, a degree-one solution
on the cone lifts to a solution of degree $q$ on $\R^2$. We first define weak solutions at the
vertex and establish local estimates and potential-energy
quantization. We then introduce the reflection and positivity
conditions under which we will study compactness in the
following sections.

Throughout, we regard $\R^2$ as $\C$. For $p>2$, let $C_p$ be the metric completion of
$(0,\infty)\times S^1$ with vertex $o$ and metric
\begin{equation*}
 ds_p^2=dr^2+p^{-2}r^2d\phi^2,\qquad
 dA_p=p^{-1}r\,dr\,d\phi,\qquad
 \Delta_p=\partial_r^2+r^{-1}\partial_r+p^2r^{-2}\partial_\phi^2.
\end{equation*}
The metric is $ds_p^2$ on $C_p$ and is Euclidean in local isometric
coordinates or on the plane.
Here $r$ is the distance from the vertex and $\phi$ has period $2\pi$.
The circle $r=R$ has length $2\pi R/p$, so the total angle at the vertex
is $2\pi/p$. Equivalently, one obtains $C_p$ by gluing the two radial
sides of the Euclidean sector $0\le\theta\le2\pi/p$, with
$\theta=\phi/p$. The cone includes the vertex and extends to arbitrarily
large $r$. Its underlying topological space is the whole plane, as
expressed by the coordinate $y=re^{i\phi}$, with $y=0$ at $o$.

We study the equation
\[
 -\Delta_pU=U(1-|U|^2)\qquad\text{on}\quad  C_p.
\]
Away from $o$, we choose a local branch of $\phi$ and setting
$z=re^{i\phi/p}$ transforms $\Delta_p$ into the Euclidean
Laplacian. Hence $U$, expressed in these coordinates, satisfies
the standard Ginzburg--Landau equation, and the usual interior
elliptic estimates apply. Euclidean disks below are understood
to lie in such coordinate neighborhoods. We use $y=re^{i\phi}$
as a global coordinate to identify the underlying domain with
$\R^2$.

To formulate the equation across the vertex, we use the global conformal
coordinate $w=r^pe^{i\phi}$. It represents $C_p$ by the whole complex plane:
\begin{gather*}
 ds_p^2=a_p(w)|dw|^2,\qquad a_p(w)=p^{-2}|w|^{\frac{2}{p}-2},\\
 \int_{C_p}|\nabla_pU|^2\,dA_p=\int_\C|\nabla g|^2\,dw,
 \qquad g(w)=U(|w|^{\frac{1}{p}},\arg w).
\end{gather*}
Here $\nabla_p$ denotes the metric gradient or covariant derivative,
while $\nabla$ denotes the Euclidean derivative and $\nabla_y$ the
Euclidean derivative in the coordinate $y$. Superscripts indicate
repeated spatial differentiation. The symbol $D$ denotes the real
differential of a map, with a subscript indicating the differentiated
variable when needed. The measure $dw$ is Lebesgue measure on the
conformal plane. For nonempty
$K\Subset(2,\infty)$, write 
$$p_- =\min K\quad\mbox{and}\quad p_+=\max K.$$
Constants may depend on a fixed compact interval
$K\Subset(2,\infty)$.
A bounded map $U$ is a \emph{weak solution} if its conformal
representative $g\in H^1_{\rm loc}(\C)\cap L^\infty(\C)$ satisfies
\begin{equation*}
 -\Delta_wg=a_p(w)F_0(g)\quad\text{in }\mathcal D'(\C),
 \qquad F_0(z)=z(1-|z|^2).
\end{equation*}
In particular, the weighted equation holds on all of $\C$, including
$w=0$. More explicitly, for every $\eta\in C_c^\infty(\C;\C)$,
\[
 \operatorname{Re}\int_\C\nabla g\cdot\nabla\bar\eta\,dw
 =\operatorname{Re}\int_\C a_p(w)F_0(g)\bar\eta\,dw.
\]
The test functions may be nonzero at the origin. This formulation
excludes a Dirac term there and imposes the equation at the vertex
through the weak identity.
Since
$a_p\in L^s_{\rm loc}(\C)$ for $1<s<\frac{p}{p-1}$, local elliptic
regularity implies $g\in W^{2,s}_{\rm loc}(\C)$ for every such $s$.
By Sobolev embedding, we then have $g\in C^{0,\alpha}_{\rm loc}(\C)$ for every
$0<\alpha<\frac{2}{p}$. Moreover,
$g\in C^\infty(\C\setminus\{0\})$. We use the continuous representative
and set $U(o)=g(0)$. The definition allows $U(o)\ne0$ and imposes
the equation without prescribing a boundary value at the vertex.
For example, the constant map $U\equiv1$ is a weak solution.
We use the local estimates in
\cite[Lecture~4, item~E]{MironescuPDE} and
\cite[Theorems~1.4--1.5, Corollaries~2.16--2.17]{FRO} throughout.
The potential energy is
\begin{equation*}
 \mathcal P_p(U)=\int_{C_p}(1-|U|^2)^2\,dA_p
       =\int_\C a_p(w)(1-|g|^2)^2\,dw.
\end{equation*}
Thus $\mathcal P$ and $\mathcal P_p$ denote scalar potential integrals
on the Euclidean plane and on $C_p$, respectively.

In the region where $r>0$ and $U\ne0$, we express $U$ locally as
$U=\rho e^{i\vartheta}$, with $\rho=|U|>0$ and $\vartheta$ real.
The function $\vartheta$ is an argument of the complex value $U$,
also called its phase. The coordinate $\phi$ describes position on
the cone, and $\vartheta(r,\phi)$ describes the direction of
$U(r,\phi)$ in the complex plane. On each connected overlap, two
choices of $\vartheta$ differ by a constant integer multiple of
$2\pi$, so their gradients agree. We denote this vector field by $\mathbf{P}$:
\begin{equation}\label{arggrad}
 \mathbf{P}=\nabla_p\vartheta
   =\frac{\operatorname{Im}(\overline U \nabla_pU)}{|U|^2},
 \qquad j(U)=\operatorname{Im}(\overline U \nabla_pU)=\rho^2\mathbf{P}.
\end{equation}
Thus $\nabla_p\arg U$ always means the field in \eqref{arggrad}, computed
using any local argument. The expression for $j(U)$ also applies at
zeros of $U$ away from the vertex. In the orthonormal frame
$e_r=\partial_r$ and $e_\theta=(p/r)\partial_\phi$, where
$\theta=\phi/p$, we have
\[
 \mathbf{P}=(\partial_r\vartheta)e_r
       +\frac p r(\partial_\phi\vartheta)e_\theta,
 \qquad |\mathbf{P}|^2=(\partial_r\vartheta)^2
       +\frac{p^2}{r^2}(\partial_\phi\vartheta)^2.
\]
These formulas measure the spatial variation of the argument using
the cone metric. The identity
$|\nabla_pU|^2=|\nabla_p\rho|^2+\rho^2|\mathbf{P}|^2$ explains the weight $\rho^2$
in the integrals of $\mathbf{P}$ below. For example, $U=f_p(r)e^{i\phi}$ satisfies
$\mathbf{P}=(p/r)e_\theta$ wherever $f_p>0$.

The real parameter $p$ determines the cone angle. The degree is instead
an integer associated with the values of $U$ along a large circle.
If $|U(r,\phi)|\to1$ uniformly as $r\to\infty$, then $U$ is nonzero
outside a sufficiently large ball, and
\[
 \deg(U,\infty)
 =\frac{1}{2\pi}\int_0^{2\pi}
   \frac{\operatorname{Im}\bigl(\overline{U(R,\phi)}
                  \partial_\phi U(R,\phi)\bigr)}{|U(R,\phi)|^2}\,d\phi
\]
is independent of the sufficiently large radius $R$. It is the winding
number of $U/|U|$ as the domain circle is traversed once. In particular,
$f_p(r)e^{i\phi}$ has degree one on $C_p$ for every $p>2$.
The factor $p$ in the metric changes lengths and energies, while this
winding number counts turns in the image circle.
Several zeros may lie inside the circle. When they are isolated, the
degree equals the sum of their local degrees. Thus global degree one
alone leaves the number and location of the zeros undetermined.
The restrictions on local degrees used later follow from additional
properties of the solutions.

At an integer parameter $p=q$, the map
$u(re^{i\theta})=U(r,q\theta)$ is single valued on the whole Euclidean
plane. As $\theta$ runs from $0$ to $2\pi$, the cone angle variable
$\phi=q\theta$ traverses its circle $q$ times. Consequently,
$\deg(u,\infty)=q\deg(U,\infty)$. For noninteger $p$, $2\pi$
periodicity of $\theta\to U(r,p\theta)$ would require an additional
symmetry of $U$. The following lemma proves smoothness of the integer
lift, including at the origin.

\begin{lemma}\label{intelift}
Let $q\ge3$ be an integer, and let $U$ be a bounded weak solution
of $-\Delta_qU=F_0(U)$ on $C_q$. Assume that
$|U(r,\phi)|\to1$ uniformly in $\phi$ as $r\to\infty$ and that
$\deg(U,\infty)=1$. Then the map
\begin{equation*}
 u(re^{i\theta})=U(r,q\theta)=g((re^{i\theta})^q),\qquad u(0)=g(0),
\end{equation*}
belongs to $C^\infty(\R^2)$, solves $-\Delta u=F_0(u)$ in $\R^2$, and
satisfies
\[
 \deg(u,\infty)=q,\qquad \lim_{|x|\to\infty}|u(x)|=1,
 \qquad u(e^{\frac{2\pi i}{q}}z)=u(z),\qquad z\in\C.
\]
Moreover,
\begin{equation}\label{potelift}
 \int_{\R^2}(1-|u|^2)^2\,dx=q\mathcal P_q(U).
\end{equation}
\end{lemma}
\begin{proof}
Away from zero, we have
\begin{gather*}
 -\Delta u=q^2|z|^{2q-2}a_q(z^q)F_0(g(z^q))=F_0(u),\\
 \int_{B_R(0)\setminus B_\delta(0)}|\nabla u|^2\,dx
 =q\int_{B_{R^q}(0)\setminus B_{\delta^q}(0)}|\nabla g|^2\,dw.
\end{gather*}
Thus $u$ is bounded and its gradient belongs locally to $L^2$. Logarithmic
cutoffs with $\|\nabla\eta_\varepsilon\|_{L^2(\R^2)}\to0$ extend both
$u\in H^1_{\rm loc}(\R^2)$ and its weak equation at zero. Smoothness follows by elliptic
bootstrapping \cite[Corollaries~2.16--2.17 and Proposition~2.18]{FRO}. The covering has degree $q$, preserves
the uniform modulus limit, and yields the stated symmetry and
\eqref{potelift} by change of variables.
\end{proof}

The radial solution satisfies estimates uniform on compact parameter intervals.

\begin{proposition}\label{radial}
For every $p>2$ there is a unique positive function
$f_p\in C^\infty((0,\infty))$ satisfying
\begin{equation}\label{radialode}
 f_p''+r^{-1}f_p'-p^2r^{-2}f_p+f_p(1-f_p^2)=0,
 \quad \lim_{r\to\infty}f_p(r)=1,
 \quad\int_0^R(r|f_p'|^2+r^{-1}f_p^2)\,dr<\infty
 \quad(R>0).
\end{equation}
It satisfies $0<f_p(r)<1$ and $f_p'(r)>0$ for every $r>0$.
For every nonempty compact set $K\Subset(2,\infty)$,
\begin{equation}\label{radiorig}
 f_p(r)=\kappa_pr^p\left(1-\frac{r^2}{4(p+1)}+O_K(r^4)\right)\quad\mbox{as}\quad r\to0,\qquad
 0<c_K\le\kappa_p\le C_K
 ,\quad \forall p\in K,
\end{equation}
with a remainder uniform in $p$. For every $0<\alpha<1$,
\begin{equation}\label{radiinfi}
 \left\|f_p-1+\frac{p^2}{2r^2}\right\|_{C^{2,\alpha}((r_*-1,r_*+1))}
 \le C_{K,\alpha}r_*^{-4}\qquad(r_*\ge R_K).
\end{equation}
In particular, $f_p'=O_K(r^{-3})$ and $f_p''=O_K(r^{-4})$ as
$r\to\infty$. For every integer $k\ge0$ and $0<a<b<\infty$, the map
$p\to f_p$ is continuous into $C^k([a,b])$.
\end{proposition}
\begin{proof}
For any degree $p$ solution, Herv\'e--Herv\'e \cite[Theorem, p.~428]{HH} prove
existence, uniqueness, strict monotonicity and the expansion at infinity.
Pacard--Rivi\`ere \cite[Theorem~3.1]{PR} state the result for $p=1$.
For real $p>2$,  let us replace the power series at zero in \cite[\S2]{HH}
by $f=\kappa r^pv$, where
\[
 v(r)=1-\int_0^rt^{-2p-1}\int_0^ts^{2p+1}
                 (1-\kappa^2s^{2p}v(s)^2)v(s)\,ds\,dt.
\]
By local contraction, we obtain the regular solution and its continuous dependence
on $(p,\kappa)$. The natural energy condition excludes the solution
$r^{-p}$ of the linear equation at zero.
The shooting and Wronskian arguments in \cite[\S3]{HH} then apply:
they use $p^2>0$, the powers $r^{\pm p}$, and the regular Bessel solution
of real order $p$, whose positive zeros exclude small positive shooting
parameters. They prove the stated qualitative conclusions for real $p$.

We verify uniformity for $p\in K$. From the integral identity, we have
$1/2\le v\le1$ on $[0,1]$, hence $\kappa_p\le2$ and
\eqref{radiorig}, apart from its lower bound.
For $A=2p_+(p_++1)$, direct differentiation shows that
$\ell_p=r^p(A+r^2)^{-p/2}$ is a subsolution of \eqref{radialode}.
The quotient identity in \cite[(8)]{HH} implies
$\ell_p\le f_p$, so $\kappa_p\ge A^{-p_+/2}$ and
$h=1-f_p=O_K(r^{-2})$.
For $L=-\partial_r^2-r^{-1}\partial_r+2$, the equation becomes
\[
 L\left(h-\frac{p^2}{2r^2}\right)
 =\frac{2p^2}{r^4}-\frac{p^2h}{r^2}+3h^2-h^3.
\]
Local ODE estimates bound the right side in $C^{0,\alpha}$ by
$C_{K,\alpha}r^{-4}$. Since $L(r^{-4})\ge r^{-4}$ for $r\ge4$,
the maximum principle and interior estimates prove
\eqref{radiinfi} and the derivative bounds.
For $p_n\to p$, these uniform bounds imply compactness at zero and infinity.
Every locally convergent subsequence has limit $f_p$ by uniqueness,
which proves the stated continuity.
\end{proof}

The standard solution $V_p=f_p(r)e^{i\phi}$ is natural: in conformal
coordinates $g_p(w)=\kappa_pw(1+O_K(|w|^{\frac2p}))$, so its weak equation
holds at zero. The estimate at infinity implies $\mathcal P_p(V_p)<\infty$.
It has conjugate reflection and positive imaginary part for $0<\phi<\pi$.

In our case, a bounded solution with $|U|\to1$ still satisfies the well known estimate $|U|\le1$.
Indeed, we have
$$\Delta_w|g|^2=2|\nabla g|^2+2a_p|g|^2(|g|^2-1).$$
Testing with $\bigl(|g|^2-(1+\varepsilon)^2\bigr)_+$, which has compact
support, and letting $\varepsilon\to0$ proves the assertion on the whole
conformal plane.

Uniform regularity at the vertex also follows from the weak formulation.

\begin{lemma}\label{tipesti}
Let $K\Subset(2,\infty)$ be nonempty. For any $p\in K$, suppose that $U$ is a 
weak solution of $-\Delta_pU=F_0(U)$ on $C_p$ with $|U|\le1$.
Set $c_0=U(o)$ and $c_2=-\frac14F_0(c_0)$. There is $C_K>0$,
independent of $p$ and $U$, such that
\begin{equation*}
 |U-c_0|\le C_Kr^2,\quad |\nabla_pU|\le C_Kr,\quad |\nabla_p^2U|\le C_K
 \quad\mbox{for}\quad 0<r<\frac14,\qquad
 \|\nabla_pU\|_{L^\infty(C_p)}\le C_K.
\end{equation*}
For $0<\sigma<\min\{p_--2,2\}$, $0<\alpha<1$, and
$A=\{y\in\R^2:1<|y|<2\}$ in the fixed coordinate $y=re^{i\phi}$,
\begin{equation}\label{tipstrong}
 \|(U-c_0-c_2r^2)(s\,\cdot)\|_{C^{2,\alpha}(A)}
 \le C_{K,\sigma,\alpha}s^{2+\sigma}
 \qquad\left(0<s<\frac1{16}\right).
\end{equation}
If $U(r,-\phi)=\overline{U(r,\phi)}$, then $c_0,c_2\in\R$.
\end{lemma}
\begin{proof}
In the fixed coordinates, one has
\[
 \Delta_p=\operatorname{div}(A_p\nabla),\qquad
 A_p=p^2I+(1-p^2)\frac{y\otimes y}{|y|^2},\qquad I\le A_p\le p_+^2I.
\]
For a solution with zero boundary values $-\Delta_pv=h$ on $B_a(o)$ with
$|h|\le Mr^m$, the radial barrier implies
\[
 |v(r,\phi)|\le\frac{M}{(m+2)^2}(a^{m+2}-r^{m+2})\qquad(m=0,2).
\]
The radial barrier above satisfies its equation weakly on the
whole disk, including $o$, so the usual weak maximum principle
applies. Let $H$ satisfy $\Delta_pH=0$ in $B_a(o)$ with
$H=U$ on its boundary. Since $|F_0(U)|\le1$, the preceding
estimate implies $\|U-H\|_{L^\infty(B_a(o))}\le a^2/4$.
In the coordinate $w=r^pe^{i\phi}$, $H$ is harmonic on the
disk of radius $a^p$. By the Poisson formula
\cite[\S1.2, (1.15)]{FRO}, we therefore obtain, with
$\omega(a)=\operatorname{diam}U(B_a(o))$,
\[
 \omega(\lambda a)
 \le \operatorname{diam}H(B_{\lambda a}(o))+\frac{a^2}{2}
 \le 8\lambda^{p_-}\omega(a)+\frac{a^2}{2}.
\]
Since $p_->2$, we can choose $0<\lambda<1/2$ such that
$8\lambda^{p_-}\le\lambda^2/2$. By iterating the preceding
inequality, we deduce $\omega(r)\le C_Kr^2$. Applying interior
estimates \cite[Corollaries~2.16--2.17 and Proposition~2.18]{FRO}
to $U(s\,\cdot)-c_0$ on the fixed annulus
$\{y\in\R^2:1/2<|y|<4\}$, we obtain
\[
 \|U(s\,\cdot)-c_0\|_{C^{2,\alpha}(A)}
 \le C_{K,\alpha}s^2.
\]
Together with interior estimates on fixed-radius Euclidean
disks away from $o$, this also proves the global gradient bound. To prove \eqref{tipstrong}, set $\eta=U-c_0-c_2r^2$.
By the choice of $c_2$ and the estimate just proved, we have
\[
 -\Delta_p\eta=F_0(U)-F_0(c_0)=O_K(r^2),\qquad \eta(o)=0.
\]
We repeat the harmonic comparison argument, using the barrier
with $m=2$. For $b(t)=\sup_{B_t(o)}|\eta|$, we thus obtain
\[
 b(\mu t)\le4\mu^{p_-}b(t)+C_Kt^4.
\]
Let $\beta=2+\sigma<\min\{p_-,4\}$ and choose $\mu>0$
small enough that $4\mu^{p_-}\le\mu^\beta/2$.
We then deduce by iteration that $b(t)\le C_{K,\beta}t^\beta$.
After scaling by $s$, the right-hand side of the equation
for $\eta$ has $C^{0,\alpha}(A)$ norm bounded by $C_{K,\alpha}s^4$.
Schauder estimates therefore yield \eqref{tipstrong}.
Finally, reflection symmetry implies $c_0,c_2\in\R$.
\end{proof}

These estimates make local compactness uniform as the parameter varies.

\begin{lemma}\label{locacomp}
Let $K\Subset(2,\infty)$ be nonempty, let $p_n\in K$ satisfy
$p_n\to p_*$, and let $U_n$ be weak solutions of
$-\Delta_{p_n}U_n=F_0(U_n)$ on $C_{p_n}$ with $|U_n|\le1$.
In the fixed coordinates, a subsequence converges to a   weak
solution $U_*$ on $C_{p_*}$, with $|U_*|\le1$, such that for every $R>0$,
\[
 \|U_n-U_*\|_{C^0(\overline{B_R(o)})}
 +\|\nabla_{p_*}(U_n-U_*)\|_{L^2(B_R(o))}\to0.
\]
Here the balls and the $L^2(B_R(o))$ norm use the metric of $C_{p_*}$.
The convergence is smooth on every compact subset of $C_{p_*}\setminus\{o\}$.
The conformal representatives satisfy
$g_n\to g_*$ in $C^0_{\rm loc}(\C)\cap H^1_{\rm loc}(\C)$.
\end{lemma}
\begin{proof}
For $1<s<p_+/(p_+-1)$, dominated convergence implies
$a_{p_n}\to a_{p_*}$ in $L^s(B_R(0))$ for every $R>0$.
By local $W^{2,s}$ estimates \cite[Theorems~9.9 and~9.11]{GT}
and compact embedding, we obtain, after passing to a subsequence,
$g_n\to g_*$ locally uniformly. Thus
$a_{p_n}F_0(g_n)\to a_{p_*}F_0(g_*)$ in $L^s_{\rm loc}$.
The limit solves the equation. Applying interior estimates
to the difference, we obtain convergence in $W^{2,s}_{\rm loc}$,
hence in $H^1_{\rm loc}$. Away from zero, smooth convergence follows from elliptic regularity.

Near the vertex, Lemma~\ref{tipesti} yields
\[
 \sup_{B_\varepsilon(o)}|U_n-U_*|
 \le |U_n(o)-U_*(o)|+C_K\varepsilon^2,\qquad
 \int_{B_\varepsilon(o)}
 |\nabla_{p_*}(U_n-U_*)|^2\,dA_{p_*}\le C_K\varepsilon^4.
\]
Since $U_n(o)\to U_*(o)$, these estimates and convergence
on annuli imply local uniform and $H^1_{\rm loc}$ convergence
in the fixed $y$ coordinate.
\end{proof}

The fact that finite potential energy is determined by the degree at infinity can also be generalized in this context. 

\begin{proposition}\label{quantize}
Let $p>2$ and let $U$ be a   weak solution of
$-\Delta_pU=F_0(U)$ on $C_p$ satisfying $|U|\le1$ and $\mathcal P_p(U)<\infty$.
Then $|U(r,\phi)|\to1$ uniformly in $\phi$ as $r\to\infty$,
the degree $c=\deg(U,\infty)\in\mathbb Z$ is defined, and
\begin{equation*}
 \mathcal P_p(U)=2\pi pc^2.
\end{equation*}
There are $R>0$ and a smooth real function $\chi$ on
$E_R=C_p\setminus\overline{B_R(o)}$ such that
\begin{equation}\label{coneprema}
 U=\rho e^{i(c\phi+\chi)},\quad \rho=|U|\ge\frac12,\quad
 \mathbf{P}=\nabla_p\arg U=\frac{pc}{r}e_\theta+\nabla_p\chi,
 \quad \nabla_p\rho,\nabla_p\chi\in L^2(E_R),
\end{equation}
where $\theta=\frac\phi p$. If $c=0$, then $U$ is constant and $|U|=1$.
\end{proposition}
\begin{proof}
We adapt Brezis--Merle--Rivi\`ere
\cite[\S2, Proposition~1 and the proof of Theorem~1]{BMR}.
Their argument on disjoint Euclidean disks implies $|U|\to1$ uniformly.
Translated limits have modulus one and are harmonic, so $\nabla_pU\to0$.
For $w_0=1-|U|^2$ and $B=|\nabla_pU|^2-w_0$, we compute directly that
\begin{equation*}
 (\Delta_p-2|U|^2)B=2|\nabla_p^2U|^2+|\nabla_pw_0|^2\ge0.
\end{equation*}
Lemma~\ref{tipesti} implies $B(o)\le0$. Since $B\to0$ at infinity,
the maximum principle yields
\begin{equation}\label{bernboun}
 |\nabla_pU|^2\le w_0,\qquad \int_{C_p}|\nabla_pU|^4\,dA_p\le \mathcal P_p(U).
\end{equation}
Choose $R_0$ with $\rho=|U|\ge1/2$ on $E_{R_0}$ and write
$\mathbf{P}=pc\,r^{-1}e_\theta+\nabla_p\chi$ there. The modulus equation is
\begin{equation}\label{modueequa}
 (-\Delta_p+\rho(1+\rho))(1-\rho)=\rho|\mathbf{P}|^2.
\end{equation}
Since $1-\rho\le w_0$, $|\mathbf{P}|^2\le4w_0$ and $\rho(1+\rho)\ge3/4$,
testing with truncated copies of $1-\rho$ proves $\nabla_p\rho\in L^2(E_{R_0})$.

Set $H(R)=\int_{B_R(o)\setminus B_{R_0}(o)}|\nabla_p\chi|^2\,dA_p$.
Equation \eqref{bernboun} implies $H(R)\le CR$.
The divergence of $\operatorname{Im}(\overline U \nabla_pU)$ is zero, and
its normal integral at $r=\varepsilon$ is $O_p(\varepsilon^2)$ by
Lemma~\ref{tipesti}. Thus
$\int_{\partial B_R(o)}\rho^2\partial_r\chi\,dS_p=0$.
Testing the phase equation with
$\chi-\langle\chi\rangle_{\partial B_R(o)}$, as in
\cite[proof of Proposition~1, Step~2]{BMR}, now yields
\[
 \frac14H(R)\le\frac{R}{2p}H'(R)+C+CH(R)^{1/2}.
\]
Here the circular Poincar\'e constant is $R/p$, the inner boundary
term is fixed, and the unweighted angular term integrates to zero.
Choose $1<\beta<p/2$. After absorption, we obtain $H(R)\le RH'(R)/\beta+C$.
An unbounded $H$ would then satisfy $H(R)-C\ge cR^\beta$ for large $R$,
contrary to its linear bound. Hence $\nabla_p\chi\in L^2(E_{R_0})$.

The Pohozaev calculation with $r\partial_r$ has inner boundary contribution
$O_p(\varepsilon^2)$ by Lemma~\ref{tipesti}, and yields
\begin{equation*}
 \int_{B_L(o)}w_0^2\,dA_p
 =L\int_{\partial B_L(o)}
 \left(L^{-2}|\partial_\theta U|^2-|\partial_rU|^2+\tfrac12w_0^2\right)\,dS_p.
\end{equation*}
Coarea provides $L_\ell\to\infty$ with
$L_\ell\int_{\partial B_{L_\ell}(o)}
(|\nabla_p\rho|^2+|\nabla_p\chi|^2+w_0^2)\,dS_p\to0$.
Substitution of \eqref{coneprema}, using
$|\partial B_L(o)|=2\pi L/p$ and $\rho\to1$, proves
$\mathcal P_p(U)=2\pi pc^2$. For $c=0$, this identity and
\eqref{bernboun} imply $|U|=1$ and $\nabla_pU=0$.
\end{proof}
Since $V_p=f_p(r)e^{i\phi}$ has degree one on $C_p$,
the preceding proposition implies $\mathcal P_p(V_p)=2\pi p.$ Next we prove a small-energy estimate, uniform for $p\in K$.

For a solution $U$ on $C_p$, set
\[
 e(U)=\frac12|\nabla_pU|^2+\frac14(1-|U|^2)^2.
\]

\begin{lemma}\label{unifhein}
Let $K\Subset(2,\infty)$ be nonempty. There are $\epsilon_K,C_K>0$
with the following property. If $p\in K$, $U$ is a   weak solution
of $-\Delta_pU=F_0(U)$ on $C_p$, $|U|\le1$, $x_0\in C_p$, $R>0$, and
\[
 |U|\ge\frac12\quad\text{on }B_R(x_0),\qquad
 E=\int_{B_R(x_0)}e(U)\,dA_p\le\epsilon_K.
\]
then
\begin{equation*}
 \sup_{B_{\frac R2}(x_0)}e(U)\le C_KR^{-2}E.
\end{equation*}
\end{lemma}
\begin{proof}
Since $|U|^2\ge\frac14$ and $|\nabla_pU|^2\le2e$, we compute by differentiation that
\begin{align*}
 \Delta_pe=|\nabla_p^2U|^2
 +|U|^2\left(1-|U|^2-\frac{|\nabla_pU|^2}{|U|^2}\right)^2+\bigl|\nabla_p(|U|^2)\bigr|^2
 -\frac{|\nabla_pU|^4}{|U|^2}\ge-16e^2.
\end{align*}
From Lemma~\ref{tipesti}, we have $|e-e(o)|\le C_Kr^2$,
$|\nabla_pe|\le C_Kr$, so the inner normal boundary integral is $O_K(\varepsilon^2)$ and the
inequality holds weakly across the vertex. Under
$\Phi_p(r,\phi)=re^{i\phi}$ it becomes
$-\operatorname{div}(A_p\nabla\widetilde e)\le16\widetilde e^2$,
with $I\le A_p\le p_+^2I$, $dA_p=p^{-1}dy$, and
\begin{equation}\label{conedcomp}
 p^{-1}|\Phi_p(x)-\Phi_p(x')|\le d_p(x,x')
 \le|\Phi_p(x)-\Phi_p(x')|.
\end{equation}
For $-\operatorname{div}(A\nabla v)\le bv$ on $B_a(y_0)$,
$v\ge0$, $I\le A\le p_+^2I$ and $ba^2\le1$,
we obtain from \cite[Theorem~8.17]{GT}, with interpolation and iteration,
\[
 \sup_{B_{\frac a2}(y_0)}v\le C_Ka^{-2}\int_{B_a(y_0)}v(y)\,dy.
\]
For the quadratic inequality on $B_1(0)$, maximize
$(\rho-|y|)^2v(y)$ at $y_*$, where $0<\rho<1$, and set
$M=v(y_*)$, $h=\frac{\rho-|y_*|}{2}$. Then $v\le4M$ on $B_h(y_*)$.
If $0<a=(64M)^{-\frac12}\le h$, the linear estimate implies
$M\le64C_KM\int_{B_1(0)}v(y)\,dy$, impossible for
$0<\int_{B_1(0)}v(y)\,dy<\frac1{128C_K}$. Otherwise
\[
 \rho^2v(0)\le4h^2M\le4C_K\int_{B_1(0)}v(y)\,dy.
\]
The case $M=0$ is immediate. Let $\rho\to1$ and apply this estimate to
$v(z)=\frac{R^2}{4}\widetilde e(\Phi_p(x)+\frac R2z)$,
$x\in B_{\frac R2}(x_0)$. Equation \eqref{conedcomp} implies
$B_{\frac R2}(\Phi_p(x))\subset\Phi_p(B_R(x_0))$ and
\[
 \int_{B_1(0)}v(z)\,dz\le p_+E,\qquad
 e(U)(x)=\frac4{R^2}v(0)\le C_KR^{-2}p_+E.
\]
Choose $\epsilon_K<\frac1{128C_Kp_+}$ and enlarge $C_K$.
\end{proof}

Using the estimates above, we now study limits obtained by
translating solutions around points tending to infinity.
These limits solve the Ginzburg--Landau equation on $\R^2$.
To relate their potential energy to their degree, we recall
the following result of Brezis--Merle--Rivi\`ere
\cite[pp.~35--38, Theorem~1, Remarks~1.1--1.2, and Proposition~1]{BMR}.
\begin{lemma}\label{planfpote}
Let $V\in C^2(\R^2)$ satisfy
\[
 -\Delta V=V(1-|V|^2)\quad\text{in }\R^2,\qquad
 |V|\le1,\qquad \int_{\R^2}(1-|V|^2)^2\,dx<\infty.
\]
Then $|V(x)|\to1$ uniformly as $|x|\to\infty$, and
\[
 d=\deg(V,\infty)\in\mathbb Z,\qquad
 \int_{\R^2}(1-|V|^2)^2\,dx=2\pi d^2.
\]
For some $R_V>0$, on $E_V=\R^2\setminus\overline{B_{R_V}(0)}$ one has
$V=\rho e^{i(d\vartheta+\chi)}$, where $\rho\ge\frac12$, $\chi$ is
a globally defined real function, and $\nabla\rho,\nabla\chi\in L^2(E_V)$.
If $d=0$, then $V$ is a constant of modulus one.
\end{lemma}

The following proposition expresses the potential integral in terms of
the central solution and the limits obtained after translation.
For a map $V$, we set $\mu_V=(1-|V|^2)^2$.
\begin{proposition}\label{poteddeco}
Let $K\Subset(2,\infty)$ be nonempty, let $p_n\in K$ with $p_n\to p_*$,
and let $U_n$ be   bounded solutions on $C_{p_n}$ such that
\[
 \lim_{r\to\infty}|U_n(r,\phi)|=1\quad\text{uniformly in }\phi,
 \qquad \deg(U_n,\infty)=1,\qquad \mathcal P_{p_n}(U_n)<\infty.
\]
After taking a subsequence, there exist an integer $J\ge0$, a  
bounded solution $U_0$ on $C_{p_*}$, and nonconstant solutions
$V_1,\ldots,V_J\in C^\infty(\R^2)$ of $-\Delta V_j=V_j(1-|V_j|^2)$,
all with finite potential integral, such that
\[
 U_n\to U_0\quad\hbox{locally as in Lemma~\ref{locacomp}},
 \qquad U_n\circ\iota_{j,n}\to V_j\quad\hbox{in }C^\infty_{\rm loc}(\R^2).
\]
Here $\iota_{j,n}$ are isometries preserving orientation centered at $a_{j,n}$,
defined on the Euclidean disks $B_{\kappa_Kr(a_{j,n})}(0)$ for some
$\kappa_K>0$ depending only on $K$, with
$r(a_{j,n})\to\infty$ and $d_{p_n}(a_{j,n},a_{k,n})\to\infty$ for $j\ne k$.
Writing $c=\deg(U_0,\infty)$ and $d_j=\deg(V_j,\infty)\ne0$, one has
\begin{equation}\label{degrprela}
 1=c+\sum_{j=1}^Jd_j,\qquad p_*=p_*c^2+\sum_{j=1}^Jd_j^2.
\end{equation}
If $J\ge1$ and
\[
 L_{*,n}=\min\left\{\min_{1\le j\le J}r(a_{j,n}),
             \min_{j\ne k}d_{p_n}(a_{j,n},a_{k,n})\right\},
 \quad \ell_n\to\infty,\quad \frac{\ell_n}{L_{*,n}}\to0,
\]
with the pairwise minimum equal to $+\infty$ when $J=1$,
let $\mu_{j,n}$ be the pushforward of
$\mu_{V_j}\mathbf1_{B_{\ell_n}(0)}$ under $\iota_{j,n}$, extended by zero,
and set $\mu_{0,n}(r,\phi)=\mu_{U_0}(r,\phi)$. Then
\begin{equation}\label{stroddeco}
 \left\|\mu_{U_n}-\mu_{0,n}-\sum_{j=1}^J\mu_{j,n}\right\|_{L^1(C_{p_n})}\to0.
\end{equation}
For $J=0$ the sum is absent. If $J\ge1$, then $c=0$, $U_0$ is a unit
constant, and $p_*\in\{3,5,7,\ldots\}$. If $J=0$, then $c=1$, and in the fixed
coordinates $\mu_{U_n}\to\mu_{U_0}$ in $L^1(C_{p_*})$ and
$1-|U_n|^2\to1-|U_0|^2$ in $L^2(C_{p_*})$.
\end{proposition}
\begin{proof}
Proposition~\ref{quantize} and
Lemma~\ref{locacomp} yield a central limit $U_0$ with
\[
 |U_n|\le1,\qquad \mathcal P_{p_n}(U_n)=2\pi p_n,\qquad
 \mathcal P_{p_*}(U_0)=2\pi p_*c^2\le2\pi p_*.
\]
We include $a_{0,n}=o$ among the centers. If $1-|U_n|^2$ has a fixed
positive lower bound at points tending away from all previous centers,
we move these points to the origin in Euclidean coordinates preserving
orientation. Interior compactness and Lemma~\ref{planfpote}
yield a nonconstant limit on $\R^2$ with potential integral
$2\pi d^2\ge2\pi$. Fatou's lemma on disjoint balls of fixed radius
implies $p_*c^2+\sum_jd_j^2\le p_*$.
Thus only finitely many centers can be selected, and
\begin{equation*}
 \lim_{R\to\infty}\limsup_{n\to\infty}
 \sup_{\substack{x\in C_{p_n}\\
 d_{p_n}(x,\{o,a_{1,n},\ldots,a_{J,n}\})\ge R}}(1-|U_n(x)|^2)=0.
\end{equation*}
For large fixed $R$ and then large $n$, set
$\Omega_{n,R}=C_{p_n}\setminus\bigcup_{j=0}^J\overline{B_R(a_{j,n})}$.
The disjoint holes have the limiting degrees and $\rho_n=|U_n|\ge\frac12$
on $\Omega_{n,R}$. By Stokes' theorem for $\mathbf{P}_n=\nabla_{p_n}\arg U_n$, we have
$1=c+\sum_jd_j$.

To estimate the remaining potential energy, we compare the phase
gradient of $U_n$ with the explicit vector field
$$\mathbf{P}_{0,n}=\nabla_{p_n}\left(c\arg w+\sum_jd_j\arg(w-w_{j,n})\right),
 \qquad w=r^{p_n}e^{i\phi},\quad w_{j,n}=w(a_{j,n}).$$
For arbitrary $p\in K$, use $w=r^pe^{i\phi}$ and take distinct
$x,a\in C_p\setminus\{o\}$. Put $r=r(x)$, $s=r(a)$,
$M=\max(r,s)$, $m=\min(r,s)$ and $\ell=d_p(x,a)$.
For the shortest developed angular separation $\alpha\in[0,\frac\pi p]$, we have
\[
 \ell^2=(r-s)^2+4rs\sin^2\frac\alpha2,\qquad
 |w(x)-w(a)|^2=(r^p-s^p)^2+4r^ps^p\sin^2\frac{p\alpha}2.
\]
When $m\ge\frac M2$, the mean value theorem and
$\sin\frac{p\alpha}2\ge\frac{2p}{\pi}\sin\frac\alpha2$
compare the radial and angular terms separately. Otherwise,
$M^p-m^p\ge(1-2^{-p})M^p$ and $\ell\le\frac32M$.
Thus, with constants depending only on $K$,
\begin{gather}
 |w(x)-w(a)|\ge c_KM^{p-1}\ell,
 \label{conepdeno}\\
 \left|\nabla_p\arg(w-w(a))\right|
 =\frac{pr^{p-1}}{|w(x)-w(a)|}\le\frac{C_K}{\ell}.
 \label{conepkern}
\end{gather}
For the central term, $|\nabla_p\arg w|=\frac p r$.
The area bound $|B_t(a)|_p\le\pi p_+t^2$ and dyadic integration yield
\begin{equation*}
 \|\mathbf{P}_{0,n}\|_{L^\infty(\Omega_{n,R})}\le\frac CR,
 \qquad\int_{\Omega_{n,R}}|\mathbf{P}_{0,n}|^4\,dA_{p_n}\le\frac C{R^2}.
\end{equation*}
The constants are independent of the configuration diameter. On fixed
annuli about the centers, $\mathbf{P}_{0,n}$ tends to $d_j\nabla\arg y$ (centrally, $c\nabla_{p_*}\arg w$),
by \eqref{conepkern} and
$w\circ\iota_{j,n}(y)-w_{j,n}=\lambda_{j,n}y(1+O_K(\frac{|y|}{r(a_{j,n})}))$.

Let $v_n=\mathbf{P}_n-\mathbf{P}_{0,n}$ and
\[
 Q_n=\rho_n^2v_n-(1-\rho_n^2)\mathbf{P}_{0,n},\qquad
 E_{n,R}=\int_{\Omega_{n,R}}\rho_n^2|v_n|^2\,dA_{p_n}.
\]
Matching periods imply $v_n=\nabla_{p_n}\psi_n$, with a real function $\psi_n$ on $\Omega_{n,R}$. The exterior phase estimate and
$\mathbf{P}_{0,n}=\frac{p_n}{r}e_\theta+O_n(r^{-p_n-1})$ imply
$v_n,Q_n\in L^2(\Omega_{n,R})$. The imaginary part of the equation,
together with the representation of $\mathbf{P}_{0,n}$ as the rotated gradient of
$c\log|w|+\sum_jd_j\log|w-w_{j,n}|$,
yields $\operatorname{div}_{p_n}Q_n=0$ and zero normal boundary integral around each hole.
Testing with the phase, with $\nu_\Omega$ the outward normal, yields
\begin{align*}
 E_{n,R}&=\int_{\Omega_{n,R}}(1-\rho_n^2)\mathbf{P}_{0,n}\cdot v_n\,dA_{p_n}
              +\sum_{j=0}^JI_{j,n}(R),\\
 I_{j,n}(R)&=\int_{\partial B_R(a_{j,n})}\psi_nQ_n\cdot\nu_\Omega\,ds,
 \qquad a_{0,n}=o.
\end{align*}
Indeed, insert a cutoff at $T,2T$ and subtract the mean of $\psi_n$ on
$A_{n,T}=\{x\in C_{p_n}:T<r(x)<2T\}$. Constant subtraction leaves the boundary integrals unchanged. The Poincar\'e inequality on the annulus bounds the outer term by
$C_K\|v_n\|_{L^2(A_{n,T})}\|Q_n\|_{L^2(A_{n,T})}\to0$ as $T\to\infty$.
The volume term is at most
$$\frac CR E_{n,R}^{\frac12}\le\frac12E_{n,R}+\frac C{R^2}. $$

For the $j$th limit let $\chi_j,\rho_j,\Gamma_{j,R},\nabla_j,ds_j$ denote its
phase correction, modulus, circle of radius $R$, gradient, and arc measure
(on $C_{p_*}$ if $j=0$, on $\R^2$ otherwise).
On each fixed circle, the tangential derivative of $\psi_n$ converges
to that of $\chi_j$. Since the limit of $\mathbf{P}_{0,n}$ is tangential,
$Q_n\cdot e_r\to\rho_j^2\partial_r\chi_j$ there.
Since each normal boundary integral is zero, we subtract the mean on each circle.
Local convergence and circular Poincar\'e inequalities therefore imply
\begin{gather*}
 I_{j,n}(R)\to-\int_{\Gamma_{j,R}}(\chi_j-\langle\chi_j\rangle_R)
 \rho_j^2\partial_r\chi_j\,ds_j,\\
 \lim_n|I_{j,n}(R)|\le C_Kb_j(R),\qquad
 b_j(R)=R\int_{\Gamma_{j,R}}|\nabla_j\chi_j|^2\,ds_j.
\end{gather*}
The exterior phase estimates and coarea imply
\[
 \int_{R_0}^\infty\sum_{j=0}^Jb_j(R)\,\frac{dR}{R}<\infty,
 \qquad \limsup_nE_{n,R}\le\frac C{R^2}+C_K\sum_{j=0}^Jb_j(R).
\]
Choose $R_k\to\infty$ with $\sum_jb_j(R_k)\to0$ and use monotonicity in
$R$ to conclude
\begin{equation}\label{phaseinfty}
 \lim_{R\to\infty}\limsup_nE_{n,R}=0.
\end{equation}
By \eqref{bernboun}, $\mathbf{P}_n,v_n$ are uniformly bounded. Hence
$\int_{\Omega_{n,R}}|\mathbf{P}_n|^4\,dA_{p_n}\le\frac C{R^2}+CE_{n,R}$.
For $h_n=1-\rho_n$, test
$(-\Delta_{p_n}+\rho_n(1+\rho_n))h_n=\rho_n|\mathbf{P}_n|^2$
with $\zeta_{n,R}^2h_n$, where $\zeta_{n,R}$ vanishes on the disks of radius $R$,
equals one on $\Omega_{n,2R}$, and $|\nabla_{p_n}\zeta_{n,R}|\le\frac CR$.
An outer cutoff is removed first, using $h_n\in L^2(\Omega_{n,R})$.
Young's inequality, $\rho_n(1+\rho_n)\ge\frac34$, and
$h_n^2\le\mu_{U_n}\le4h_n^2$ yield
\[
 \int_{\Omega_{n,2R}}\mu_{U_n}\,dA_{p_n}
 \le C\int_{\Omega_{n,R}}|\mathbf{P}_n|^4\,dA_{p_n}
       +\frac C{R^2}\mathcal P_{p_n}(U_n)
 \le\frac C{R^2}+CE_{n,R}\to0
\]
in the iterated upper limit $n\to\infty$, $R\to\infty$.
For each fixed $R$ and every sequence $\ell_n$ satisfying the stated conditions, we deduce from local convergence that
\begin{align*}
 \limsup_n\left\|\mu_{U_n}-\mu_{0,n}-\sum_{j=1}^J\mu_{j,n}
                  \right\|_{L^1(C_{p_n})}
 &\le\limsup_n\int_{\Omega_{n,R}}\mu_{U_n}\,dA_{p_n}
       +\int_{C_{p_*}\setminus B_R(o)}\mu_{U_0}\,dA_{p_*}\\
 &\quad+\sum_{j=1}^J\int_{\R^2\setminus B_R(0)}\mu_{V_j}\,dx\to0.
\end{align*}
This proves \eqref{stroddeco}. Integration and
$dA_{p_n}=\frac{p_*}{p_n}dA_{p_*}$ imply \eqref{degrprela}.
For $J\ge1$, $p_*(1-c^2)=\sum_jd_j^2>0$ implies $c=0$, so $U_0$ is a unit
constant, and $p_*-1=\sum_jd_j(d_j-1)\in2\mathbb Z$.
For $J=0$, $c=1$ and
\[
 \|(1-|U_n|^2)-(1-|U_0|^2)\|_{L^2(C_{p_*})}^2
 \le\|\mu_{U_n}-\mu_{U_0}\|_{L^1(C_{p_*})}\to0,
\]
by $(a-b)^2\le|a^2-b^2|$ for $a,b\ge0$.
\end{proof}
 
For nonempty $K\Subset(2,\infty)$ let $\mathscr R_K$ be the set of pairs
$(p,U)$ such that $p\in K$, $U$ is a   bounded solution on $C_p$,
$|U|\to1$ uniformly at infinity, $\deg(U,\infty)=1$, $\mathcal P_p(U)<\infty$, and
\begin{equation*}
 U(r,-\phi)=\overline{U(r,\phi)},\qquad
 \|U(r,\cdot)-e^{i\cdot}\|_{C^1(S^1)}\to0\quad(r\to\infty).
\end{equation*}
Write $\Omega_p^+=\{(r,\phi)\in C_p:r>0,\ 0<\phi<\pi\}$ and
$\mathscr P_K=\{(p,U)\in\mathscr R_K:\operatorname{Im}U>0\text{ on }\Omega_p^+\}$.
The topology consists of parameter convergence, uniform convergence on
closed vertex balls, convergence in
$C^1(B_t(0)\setminus\overline{B_s(0)})$ for every $0<s<t<\infty$
in $y=re^{i\phi}$, and
\begin{equation}\label{posictopo}
 \tau(U_n,U):=\sup_{r\ge1}
       \|U_n(r,\cdot)-U(r,\cdot)\|_{C^1(S^1)}\to0.
\end{equation}
For example, with $A_j=\{y\in\C:\frac1j\le|y|\le j\}$ it is metrized by
\[
 d_K((p,U),(q,V))=|p-q|+\min\{1,\tau(U,V)\}
 +\sum_{j\ge2}2^{-j}\min\{1,\|U-V\|_{C^0(\overline{B_j(o)})}
                                  +\|U-V\|_{C^1(A_j)}\}.
\]

The following proposition shows that the positivity condition is preserved under limits and sufficiently small perturbations within $\mathscr R_K$.
\begin{proposition}
\label{posiclop}
Let $K\Subset(2,\infty)$ be nonempty, and equip $\mathscr R_K$ with the
topology defined above. Then $\mathscr P_K$ is both open and closed in
$\mathscr R_K$. For every $(p,U)\in\mathscr P_K$,
\[
 \{x\in C_p:U(x)=0\}
 \subset\{o\}\cup\{(r,\phi)\in C_p:r>0,\ \phi\in\{0,\pi\}\}.
\]
\end{proposition}
\begin{proof}
For $b=\operatorname{Im}U$, $-\Delta_pb=(1-|U|^2)b$ and $|U|\le1$.
A limit of positive solutions has $b\ge0$ and
$b(r,\frac\pi2)\to1$. The strong minimum principle proves closedness.

For openness, Hopf's lemma and the normalization imply
$m=\inf_{r\ge1,\,0<\phi<\pi}\frac{b(r,\phi)}{\sin\phi}>0$.
Indeed, Hopf controls each compact annulus, while at infinity use
\[
 \frac{|h(\phi)|}{\sin\phi}\le\frac\pi2\|h'\|_{L^\infty((0,\pi))},
 \qquad h\in C^1([0,\pi]),\quad h(0)=h(\pi)=0.
\]
Thus $\tau(U,V)<\frac m\pi$ implies
$\operatorname{Im}V(r,\phi)>\frac m2\sin\phi$ for $r\ge1$.
On $B_q^+=\Omega_q^+\cap B_1(o)$, the   weak test
$v=(-\operatorname{Im}V)_+$ has zero boundary trace. By the angular Poincar\'e inequality, we obtain
\[
 \int_{B_q^+}|\nabla_qv|^2\,dA_q
 =\int_{B_q^+}(1-|V|^2)v^2\,dA_q
 \le\int_{B_q^+}v^2\,dA_q
 \le\frac1{q^2}\int_{B_q^+}|\nabla_qv|^2\,dA_q.
\]
Since $q>2$, we have $v=0$. The strong minimum principle proves positivity
throughout $\Omega_q^+$. The stated inclusion for the zero set then follows by reflection.
\end{proof}

To identify the nonconstant translated limits on $\R^2$,
we first show that reflection symmetry and positivity force
their degrees to be $\pm1$. Since these limits have finite
potential energy, Mironescu's theorem
\cite[Theorem~1, p.~594]{Mironescu96} then implies that each
is a translate and constant phase rotation of
$V_{\pm1}=f_1(r)e^{\pm i\vartheta}$.
This argument requires neither the phase normalization at
infinity nor convergence in \eqref{posictopo}.

\begin{lemma}\label{unitdlimi}
Let $K\Subset(2,\infty)$ be nonempty, let $p_n\in K$, and let $U_n$ be
  bounded solutions on $C_{p_n}$ with limiting modulus one,
degree one at infinity, finite potential integral, and
\[
 U_n(r,-\phi)=\overline{U_n(r,\phi)},\qquad
 \operatorname{Im}U_n>0\quad\text{on }\Omega_{p_n}^+.
\]
Suppose that $\iota_n:B_{R_n}(0)\to C_{p_n}$ are isometries preserving orientation from Euclidean disks,
$R_n\to\infty$, $r(\iota_n(0))\to\infty$, and
$U_n\circ\iota_n\to V$ in $C^\infty_{\rm loc}(\R^2)$, where $V$ is a
nonconstant solution on $\R^2$ with
$\int_{\R^2}(1-|V|^2)^2\,dx<\infty$. Then
\[
 \deg(V,\infty)=d\in\{-1,1\},\qquad
 V(z)=e^{i\beta}V_d(z-z_0)
\]
for some $\beta\in\R$ and $z_0\in\C$.
\end{lemma}
\begin{proof}
By Lemma~\ref{planfpote}, we have $d=\deg(V,\infty)\ne0$.
Degree preservation on a large fixed circle provides zeros
$y_n\to y_0$ of $U_n\circ\iota_n$. Their images in the original variables escape along the reflection
boundaries. Rotate the coordinates so upper half disks lie in the positive angular sector. A subsequence converges to $W(z)=V(y_0+e^{i\theta}z)$ with
\[
 W(0)=0,\qquad W(\bar z)=\overline{W(z)},\qquad
 \operatorname{Im}W\ge0\quad(\operatorname{Im}z>0),\qquad \deg(W,\infty)=d.
\]
The strong minimum principle implies $\operatorname{Im}W>0$ above the real
axis, since a real limit has degree zero.
On a large upper semicircle write $W=|W|e^{ia}$ with
$0<a(\vartheta)<\pi$ and $a(0),a(\pi)\in\{0,\pi\}$. Reflection yields
\[
 2\pi d=2\int_0^\pi a'(\vartheta)\,d\vartheta
       =2(a(\pi)-a(0)),\qquad d\in\{-1,1\}.
\]
By \cite[Theorem~1, p.~594]{Mironescu96},
$V=e^{i\beta}V_d(\,\cdot-z_0\,)$ for some $\beta\in\R$ and $z_0\in\C$.
\end{proof}

Consequently an escaping positive reflection sequence satisfies
\begin{equation}\label{posiearit}
 c=0,\qquad p_*=J,\qquad \sum_{j=1}^Jd_j=1,\qquad
 d_j\in\{-1,1\},\qquad J\in\{3,5,7,\ldots\}.
\end{equation}
\section{Potential energy estimates and vortex escape}\label{separation}

 This section studies the potential energy of sequences whose vortices
separate from one another and move away from the vertex.
For a sequence in the positive reflection class with $p_n\to p_*>2$,
Proposition~\ref{poteddeco} and Lemma~\ref{unitdlimi} describe every
possible escape after taking a subsequence. The limit on bounded sets
is a constant of modulus one, and the translated limits are unit vortices.
If there are $N$ such
limits, each has potential energy $2\pi$, and their degrees $d_j=\pm1$
add up to one. The total potential energy $2\pi p_n$ therefore tends
to $2\pi N$, so $p_*=N$. There is one more positive vortex than
negative vortices, and hence $N$ is odd. Thus the limiting parameter
of any escaping sequence must be an odd integer. This already excludes
escape when $8<p_*<9$, in particular at any fixed $p\in(8,9)$.

For varying parameters in $[8,9]$, the remaining possibility is
$p_n\to9$. The main task of this section is to exclude this possibility.
More generally, Proposition~\ref{uniteside} shows that an escaping
sequence of $N$ unit vortices with $p_n\to N$ satisfies
$\mathcal P_{p_n}(U_n)>2\pi N$ for all sufficiently large $n$.
The proof combines an identity for the potential energy with a
logarithmic lower bound for the phase correction, whose coefficient
is bounded away from zero. Since $\mathcal P_{p_n}(U_n)=2\pi p_n$, we obtain
$p_n>N$: an escaping sequence can approach its limiting odd integer
only from above. When $N=9$, this contradicts $p_n\le9$.
Thus vortex escape is excluded throughout $[8,9]$.

Fix a nonempty compact set $K\Subset(2,\infty)$. For $p\in K$, an
integer $N\ge2$, distinct points $c_1,\ldots,c_N\in C_p\setminus\{o\}$,
and signs $d_j\in\{-1,1\}$, put
\begin{equation}\label{poingeom}
 Q=\sum_jd_j,\quad t=\min_jr(c_j),\quad
 L=\min_{j\ne k}\operatorname{dist}_p(c_j,c_k),\quad
 L_*=\min(t,L),\quad r_b(x)=\min_j\operatorname{dist}_p(x,c_j).
\end{equation}
Write $a_j=w(c_j)$, where $w=r^pe^{i\phi}$, and in a local oriented
Euclidean coordinate $z=re^{\frac{i\phi}p}$ define
\begin{equation}\label{poinfdefi}
 K_j=\frac{pz^{p-1}}{z^p-a_j},\qquad F=\sum_jd_jK_j,\qquad
 T_1=\sum_jK_j,\qquad \mathbf{P}_{0,x}-i\mathbf{P}_{0,y}=-iF.
\end{equation}
Thus $\mathbf{P}_0=\nabla_p\sum_jd_j\arg(w-a_j)$ is a global closed one form with
periods $2\pi d_j$.

The functions $K_j$ encode the derivatives of the local arguments
$\arg(w-a_j)$. The next lemma bounds these functions and their derivatives
in terms of distance on the cone, with constants independent of the
center positions. These estimates control the explicit field $\mathbf{P}_0$
away from the centers, including near the vertex.

\begin{lemma}
\label{poinfboun}
Let $p\in K$, let $c_1,\ldots,c_N\in C_p\setminus\{o\}$ be distinct,
and let $d_1,\ldots,d_N\in\{-1,1\}$. Define $K_j,F,T_1$ by
\eqref{poinfdefi}. For each $j$ and $x\in C_p\setminus\{o,c_j\}$,
write $r=r(x)$, $s=r(c_j)$ and $\ell=\operatorname{dist}_p(x,c_j)$.
There are constants $c_K,C_K>0$, independent of the centers, such that
\begin{align}
 |w(x)-a_j|&\ge c_K\max(r,s)^{p-1}\ell,\label{poindeno}\\
 |K_j|&\le C_K\ell^{-1},\quad |K_j'|\le C_K\ell^{-2},\quad
 \frac{r^{p-2}s}{|w(x)-a_j|}\le C_K\ell^{-1}.
 \label{poinkboun}
\end{align}
On $C_p\setminus\{o,c_1,\ldots,c_N\}$, consequently $|F|\le C_{K,N}r_b^{-1}$, $|T_1'|\le C_{K,N}r_b^{-2}$,
and
\begin{equation}\label{poinsderi}
 \left|\nabla_p^k(|F|^2)\right|\le C_{K,N}r_b^{-2-k},\qquad k=0,1,2.
\end{equation}
\end{lemma}
\begin{proof}
Put $M=\max(r,s)$. The general estimate
\eqref{conepdeno} implies \eqref{poindeno},
while \eqref{conepkern} yields $|K_j|\le C_K\ell^{-1}$.
Differentiating the kernel and using \eqref{poindeno}, we obtain
$|K_j'|\le C_K((M\ell)^{-1}+\ell^{-2})$.
Since $\ell\le2M$ and $r^{p-2}s\le M^{p-1}$, the remaining
bounds in \eqref{poinkboun} follow.

Put $b=r_b(x)$. For $r\ge\frac b4$, Cauchy's estimates on a Euclidean disk
$B_{c_Kb}(x)$ yield $|\nabla_p^kF|\le Cb^{-1-k}$. For $r<\frac b4$,
$r(c_j)>\frac{3b}4$ and the convergent expansion
$K_j=-p\sum_{m\ge0}z^{p(m+1)-1}a_j^{-m-1}$ implies
\begin{equation*}
 |\nabla_p^kF|\le C r^{p-1-k}b^{-p},\qquad k=0,1,2.
\end{equation*}
The product rule proves \eqref{poinsderi}.
In the $y$ coordinates the Christoffel symbols satisfy
$|\Gamma_p|\le\frac{C_K}{r}$. Near the vertex, the same expansion yields
\[
 r^{-2}|F|^2+r^{-1}|\nabla_y(|F|^2)|+|\nabla_y^2(|F|^2)|
 \le Ct^{-2p}r^{2p-4}\to0 .
\]
Also $K_j=O(r^{p-1})$ and $K_j'=O(r^{p-2})$. This finishes the proof.
\end{proof}

Near each center, $F$ consists of the singular term of an isolated unit
vortex and a holomorphic remainder. The next lemma describes the
constant and linear terms of this remainder and derives algebraic
identities for $F^2$ and the constant terms. These formulas connect
the local expansions with the integral estimates below and with the
distance of the centers from the vertex.

\begin{lemma}
\label{poinalge}
Let $p\in K$, let $c_1,\ldots,c_N\in C_p\setminus\{o\}$ be distinct,
and let $d_j\in\{-1,1\}$. Use the definitions in
\eqref{poingeom}--\eqref{poinfdefi}. In a
local Euclidean coordinate chart $z=z_j+y$ centered at $c_j$, where $z_j^p=a_j$, set
\begin{equation}\label{poininci}
 E_j=\sum_{k\ne j}\frac{d_k}{a_j-a_k},\qquad
 A_j=pz_j^{p-1}E_j+\frac{d_j(p-1)}{2z_j}.
\end{equation}
There are coefficients $\beta_j\in\C$ and constants $c_K,C_{K,N}>0$ such that
\begin{equation}\label{poinlexpa}
 F(z_j+y)=\frac{d_j}{y}+A_j+\beta_jy+O(L_*^{-3}|y|^2),
 \quad |A_j|\le CL_*^{-1},\quad|\beta_j|\le CL_*^{-2},
 \quad 0<|y|<c_KL_*.
\end{equation}
Moreover,
\begin{align}
 \label{pointB}
 F^2 =-T_1'+B,\qquad 
 B =2pz^{p-2}\sum_j\frac{d_jz_jA_j}{z^p-a_j},
  \qquad 
 2\sum_jd_jz_jA_j=pQ^2-N.
\end{align}
The differential $B(z)(\,dz)^2$ and the numbers $z_jA_j$ are
independent of the branch of the local Euclidean coordinate.
\end{lemma}
\begin{proof}
The expansion $K_j(z_j+y)=\frac1y+\frac{p-1}{2z_j}+O(t^{-2}|y|)$
determines $A_j$. Cauchy's estimates for the regular part
$F-\frac{d_j}y=O(L_*^{-1})$ on $B_{2c_KL_*}(c_j)$ prove
\eqref{poinlexpa}.

Using the method of partial fractions of \cite[\S2]{Kurzke}, we obtain, with ordered sums,
\begin{gather*}
 \left(\sum_j\frac{d_j}{w-a_j}\right)^2
 =\sum_j\frac1{(w-a_j)^2}+2\sum_j\frac{d_jE_j}{w-a_j},\\
 \sum_jd_jE_j=0,\quad
 \sum_jd_ja_jE_j=\frac{Q^2-N}{2}.
\end{gather*}
Using \eqref{poinfdefi} and
$z_jA_j=pa_jE_j+\frac12d_j(p-1)$, we deduce the identities, including
$2\sum_jd_jz_jA_j=p(Q^2-N)+(p-1)N=pQ^2-N$.
Under $\widetilde z=\omega z$, $\omega=e^{\frac{2\pi im}p}$,
$F,A_j$ acquire $\omega^{-1}$ and $B$ acquires $\omega^{-2}$,
which proves the invariance.
\end{proof}

\Needspace{6\baselineskip}The identity $F^2=-T_1'+B$ reduces the study of $|F|^4$ to integrals
involving $T_1'$ and $B$. We proceed to bound the integral of $|B|^2$
in terms of the constants $A_j$, uniformly in the configuration
diameter. This estimate controls the $|B|^2$ term when we expand
$|F^2|^2$ in Proposition~\ref{conefpart}. For the quantities in
\eqref{poingeom}, \eqref{poininci} and \eqref{pointB}, set
$A_*=\max_{1\le j\le N}|A_j|$ and
$\Omega_R=\{x\in C_p:r_b(x)>R\}$.

\begin{lemma}\label{pointBdual}
For every $p\in K$, distinct centers $c_j\in C_p\setminus\{o\}$,
signs $d_j\in\{-1,1\}$, and $R>0$, one has, with $C=C_{K,N}$,
\begin{equation*}
 |B|\le C\min\left(\frac{A_*}{r_b},r_b^{-2}\right),\qquad
 \int_{\Omega_R}|B|^2\,dA_p
 \le CA_*^2\left(1+\log_+\frac1{A_*R}\right).
\end{equation*}
The right side is zero if $A_*=0$. For fixed $R,M>0$ and any sequence of such configurations with
$L\to\infty$ and $A_*\le ML^{-3}$, the integral is
$O_{K,N,R,M}(L^{-6}\log L)$.
\end{lemma}
\begin{proof}
The two formulas for $B$ yield $|B|\le C\min(A_*r_b^{-1},r_b^{-2})$.
We cut opposite each center and use Euclidean coordinates on the resulting sectors. Enlarging
to the plane, we obtain
\[
 \int_{\Omega_R}|B|^2\,dA_p
 \le C\int_R^\infty\min(A_*^2s^{-1},s^{-3})\,ds
 \le CA_*^2\left(\frac12+\log_+\frac1{A_*R}\right).
\]
For $A_*>0$, we split the integral at $\frac1{A_*}$.
If $A_*=0$, then $B=0$.
The last assertion follows from monotonicity of $a^2(1+\log\frac1{aR})$ for $0<a<\frac1R$.
\end{proof}
 
Fix an odd $N\ge3$ and $p_n\to N$ in $K\Subset(2,\infty)$.
Throughout the remainder of this section, $U_n$ are bounded  
solutions of $-\Delta_{p_n}U_n=U_n(1-|U_n|^2)$ on $C_{p_n}$, with
$\deg(U_n,\infty)=1$, $\mathcal P_{p_n}(U_n)<\infty$ and $|U_n|\to1$ at infinity. Assume there are $N$ translated limits with local degrees $d_j=\pm1$, centers
$c^0_{j,n}$, and phases $\gamma^0_{j,n}$, $1\le j\le N$. Here $t_n^0$
and $L_n^0$ are the vertex distance and minimum pairwise separation
in \eqref{poingeom} for these unmodulated centers. Assume that
\begin{equation}\label{uniteassu}
 t_n^0\to\infty,\qquad L_n^0\to\infty,\qquad
 e^{-i\gamma^0_{j,n}}U_n(c^0_{j,n}+y)\to
 V_{d_j}(y)=f(|y|)e^{id_j\arg y}
 \quad\hbox{in }C^\infty_{\rm loc}(\R^2).
\end{equation}
The translated limits in \eqref{posiearit} satisfy
these assumptions. The argument below applies independently of reflection
symmetry and of any upper bound for the diameter. By
\cite[pp.~428, 433, 435--437]{HH}, the radial function $f$ satisfies $f'>0$,
$f(r)=\kappa r+O(r^3)$, and
\begin{equation*}
 f=1-\frac1{2r^2}+O(r^{-4}),\qquad f'=O(r^{-3}),\qquad f''=O(r^{-4}).
\end{equation*}
Also $\mathcal P(V_{\pm1})=2\pi$ by
Lemma~\ref{planfpote}. Put $\rho_n=|U_n|$ and, wherever
$U_n\ne0$, let $\mathbf{P}_n=\nabla_{p_n}\arg U_n$. Let $\mathbf{P}_{0,n}$ be the explicit phase gradient
\eqref{poinfdefi} associated with $c^0_{j,n},d_j$, and
write $v_n=\mathbf{P}_n-\mathbf{P}_{0,n}$.

The next lemma confines all zeros of $U_n$ to fixed disks about the
chosen centers and proves that the modulus tends to one whenever the
distance from every center tends to infinity. Outside these disks,
$v_n=\mathbf{P}_n-\mathbf{P}_{0,n}$ is the gradient of a single valued real function,
and the integral of $\rho_n^2|v_n|^2$ tends to zero.

\begin{lemma}\label{unitpener}
For the sequences, centers, degrees and phases fixed above and satisfying
\eqref{uniteassu}, there exist $R_0>0$ and $n_0$ such that, for $n\ge n_0$,
\[
 \Omega_{n,R_0}=C_{p_n}\setminus\bigcup_{j=1}^N\overline{B_{R_0}(c^0_{j,n})},
 \qquad \rho_n\ge\frac12\ \hbox{on }\Omega_{n,R_0},\qquad
 \sum_{j=1}^Nd_j=1.
\]
All zeros of $U_n$ lie in the deleted disks. There is a single valued
real function $\psi_n$ on $\Omega_{n,R_0}$ such that
$v_n=\nabla_{p_n}\psi_n\in L^2(\Omega_{n,R_0})$, and
\begin{equation*}
 E_n=\int_{\Omega_{n,R_0}}\rho_n^2|v_n|^2\,dA_{p_n}\to0.
\end{equation*}
Moreover, every sequence $x_n\in C_{p_n}$ satisfying
$r_{b,n}(x_n)\to\infty$ satisfies $\rho_n(x_n)\to1$.
\end{lemma}
\begin{proof}
We adapt the modulus and phase argument of
\cite[\S2, Step~1 and Proposition~1]{BMR} to the present domains.
We subtract the convergent integrals over the fixed disks from $\mathcal P_{p_n}(U_n)=2\pi p_n$.
Using $\mathcal P(V_{\pm1})=2\pi$, we then obtain
\begin{equation*}
 \limsup_n\int_{\Omega_{n,R}}(1-\rho_n^2)^2\,dA_{p_n}
 \le \tau(R):=2\pi N\int_R^\infty(1-f^2)^2s\,ds\to0.
\end{equation*}
If $1-\rho_n(x)\ge\delta$, the uniform Lipschitz bound $M_K$ implies
\[
 \int_{B_{\frac{\delta}{2M_K}}(x)}(1-\rho_n^2)^2\,dA_{p_n}
 \ge\frac{\pi\delta^4}{16p_+M_K^2}.
\]
The preceding integral bound proves the assertions about the modulus.
Degree additivity implies $\sum_jd_j=1$ and zero periods for $v_n=\nabla_{p_n}\psi_n$.
The estimates at the vertex imply $v_n=O_n(r)$. Proposition~\ref{quantize}
and $\mathbf{P}_{0,n}-\nabla_{p_n}\phi=O_n(r^{-p_n-1})$ yield
$v_n\in L^2(\Omega_{n,R_0})$.

Suppress $n$ and the superscript $0$. The field
$Q=\rho^2v-(1-\rho^2)\mathbf{P}_0=j(U)-\mathbf{P}_0\in L^2(\Omega_R)$
has zero divergence and has zero normal boundary integral on every hole.
Testing with $\psi$ and denoting its circle means by $m_j$, we obtain
\begin{equation}\label{unitpiden}
 \int_{\Omega_R}\rho^2|v|^2\,dA_p
 =\int_{\Omega_R}(1-\rho^2)\mathbf{P}_0\cdot v\,dA_p
 +\sum_j\int_{\partial B_R(c_j)}
       (\psi-m_j)Q\cdot\nu_{\Omega_R}\,ds_p .
\end{equation}
Use a radial cutoff $\eta_T$ between $B_T(o)$ and $B_{2T}(o)$,
with $|\nabla_p\eta_T|\le\frac CT$. We subtract the mean $m_T$ on
$A_T=B_{2T}(o)\setminus\overline{B_T(o)}$ and apply the Poincar\'e inequality to obtain
\[
 \left|\int_{A_T}(\psi-m_T)Q\cdot \nabla_p\eta_T\,dA_p\right|
 \le C_K\|v\|_{L^2(A_T)}\|Q\|_{L^2(A_T)}\to0 .
\]
At the vertex, $Q=O_n(r)$ and $\psi-\psi(o)=O_n(r^2)$ make the
boundary term $O_n(r^4)$. Thus the test is justified.

For fixed $R$, local convergence in translated coordinates makes the boundary sum $o_n(1)$.
Using $\rho_n\ge\frac12$ and $|\mathbf{P}_{0,n}|\le\frac CR$,
we absorb half the left side by Young's inequality and obtain
\[
 \limsup_{n\to\infty}
 \int_{\Omega_{n,R}}\rho_n^2|v_n|^2\,dA_{p_n}
 \le\frac C{R^2}\tau(R).
\]
The integrals on the fixed annuli about the centers between $R_0$ and $R$ tend to zero.
Consequently $\limsup_nE_n\le\frac C{R^2}\tau(R)\to0$ as $R\to\infty$.
\end{proof}

Pointwise estimates throughout the cone require a gradient bound on
disks that contain the vertex. The next lemma estimates the gradient
of a solution to the Poisson equation in terms of the boundary
oscillation and the bounded right side. 

\begin{lemma}\label{tippgrad}
Let $p\in K\Subset(2,\infty)$, let $h\in L^\infty(B_1(o))$ with
$\|h\|_{L^\infty(B_1(o))}\le M$, and let
a real function $\chi\in H^1_{\rm loc}(B_1(o))$ solve $\Delta_p\chi=h$ naturally
across the vertex.
Suppose the conformal representative $\widehat\chi$ has distribution
boundary trace $g\in L^\infty(S^1)$ of essential oscillation at most $B$.
Then
\begin{equation*}
 |\nabla_p\chi(r,\phi)|\le C_K(Br^{p-1}+Mr),\qquad 0<r\le\frac12.
\end{equation*}
For every $0<r<1$ and $0<\alpha<1$,
$\chi|_{\partial B_r(o)}\in C^{1,\alpha}(\partial B_r(o))$.
\end{lemma}
\begin{proof}
In conformal coordinates the right side
$p^{-2}|w|^{\frac2p-2}\widehat h$ belongs to $L^s(B_1(0))$
for $1<s<\frac{p_+}{p_+-1}$. Its Green potential with zero boundary values $v$ satisfies
\[
 |\nabla v(w)|\le C_KM\left(1+\int_{B_1(0)}
       \frac{|\zeta|^{\frac2p-2}}{|w-\zeta|}\,d\zeta\right)
 \le C_KM|w|^{\frac2p-1}\quad\left(0<|w|\le\frac14\right).
\]
Split at $\frac{|w|}2$ and $2|w|$ for the last inequality.
By \cite[Corollary~1.12]{FRO}, the gradient of the harmonic remainder
is bounded by $CB$ on $B_{\frac14}(0)$. Multiplication by $pr^{p-1}$
proves the estimate.
The traces follow from interior elliptic regularity on annuli \cite[Proposition~2.18]{FRO}.
\end{proof}

We now turn the integral convergence in Lemma~\ref{unitpener} into
pointwise estimates outside the vortex disks. The following result quantifies
the approach of the modulus to one and bounds
$\mathbf{P}_n=\nabla_{p_n}\arg U_n$ in terms of the distance to the nearest center.
It also proves that
$r_{b,n}|v_n|$ tends to zero uniformly, which will control the boundary
integrals in the following proposition.
\begin{lemma}
\label{unitfregi}
Assume the hypotheses of Lemma~\ref{unitpener} for the sequences
$p_n,U_n,c^0_{j,n},d_j,\gamma^0_{j,n}$. Use its notation
$R_0,\rho_n,\mathbf{P}_n,v_n,E_n$.
There are $C>0$ and $n_0$ such that, for every $n\ge n_0$ and
$x\in\Omega_{n,R_0}$,
\begin{equation}\label{unitfboun}
 \begin{gathered}
 |\mathbf{P}_n|\le\frac C{r_{b,n}},\qquad
 1-\rho_n+|\nabla_{p_n}\rho_n|\le\frac C{r_{b,n}^2},\\
 r_{b,n}|v_n|\le C\left(\sqrt{E_n}+\frac1{r_{b,n}}\right).
 \end{gathered}
\end{equation}
The constant is independent of the configuration diameter. Moreover,
\begin{equation*}
 \delta_{{\rm ph},n}:=\sup_{x\in\Omega_{n,R_0}}r_{b,n}(x)|v_n(x)|\to0.
\end{equation*}
\end{lemma}
\begin{proof}
We suppress $n$ and put $s=r_b(x)\ge2R_0$. We fix $0<\kappa<\frac18$.
Since distance to a set is Lipschitz with constant one,
$(1-2\kappa)s\le r_b\le(1+2\kappa)s$ on $B_{2\kappa s}(x)$.
The bounds $|\mathbf{P}_0|\le\frac Cs$ and $\rho\ge\frac12$ yield
\begin{equation*}
 \int_{B_{2\kappa s}(x)}|\mathbf{P}|^2\,dA_p\le C\kappa^2+8E_n.
\end{equation*}
The equation satisfied by $h=1-\rho$ is
\begin{equation*}
 (-\Delta_p+\rho(1+\rho))h=\rho|\mathbf{P}|^2,\qquad
 \frac34\le\rho(1+\rho)\le2.
\end{equation*}
We test with $\eta^2h$, for a cutoff between $B_{\kappa s}(x)$ and
$B_{2\kappa s}(x)$ with $|\nabla_p\eta|\le\frac C{\kappa s}$, and use
$\int_{C_p}h^2\,dA_p\le2\pi p$ to obtain
\[
 \int_{B_{\kappa s}(x)}e(U)\,dA_p
 \le C\left(\kappa^2+E_n+\frac1{\kappa^2s^2}\right).
\]
Choose $\kappa$ small and then $n,s$ large. By Lemma~\ref{unifhein},
we then have $|\mathbf{P}(x)|^2\le8e(U)(x)\le Cs^{-2}$.

For   $\chi$ on $B_{8\ell}(x)$ we use
\begin{equation}\label{unitlpgrad}
 \ell|\nabla_p\chi(x)|\le C_K\left[
 \|\nabla_p\chi\|_{L^2(B_{8\ell}(x))}
 +\ell^2\|\Delta_p\chi\|_{L^\infty(B_{8\ell}(x))}\right].
\end{equation}
Choose a Euclidean disk or a disk centered at the vertex,
$B_R(x_*)\subset B_{8\ell}(x)$, with
$c_K\ell\le R\le C_K\ell$ and $x\in B_{\frac R4}(x_*)$. Let $q$ solve
$\Delta_pq=\Delta_p\chi$ with zero boundary values. Then $\chi-q$
is a   harmonic function. The quadratic barrier,
Lemma~\ref{tippgrad}, the energy estimate for $q$,
and the harmonic gradient bound prove \eqref{unitlpgrad}.
At the vertex, this bound follows from the expansion in
$r^{kp}\cos k\phi,r^{kp}\sin k\phi$. All constants are uniform for $p\in K$.

We now choose the disks with $c_Ks\le R\le C_Ks$. For the radial distance $t$,
$\Delta_pI_0(\frac t2)=\frac14I_0(\frac t2)$.
Since $\rho(1+\rho)\ge\frac34$, the function
$\frac{C_1}{s^2}+\frac{I_0(\frac t2)}{I_0(\frac R2)}$
has boundary values at least one. For large $C_1$, applying
$-\Delta_p+\rho(1+\rho)$ to it yields at least $\rho|\mathbf{P}|^2$.
The maximum principle and $\frac13e^{\frac b2}\le I_0(b)\le e^b$ imply
\[
 h(x)\le\frac{C_1}{s^2}+3e^{-\frac R8}\le Cs^{-2}.
\]
Cutoff estimates on disks of radius one and \eqref{unitlpgrad}
yield $|\nabla_p\rho(x)|\le Cs^{-2}$. The imaginary part of the equation reads
\begin{equation*}
 \Delta_p\psi=-2\nabla_p\log\rho\cdot \mathbf{P}.
\end{equation*}
Its right side is $O(s^{-3})$. Apply \eqref{unitlpgrad}
with $\ell=\frac{s}{32}$ and
$\|\nabla_p\psi\|_{L^2(\Omega_{n,R_0})}\le2\sqrt{E_n}$
to obtain the phase bound. Local convergence in translated coordinates handles fixed annuli.
On $r_b\ge R$, the supremum of $r_b|v_n|$ has limsup at most $\frac CR$.
Let $R\to\infty$.
\end{proof}

The estimates so far involve both the minimum distance to the vertex
and the minimum separation between centers. The following proposition
shows that the former is bounded below by a fixed multiple of the
latter. Consequently, the local estimates below can be carried out
on disjoint Euclidean disks whose radii are proportional to the
separation between centers.
\begin{proposition}
\label{distfvert}
Let $p_n,U_n,c^0_{j,n},d_j,\gamma^0_{j,n}$ satisfy the hypotheses of
Lemma~\ref{unitpener}. Then there are $c_*>0$ and $n_0$
such that $t_n^0\ge c_*L_n^0$ for every $n\ge n_0$.
\end{proposition}
\begin{proof}
Suppose $\frac tL\to0$, suppressing $n$ and the superscript $0$.
Choose $|z_j|=t$. For $k\ne j$, $r(c_k)\ge L-t$, so
\eqref{poininci} yields
\begin{equation}\label{tipsinci}
 A_j=\frac{d_j(p-1)}{2z_j}
 +O\!\left(\frac{t^{p-1}}{(L-t)^p}\right),\qquad
 t|A_j|\to\frac{N-1}2.
\end{equation}
We choose a Euclidean disk $B_{2R}(c_j)$, $R=\kappa_Kt$, disjoint from all other centers.
Put $S(X)=X\otimes X-\frac12|X|^2I$ and
\begin{equation*}
 T_0=S(\mathbf{P}_0),\qquad
 T(U)=\nabla U\mathbin{\otimes}\nabla U
 -\left(\frac12|\nabla U|^2+\frac14(1-|U|^2)^2\right)I.
\end{equation*}
The tensor $T(U)$ has zero divergence, so $\int_{\partial B_R(c_j)}T(U)\nu\,ds=0$.
We set $r_j=\operatorname{dist}_p(\cdot,c_j)$.
Write $\mathbf{P}_0=\frac{d_j}{r_j}e_\vartheta+H_j$ and let $J$ rotate by
$\frac\pi2$. The mixed term in $T_0\nu$ is $\frac{d_j}{r_j}JH_j$.
As the circle shrinks to the center, the integrals of both pure terms
tend to zero. Thus
\[
 \int_{\partial B_R(c_j)}T_0\nu\,ds=2\pi d_jJH_j(0),
 \qquad |H_j(0)|=|A_j|.
\]
On this circle, $r_b=R$, and
$T(U)=S(\nabla\rho)+\rho^2S(\mathbf{P})-\frac14(1-\rho^2)^2I$.
From \eqref{unitfboun}, we thus deduce
$|T(U)-T_0|\le C(R^{-4}+\delta_{{\rm ph},n}R^{-2})$, whence
\[
 2\pi|A_j|\le
 \int_{\partial B_R(c_j)}|T(U)-T_0|\,ds
 \le C(R^{-3}+\delta_{{\rm ph},n}R^{-1})=o(t^{-1}),
\]
contrary to \eqref{tipsinci}.
\end{proof}
 
We now fix the local phase and translation. Choose a nonnegative
nonzero radial weight $\omega\in C_c^\infty(B_{R_0}(0))$ and
set $Z_0=iV_d$, $Z_a=\partial_aV_d$ for $a=1,2$. Let us adjust each center and phase
by $o(1)$ so that
\begin{equation}\label{unitgram}
 \int_{B_{R_0}(0)}\omega\operatorname{Re}(W_j\overline{Z_\ell})\,dy=0,\quad
 W_j(y)=e^{-i\gamma_j}U(c_j+y)-V_{d_j}(y),\quad \ell=0,1,2.
\end{equation}
For the adjusted centers, we set $r_j(x)=\operatorname{dist}_p(x,c_j)$
and $f_j=f(r_j)$. In the Euclidean chart about $c_j$, let
$b_j=d_je_\vartheta/r_j$, where $e_\vartheta$ is the angular unit vector
centered at $c_j$.
For a small fixed $c>0$, take disjoint Euclidean disks $B_{8cL}(c_j)$ and
radial cutoffs $\chi_j=1$ on $B_{4cL}(c_j)$ with
$|\nabla^k\chi_j|\le C_kL^{-k}$. Define
\begin{equation*}
 \mathcal R=1+\sum_j\chi_j(f(r_j)-1),\qquad
 \Omega=C_p\setminus\bigcup_j\overline{B_{R_0}(c_j)},
\end{equation*}
and, using the modulated explicit phase gradient,
\begin{equation}\label{unitenorm}
 \delta=\max_j\|W_j\|_{C^0(B_{2R_0}(0))}
       +\|\rho-\mathcal R\|_{L^\infty(\Omega)}
       +\sup_{x\in\Omega}r_b(x)|v(x)|.
\end{equation}
The next proposition quantifies convergence in these coordinates and
bounds the regular phase gradient at each center.
\begin{proposition}\label{unitrates2}
Under the assumptions of Lemma~\ref{unitpener}, there are phase
and center corrections tending to zero for which
\eqref{unitgram} holds. For the resulting quantities in
\eqref{poingeom}, \eqref{poininci}, and
\eqref{unitenorm},
\begin{equation*}
 \delta_n=O(L_n^{-2}),\qquad
 A_{*,n}=O(L_n^{-3}),\qquad
 \|W_{j,n}\|_{C^k(\overline{B_R(0)})}\le C_{R,k}L_n^{-2}.
\end{equation*}
\end{proposition}

We remark that the last estimate holds for every fixed $R>0$, every integer $k\ge0$,
every $1\le j\le N$, and all sufficiently large $n$. The constants may
depend on the sequence and the fixed normalization, but are independent
of the largest center radius and the configuration diameter.
\begin{proof}
We first impose the orthogonality conditions, then prove
$\delta_n=O(L_n^{-2})$ by contradiction using the local limits and
estimates on the whole cone.

\emph{Step 1. Normalization and preliminary estimates.}
The derivative of the modulation equations \eqref{unitgram} tends to
$G\operatorname{diag}(-1,1,1)$, where, for $q\in\R^3$,
\begin{gather*}
 G_{k\ell}=\int_{\R^2}\omega\operatorname{Re}(Z_k\overline{Z_\ell})\,dy,
 \quad q^tGq=\int_{\R^2}\omega\Big|\sum_{\ell=0}^2q_\ell Z_\ell\Big|^2\,dy,\\
 \operatorname{Re}\Big(e^{-id\vartheta}\sum_{\ell=0}^2q_\ell Z_\ell\Big)
 =f'(r)(q_1\cos\vartheta+q_2\sin\vartheta).
\end{gather*}
Since $f'>0$ and $\omega\ge0$ is radial and nonzero, $G$ is positive
definite. By the implicit function theorem, we obtain corrections of size
$h_n=o(1)$. They satisfy $|L_n-L_n^0|+|t_n-t_n^0|\le3h_n$ and
\begin{gather*}
 |\mathbf{P}_{0,n}-\mathbf{P}_{0,n}^0|\le Ch_n\sum_j\operatorname{dist}_{p_n}(x,c_{j,n})^{-2},\\
 \|\mathbf{P}_{0,n}-\mathbf{P}_{0,n}^0\|_{L^2(\{x\in C_{p_n}:r_{b,n}(x)\ge2R_0\})}
 \le\frac{Ch_n}{R_0}.
\end{gather*}
The analogous weighted supremum bound and local convergence in translated coordinates on the
changed annuli imply $\delta_n\to0$ and $t_n\ge c_*L_n$.

Suppress $n$. On $B_{8cL}(c_j)$ write $\mathbf{P}_0=b_j+H_j$. Cauchy's estimate implies
\begin{equation}\label{unitrfiel}
 \begin{gathered}
 |\nabla^kH_j|\le C_kL^{-1-k}\ (k\ge1),\qquad
 H_j(y)=H_j(0)+L^{-2}M_jy+O(L^{-3}|y|^2),\\
 |H_j(0)|=|A_j|\le CL^{-1},\qquad
 M_j^t=M_j,\quad\operatorname{tr}M_j=0,\quad|M_j|\le C.
 \end{gathered}
\end{equation}
The preceding boundary identity and \eqref{unitfboun} yield
$2\pi|A_j|\le CL(\delta L^{-2}+L^{-4})$, hence
\begin{equation*}
 A_*\le CL^{-1}(\delta+L^{-2}).
\end{equation*}
Define
\begin{align*}
 \mathcal A=\Delta_p\mathcal R-\mathcal R|\mathbf{P}_0|^2
                   +\mathcal R(1-\mathcal R^2),
                   \quad\mbox{and}\quad
 \mathcal F=(1-\mathcal R^2)\mathbf{P}_0-\sum_j\chi_j(1-f_j^2)b_j.
\end{align*}
For $R_0\le r_j\le4cL$, we infer from the radial equation that
$\mathcal A=-f_j(2b_j\cdot H_j+|H_j|^2)$ and
$\mathcal F=(1-f_j^2)H_j$, so
\begin{equation*}
 |\nabla^k\mathcal A|\le C(A_*r_j^{-1-k}+L^{-2}r_j^{-k})\ (k\le2),
 \qquad |\mathcal F|\le C(A_*r_j^{-2}+L^{-2}r_j^{-1}).
\end{equation*}
We use the radial equation before differentiating the cutoff.
On $4cL\le r_j\le8cL$, we obtain
\begin{equation*}
 |\nabla^k\mathcal A|\le CL^{-2-k}\ (k\le2),\qquad
 |\mathcal F|\le CL^{-3}.
\end{equation*}
Outside the supports, $\mathcal F=0$, $\mathcal A=-|\mathbf{P}_0|^2$, and
\begin{equation*}
 |\nabla_p^k\mathcal A|\le Cr_b^{-2-k},\qquad k\le2.
\end{equation*}
The vector field
\begin{equation*}
 \mathcal J=\rho^2v+(\rho^2-\mathcal R^2)\mathbf{P}_0-\mathcal F
 =j(U)-\mathbf{P}_0+\sum_j\chi_j(1-f_j^2)b_j
\end{equation*}
has zero divergence away from the centers and vertex, since the cutoffs are radial.
Near each center, $\mathcal J=j(U)-f_j^2b_j-H_j$ is smooth. At the vertex,
$\int_{\partial B_\eta(o)}|\mathcal J|\,ds_p\le C_n(\eta^2+\eta^p)\to0$.
By excision, we therefore have
\[
 \operatorname{div}_p\mathcal J=0\text{ in }\mathcal D'(C_p),\qquad
 \int_\Gamma\mathcal J\cdot\nu\,ds_p=0
\]
for every smooth closed contour bounding a domain.

\emph{Step 2. Local limits and the modulus estimate.}
Suppose $\delta_nL_n^2\to\infty$ along a subsequence. Set
\begin{equation*}
 a=\frac{\rho-\mathcal R}{\delta},\quad V=\frac v\delta,\quad
 g=\frac{\mathcal A}{2\delta},\quad
 \tau=\frac{L^{-2}}\delta\to0,\quad
 \kappa=\frac{A_*}\delta\le\frac C L(1+\tau).
\end{equation*}
Then $|a|\le1$, $r_b|V|\le1$, and $\frac{H_j}{\delta}\to0$ on fixed disks.
Let $\Psi_j$ be the local phase correction relative to
$e^{i\gamma_j}V_{d_j}$ in the chart centered at $c_j$, with branch
fixed by the orthogonality conditions.
Integrating in the radial variable, we find
\[
 \sup_{R_0\le|y|\le R}|\Psi_j(y)|
 \le C\delta+\int_{R_0}^R C(\delta s^{-1}+A_*+L^{-2}s)\,ds\le C_R\delta.
\]
Together with the bound on the inner disk in \eqref{unitenorm}, this
bounds $\frac{W_j}{\delta}$ on fixed disks. Its exact equation
$-\Delta\frac{W_j}{\delta}=\int_0^1DF_0(V_{d_j}+tW_j)\frac{W_j}{\delta}\,dt$
and interior estimates yield a subsequential limit in
$C^\infty_{\rm loc}(\R^2)$ satisfying
\begin{equation*}
 -\Delta W+(|V_d|^2-1)W+
            2V_d\operatorname{Re}(\overline{V_d}W)=0.
\end{equation*}
Write $W=e^{id\vartheta}(A+if\eta)$. Outside a fixed disk,
$|A|\le1$, $r|\nabla\eta|\le1$, and the vector field obtained from the imaginary
part of the linearized equation extends smoothly across zero. The divergence theorem yields
\begin{gather*}
 0=\int_{\partial B_r(0)}\left(f^2\nabla\eta+\frac{2dfA}{r}e_\vartheta\right)\cdot\nu\,ds
 =2\pi rf^2\langle\eta\rangle',\\
 |\eta-\langle\eta\rangle|\le\int_0^{2\pi}|\partial_\vartheta\eta|\,d\vartheta\le2\pi.
\end{gather*}
Thus $W$ is bounded. By \cite[Theorem~3.2, p.~59 of the author version]{PR},
conjugated by $e^{i\vartheta}$ for $d=1$ and by complex conjugation for $d=-1$,
$W=\sum_{\ell=0}^2q_\ell Z_\ell$. The limiting Gram equations yield
$Gq=0$, hence $W=0$ and $\frac{W_{j,n}}{\delta_n}\to0$ in
$C^\infty_{\rm loc}(\R^2)$.

Outside the fixed disks, we subtract the equations for the modulus exactly and obtain
\begin{equation}\label{unitnampl}
 \begin{aligned}
 (-\Delta_p+2)(a-g)
 &=(3-|\mathbf{P}_0|^2-\rho^2-\rho\mathcal R-\mathcal R^2)a+\Delta_pg\\
 &\qquad-\rho(2\mathbf{P}_0\cdot V+\delta|V|^2).
 \end{aligned}
\end{equation}
The bounds for $\mathcal A$ and $\mathcal F$, with $s=r_b$, yield
\begin{equation*}
 |g|+s|\nabla_pg|+s^2|\nabla_p^2g|
 \le C
 \begin{cases}
  \frac\kappa s+\tau,&R_0\le r_j\le4cL,\\
  \tau,&4cL\le r_j\le8cL,\\
  \tau\left(\frac L s\right)^2,&\text{outside all cutoff supports}.
 \end{cases}
\end{equation*}
Here $s\ge cL$ in the last region. The coefficient of $a$ in
\eqref{unitnampl} is $O(s^{-2})$, so its right
side is $O(s^{-2})$.
On Euclidean disks or disks centered at the vertex, $B_R(z_0)$,
$c_Ks\le R\le C_Ks$, with the chosen point
in $B_{\frac R4}(z_0)$, the barrier
$b_R(t)=\exp\frac{t^2-R^2}{4R}$, $t=\operatorname{dist}_p(\cdot,z_0)$,
satisfies $(-\Delta_p+2)b_R\ge\frac34b_R$ for $R\ge1$.
The maximum principle with $Cs^{-2}+Cb_R$, followed by the Poisson estimate on disks of radius one, yields
\begin{equation*}
 |a-g|+|\nabla_p(a-g)|\le Cs^{-2}+Ce^{-c_0s}.
\end{equation*}
The condition $\rho\ge\frac12$ and the
$C^2(B_{\frac t2}(0))$ extension of $|\mathbf{P}_0|^2$ in the $y$ coordinate
justify the same argument at the vertex.
Consequently
\begin{equation*}
 |a|+s|\nabla_pa|\le
 C\left(\tau+\frac\kappa s+\frac1s+se^{-c_0s}\right),
\end{equation*}
and, if $s_n=r_{b,n}(x_n)\to\infty$,
\begin{equation*}
 a_n(x_n)\to0,\qquad s_n|\nabla_{p_n}a_n(x_n)|\to0.
\end{equation*}
Combining this with the limit on fixed disks yields
$\|a_n\|_{L^\infty(\Omega_n)}\to0$.

\emph{Step 3. The weighted gradient estimate.}
It remains to prove $\sup_{x\in\Omega_n}r_{b,n}(x)|V_n(x)|\to0$.
If the phase term has a positive limsup, we choose a subsequence and points
$x_n$ such that $s_n=r_{b,n}(x_n)$ and $s_n|V_n(x_n)|\ge\epsilon>0$.
The convergence on fixed disks implies $s_n\to\infty$.
After rescaling by $s_n$, we obtain a pointed limit $X=\R^2$ if
$\frac{r(x_n)}{s_n}\to\infty$, and $X=C_N$ otherwise. Let
$\mathcal S=\{a\in X:a\text{ is a finite limit of a rescaled }c_{j,n}\}$.
Centers with the same rescaled limit yield one point, and
$1\le|\mathcal S|\le N$, $\operatorname{dist}_X(x_\infty,\mathcal S)=1$,
$\frac{r_{b,n}}{s_n}\to\operatorname{dist}_X(\cdot,\mathcal S)$ locally.
Write $\widehat v_n=s_nV_n$ in the rescaled coordinates. Using the radial cancellation $\nabla_pf_j\cdot b_j=0$, we deduce
\begin{equation*}
 \frac{s^2}{\delta}|\nabla_p\mathcal R\cdot \mathbf{P}_0|
 +\frac{s}{\delta}|\mathcal F|
 \le C\left(\frac\kappa s+\tau\right),\qquad s|\nabla_p\mathcal R|\le Cs^{-2}.
\end{equation*}
The first left side is zero outside the supports. On every compact subset
of $X\setminus\mathcal S$, there are constants $c,C>0$, depending on that
subset and independent of $n$, such that $cs_n\le s\le Cs_n$. Consequently,
\begin{align*}
 \frac{s_n^2}{\delta}\nabla_p\log\rho\cdot \mathbf{P}
 &=\frac{s_n^2}{\rho}\left(
 \nabla_pa\cdot \mathbf{P}+\delta^{-1}\nabla_p\mathcal R\cdot \mathbf{P}_0+
 \nabla_p\mathcal R\cdot V\right),\\
 \left|\frac{s_n^2}{\delta}\nabla_p\log\rho\cdot \mathbf{P}\right|
 &\le C\left(s_n|\nabla_pa|+\frac\kappa{s_n}+\tau+s_n^{-2}\right)\to0.
\end{align*}
Interior Poisson estimates for $\psi_n/\delta_n$, after subtracting
its local mean in the rescaled coordinates, yield locally uniform
convergence of $\widehat v_n$ away
from $\mathcal S\cup\{o\}$ to a harmonic closed field satisfying
$|v_\infty|\le\operatorname{dist}_X(\cdot,\mathcal S)^{-1}$.
If $o\notin\mathcal S$, Lemma~\ref{tippgrad} implies
$\limsup_n\sup_{B_\eta(o)}|\widehat v_n|\le C\eta^{p_--1}$.
Thus $|v_\infty(x_\infty)|\ge\epsilon$.

The existence of a single valued potential implies that all tangential line integrals vanish.
For the normal boundary integrals, we use
\begin{equation*}
 \frac{s_n}{\delta_n}\mathcal J_n
 =\rho_n^2\widehat v_n
   +(\rho_n+\mathcal R_n)a_n(s_n\mathbf{P}_{0,n})
   -\frac{s_n}{\delta_n}\mathcal F_n.
\end{equation*}
The last two terms vanish locally away from $\mathcal S$. For every
positively oriented simple contour $\Gamma$ disjoint from $\mathcal S$
and, when $X=C_N$, from the vertex, the rescaled distance to the
centers is bounded below on $\Gamma$. Its pullback $\Gamma_n$ is therefore disjoint from the holes
for large $n$, and the two limiting boundary integrals below follow from local uniform convergence. With $\nu$ the outward normal,
\[
 \int_\Gamma v_\infty\cdot\tau\,ds_X=0,\qquad
 \int_\Gamma v_\infty\cdot\nu\,ds_X
 =\lim_n\frac1{\delta_n}\int_{\Gamma_n}\mathcal J_n\cdot\nu\,ds_{p_n}=0.
\]
These contours may surround the vertex and points arising from several centers. The
holomorphic one form $(v_{\infty,1}-iv_{\infty,2})\,dz$ is globally
defined and satisfies
\[
 \int_\Gamma(v_{\infty,1}-iv_{\infty,2})\,dz
 =\int_\Gamma v_\infty\cdot\tau\,ds_X
       +i\int_\Gamma v_\infty\cdot\nu\,ds_X=0.
\]
Its at most simple poles away from the vertex therefore have zero
residues. When $X=C_N$, use $w=z^N$ and write
$(v_{\infty,1}-iv_{\infty,2})\,dz=q(w)\,dw$.
At the vertex,
\[
 |q(w)|\le
 \begin{cases}C|w|^{\frac1N-1},&o\notin\mathcal S,\\
 C|w|^{-1},&o\in\mathcal S,\end{cases}
 \qquad\operatorname{Res}_{w=0}q=0.
\]
All finite singularities are removable. Since $\mathcal S$ is finite,
$|v_\infty|=O(r^{-1})$ at infinity. Thus
\[
 |q(w)|=\frac{|v_\infty|}{Nr^{N-1}}
          =O(r^{-N})=O(|w|^{-1})
\]
when $X=C_N$. Liouville's theorem implies $q=0$.
For $X=\R^2$, the choice $w=z$ and the same $O(|w|^{-1})$ bound
yield this conclusion.
This contradicts $|v_\infty(x_\infty)|\ge\epsilon$. Therefore
\[
 1=\max_j\left\|\frac{W_{j,n}}{\delta_n}\right\|_{C^0(B_{2R_0}(0))}
 +\|a_n\|_{L^\infty(\Omega_n)}
 +\sup_{x\in\Omega_n}r_{b,n}(x)|V_n(x)|\to0,
\]
a contradiction. Thus $\delta_n\le CL_n^{-2}$. Reapplying the boundary
integral estimate above and the elliptic estimates on fixed disks about the centers, we have
\begin{gather*}
 A_{*,n}\le\frac C{L_n}(\delta_n+L_n^{-2})\le CL_n^{-3},\\
 \|W_{j,n}\|_{C^k(\overline{B_R(0)})}
 \le C_{R,k}(\delta_n+A_{*,n}+L_n^{-2})\le C_{R,k}L_n^{-2}.
\end{gather*}This finishes the proof.
\end{proof}

On each annulus $R_0\le r_j\le4cL$ we now have
\begin{equation*}
 |\rho-f_j|\le CL^{-2},\qquad |v|\le\frac{CL^{-2}}{r_j},\qquad
 |H_j|\le C(L^{-3}+L^{-2}r_j).
\end{equation*}
The local phase difference, initially $O(L^{-2})$, consequently
satisfies
$|\Psi_j|\le CL^{-2}(1+\log r_j+L^{-1}r_j+r_j^2)$.
Together with the estimate on fixed disks, this implies 
\begin{equation*}
 |W_j(y)|\le CL^{-2}(1+|y|^2),\qquad |y|\le4cL.
\end{equation*}

To proceed, let us set $\varepsilon=L^{-1}$. In a fixed centered Euclidean chart suppress $j$ and write
$d=\pm1$, $b=\frac d r e_\vartheta$, $H=\mathbf{P}_0-b$, and
$W=e^{-i\gamma}U(c_j+y)-V_d(y)$. Then
\begin{equation*}
 H(y)=H(0)+\varepsilon^2M_ny+O(\varepsilon^3|y|^2),\qquad
 H(0)=O(\varepsilon^3),\qquad
 \sigma_n^2=\tfrac12|M_n|_{\mathrm F}^2=L^4|\beta_j|^2.
\end{equation*}
The matrices are bounded, symmetric, and of zero trace.
Angular averages are denoted by $\langle\cdot\rangle$.

On each fixed disk about a center, the angular mean of the linear
modulus error is of order $L^{-4}$.
\begin{lemma}\label{fixehforde}
Let $U_n$, $p_n$, and the modulated centers and phases satisfy the
hypotheses and normalization of Proposition~\ref{unitrates2}.
Fix a center of degree $d=\pm1$ and use the preceding definitions of
$W$, $b$ and $\varepsilon$ in its Euclidean coordinate $y$.
For every fixed $R>0$,
\begin{equation}\label{fixedpote}
 \int_{B_R(0)}\bigl[(1-|V_d+W|^2)^2-(1-f^2)^2\bigr]\,dy
 =O_R(\varepsilon^4).
\end{equation}
For every $0<R_1<R_2<\infty$, the annulus
$\{y\in\R^2:R_1\le|y|\le R_2\}$ is disjoint from the zero set of $V_d+W$
for all sufficiently large $n$. On this annulus, set $\rho=|V_d+W|$, $\mathbf{P}=\nabla\arg(V_d+W)$, and
$\mu=\langle\rho-f\rangle$. Then
\[
 \|\mu\|_{C^1([R_1,R_2])}
 +\|\langle|\mathbf{P}|^2-|b|^2\rangle\|_{C^0([R_1,R_2])}
 \le C_{R_1,R_2}\varepsilon^4.
\]
Constants may depend on the sequence, fixed radii, and modulation.
They are independent of the largest center radius and the configuration diameter.
\end{lemma}
\begin{proof}
Write $e^{-id\vartheta}W=A+iB$ and $m=\langle A\rangle$.
The angular derivatives disappear on averaging the exact difference equation:
\begin{equation*}
 \mathcal M m=-f\langle3A^2+B^2\rangle
                    -\langle A(A^2+B^2)\rangle,\qquad
 \mathcal M=-\partial_r^2-r^{-1}\partial_r+r^{-2}+3f^2-1.
\end{equation*}
We set $T=\frac c\varepsilon$. From the bound
$|W|\le C\varepsilon^2(1+r^2)$, we obtain
\[
 |\mathcal Mm|\le C\varepsilon^4f(1+r^4),\qquad
 |m(T)|\le C\varepsilon^2(1+T^2).
\]
At zero the term of degree three obeys the same bound because its constant
projection onto frequency $d$ vanishes and
$|\langle A|W|^2\rangle|\le r\|\nabla(W|W|^2)\|_{L^\infty(B_1(0))}
\le C\varepsilon^6r$. By smoothness, $m/r$ and $f/r$ extend smoothly as functions of $r^2$
at zero. Since $f/r$ has a positive limit, the quotient $m/f$ is regular there. The radial equation yields
\begin{equation}\label{annuaconj}
 \mathcal M(fh)
 =f\left[-\frac{(rf^2h')'}{rf^2}+2f^2h\right].
\end{equation}
For $\Phi=f(C_0-r^2+r^4)$, we compute from \eqref{annuaconj} that
\[
 \frac{\mathcal M\Phi}{f}
 =4-16r^2+4r(1-2r^2)\frac{f'}f
       +2f^2(C_0-r^2+r^4).
\]
The ratio $\mathcal M\Phi/[f(1+r^4)]$ tends to $8$ at zero and $2$
at infinity. Starting with $C_0\ge1$, increasing $C_0$ increases this
ratio by a positive multiple of $2f^2/(1+r^4)$. Thus one fixed $C_0$
makes $\Phi>0$ and $\mathcal M\Phi\ge f(1+r^4)$ on $(0,\infty)$,
by treating the remaining compact interval. Since $T=c/\varepsilon$
and $f(T)\to1$,
\[
 |m(T)|\le C\varepsilon^2(1+T^2)
           \le C_c\varepsilon^4\Phi(T)
\]
for all sufficiently small $\varepsilon$. Indeed, $\varepsilon^2T^2=c^2$,
and $C^{-1}T^4\le\Phi(T)\le CT^4$ for all sufficiently large $T$, with
$C$ independent of $\varepsilon$.
Consequently, a constant independent of $\varepsilon$ can be chosen so that
$\mathcal M(C\varepsilon^4\Phi\pm m)\ge0$ and the data at $T$ are nonnegative.
For $h=\frac{C\varepsilon^4\Phi\pm m}{f}$, we test against
$h_-=\max\{-h,0\}$ and obtain
\[
 0\ge\int_0^T\bigl(rf^2|(h_-)'|^2+2rf^4h_-^2\bigr)\,dr\ge0.
\]
The boundary term at zero vanishes because $rf^2=O(r^3)$ and $h$ is regular.
Thus
\begin{equation*}
 |m(r)|\le C\varepsilon^4\Phi(r)
          \le C\varepsilon^4f(r)(1+r^4),\qquad 0\le r\le T.
\end{equation*}
The interior estimate in one dimension also implies
$\|m\|_{C^1([R_1,R_2])}\le C_{R_1,R_2}\varepsilon^4$.
Now use the exact identity
\begin{equation*}
 (1-|V_d+W|^2)^2-(1-f^2)^2
 =-4(1-f^2)fA-2(1-f^2)|W|^2
   +4f^2A^2+4fA|W|^2+|W|^4.
\end{equation*}
The linear term has mean $m=O_R(\varepsilon^4)$, and every remaining
term is at least quadratic in $W=O_R(\varepsilon^2)$.
This proves \eqref{fixedpote}.
On fixed annuli $f+A\ge\frac12f(R_1)>0$, and
\begin{gather*}
 \rho-f=A+\frac{B^2}{\rho+f+A},\quad
 \|\mu-m\|_{C^1([R_1,R_2])}\le C_{R_1,R_2}\varepsilon^4,\\
 \langle|\mathbf{P}|^2-|b|^2\rangle=\langle|\nabla\psi|^2\rangle=O(\varepsilon^4),
\end{gather*}
where $\psi=\arctan\frac B{f+A}$ is periodic,
$\nabla\psi=O(\varepsilon^2)$, and $\langle b\cdot \nabla\psi\rangle=0$.
\end{proof}

We next estimate the modulus and its derivatives on the expanding annuli
about the centers.
We use the sequence and modulation of Proposition~\ref{unitrates2}
and Lemma~\ref{fixehforde}. Thus $p_n\to N$, $N\ge3$ is
odd, $\varepsilon_n=L_n^{-1}$, $t_n\ge c_*L_n$, and
$A_{*,n}=O(L_n^{-3})$. Fix a sufficiently small $c>0$ so that the disks
$B_{16cL_n}(c_{j,n})$ are disjoint Euclidean disks. In one such disk write
\[
 a=\rho-f,\qquad b=\frac d r e_\vartheta,\qquad H=\mathbf{P}_0-b,
 \qquad v=\mathbf{P}-\mathbf{P}_0,\qquad \mu=\langle a\rangle,
 \qquad h=\langle|H|^2\rangle,
\]
where $d=\pm1$ and
$\langle g\rangle=\frac1{2\pi}\int_0^{2\pi}g\,d\vartheta$.
We differentiate in these Euclidean coordinates. Constants may depend on
fixed $R,c$ and the sequence. They are independent of its largest center radius and diameter.

\begin{lemma}\label{annuampl}
For the sequence, modulation and constant $c$ fixed above, use the
preceding notation in each centered Euclidean chart and set
$\delta_H=2b\cdot H+|H|^2$ and $g=-\frac12f\delta_H$.
There is $R_*>0$ such that, for every fixed $R\ge R_*$ and all
sufficiently large $n$, the following estimates hold on
$\{y\in\R^2:R\le|y|\le cL_n\}$:
\begin{gather}
 |a-g|+|\nabla(a-g)|+|\Delta(a-g)|\le C_R\varepsilon^2r^{-2},\notag
\\
 \left|\mu+\frac{fh}{2}\right|
 +\left|\mu'+\left(\frac{fh}{2}\right)'\right|
 \le C_R\varepsilon^4,
 \label{annumappr}\\
 |\mu|\le C_R\varepsilon^4(1+r^2),\quad
 |\mu'|\le C_R\varepsilon^4(1+r),\quad
 |\langle\Delta a\rangle|\le C_R\varepsilon^4,
 \label{annumboun}\\
 |\nabla a|\le C_R\varepsilon^2r^{-1},\qquad
 |\Delta a|\le C_R\varepsilon^2r^{-2}.
 \label{annuaderi}
\end{gather}
\end{lemma}
\begin{proof}
On the larger annulus,
$|a|+r|v|+r^{-1}|H|+|\nabla H|\le C\varepsilon^2$ and
$|\nabla^2H|\le C\varepsilon^3$.
The vanishing periods imply $\langle b\cdot H\rangle=\langle b\cdot v\rangle=0$.
Subtracting and averaging the equations for the modulus yields
\begin{align*}
 \mathcal M\mu
 &=-f\langle|H+v|^2\rangle
 -\langle a[2b\cdot(H+v)+|H+v|^2]\rangle\\
 &\quad-3f\langle a^2\rangle-\langle a^3\rangle
 =-fh+O(\varepsilon^4),\\
 \mathcal M&=-\partial_r^2-r^{-1}\partial_r+q(r),
 \qquad q=r^{-2}+3f^2-1=2+O(r^{-2}).
\end{align*}
The error is uniform since
$|\mathcal M\mu+fh|\le C(\varepsilon^4+\varepsilon^4r^{-2}
+\varepsilon^6r^2)\le C\varepsilon^4$.
Harmonicity of $H$ and the radial equation yield
\begin{gather*}
 h\le C\varepsilon^4r^2,\quad |h'|\le C\varepsilon^4r,\\
 h''+\frac{h'}r=2\langle|\nabla H|_{\mathrm F}^2\rangle\le C\varepsilon^4,
 \\
 \mathcal M\left(-\frac{fh}2\right)+fh
 =\frac f2\left(h''+\frac{h'}r\right)+f'h'+f(1-f^2)h=O(\varepsilon^4).
\end{gather*}
For $e=\mu+\frac{fh}2$, Lemma~\ref{fixehforde} implies
$|e(R-1)|\le C_R\varepsilon^4$, while
$|e(2cL)|\le C\varepsilon^2$ and $|\mathcal Me|\le C\varepsilon^4$.
We choose $R$ so that $q\ge1$ on $[R-1,2cL]$. The maximum principle,
applied with a constant multiple of $\varepsilon^4$ and the supersolutions
$e^{-\frac{r-R+1}2}$, $e^{-\frac{2cL-r}2}$, yields
\[
 |e(r)|\le C_R\varepsilon^4(1+e^{-\frac{r-R+1}2})
 +C\varepsilon^2e^{-\frac{2cL-r}2}.
\]
The interior estimate in one dimension yields
$|e|+|e'|\le C_R\varepsilon^4$ for $R\le r\le cL$.
This proves \eqref{annumappr}. Moreover,
$\Delta\mu=q\mu+fh+O(\varepsilon^4)
=fh(1-\frac q2)+O_R(\varepsilon^4)=O_R(\varepsilon^4)$,
which proves \eqref{annumboun}.

Exact subtraction of the equations for $\rho$ and $f$ yields
\begin{gather*}
 \begin{aligned}
 (-\Delta+2)(a-g)
 &=(3-|b|^2-\rho^2-\rho f-f^2)a-a\delta_H\\
 &\quad-\rho\bigl(2(b+H)\cdot v+|v|^2\bigr)-\tfrac12\Delta(f\delta_H),
 \end{aligned}
 \\
 \Delta \delta_H=4\nabla b:\nabla H+2|\nabla H|_{\mathrm F}^2.
\end{gather*}
Here $\bigl|3-|b|^2-\rho^2-\rho f-f^2\bigr|\le Cr^{-2}$,
$|\delta_H|+r|\nabla \delta_H|\le C\varepsilon^2$, and
$\bigl|2(b+H)\cdot v+|v|^2\bigr|+|\Delta \delta_H|\le C\varepsilon^2r^{-2}$,
so the right side is
$O(\varepsilon^2r^{-2})$ without differentiating $v$.
For $|y|=r\ge2R$, the maximum principle on $B_{\frac r4}(y)$ with
\[
 \frac{C_R\varepsilon^2}{r^2}
 +C\varepsilon^2\exp\left(\frac{|x-y|^2-\frac{r^2}{16}}r\right)
\]
applies because $|a-g|\le C\varepsilon^2$ on its boundary. It yields
$|a-g|(y)\le C_R\varepsilon^2r^{-2}$.
The equation and the Poisson estimate on disks of radius one yield the same bound for
$|\nabla(a-g)|+|\Delta(a-g)|$. The estimates on fixed disks handle $R\le r\le2R$.
Finally $|\nabla g|\le C\varepsilon^2r^{-1}$ and
$|\Delta g|\le C\varepsilon^2r^{-2}$ prove
\eqref{annuaderi}.
\end{proof}

Define the nonnegative function
\begin{equation*}
 \mathcal I(\rho)=2(1+\rho^{-2})|\nabla_p\rho|^2
                         +\left(\frac{\Delta_p\rho}\rho\right)^2.
\end{equation*}
For functions on Euclidean domains, the same notation means the
expression above with Euclidean derivatives.
Let $E\Subset\{x\in C_p:\rho(x)>0\}$ be a bounded piecewise $C^1$
domain, and let $\nu$ be its outward unit normal for the cone metric.
We write $\partial_\nu\rho=\langle\nabla_p\rho,\nu\rangle_p$.
With $w=1-\rho^2$,
we infer from the equation for the modulus that
\[
 w=|\mathbf{P}|^2-\frac{\Delta_p\rho}\rho,\qquad
 w^2=|\mathbf{P}|^4-2\frac w\rho\Delta_p\rho
              -\left(\frac{\Delta_p\rho}\rho\right)^2,\qquad
 \nabla_p\left(\frac w\rho\right)=-(1+\rho^{-2})\nabla_p\rho.
\]
Integration by parts therefore yields
\begin{equation}\label{exacpiden}
 \int_E(1-\rho^2)^2\,dA_p
 =\int_E\bigl(|\mathbf{P}|^4-\mathcal I(\rho)\bigr)\,dA_p
       -2\int_{\partial E}\frac{1-\rho^2}{\rho}
                                      \partial_\nu\rho\,ds_p.
\end{equation}
If $E$ contains the vertex, apply the identity on
$E\setminus\overline{B_\eta(o)}$ and let $\eta\to0$.
The   vertex estimates make the additional boundary term vanish.
The resulting expression involves only the phase gradient itself.
The next estimate compares the nonnegative correction integral outside
the fixed disks with the corresponding sum for isolated unit vortices.

\begin{lemma}\label{modudefe}
Under the assumptions and normalization of Proposition~\ref{unitrates2},
use its modulated centers $c_{j,n}$, separation $L_n$ and modulus
$\rho_n=|U_n|$. For every sufficiently large fixed $R$,
set $\Omega_{R,n}=C_{p_n}\setminus\bigcup_{j=1}^N\overline{B_R(c_{j,n})}$.
Then, as $n\to\infty$,
\begin{equation}\label{modudefe2}
 \int_{\Omega_{R,n}}\mathcal I(\rho_n)\,dA_{p_n}
 -N\int_{\R^2\setminus\overline{B_R(0)}}\mathcal I(f)\,dy
 =O_{R,c}(L_n^{-4}).
\end{equation}
\end{lemma} 

\begin{proof}
We suppress $n$ and use $\varepsilon=L^{-1}$ and the notation of
Lemma~\ref{annuampl}. In Euclidean coordinates $y$ about a center
$c_j$, let $r=|y|$ and $\mathcal A=\{R<r<cL\}$. In particular,
$a=\rho-f$, $\delta_H=|\mathbf{P}_0|^2-r^{-2}$, and $g=-f\delta_H/2$.

\emph{Step 1. The integral of $|\nabla\rho|^2-f'^2$ on $\mathcal A$.}
Lemma~\ref{annuampl} and the estimates for $H$ imply
\begin{gather*}
 |\delta_H|+r|\nabla \delta_H|\le C\varepsilon^2,\qquad
 |a-g|+|\nabla(a-g)|\le C_R\varepsilon^2r^{-2},\\
 |\langle \delta_H\rangle|\le C\varepsilon^4r^2,\qquad
 |\langle\partial_r\delta_H\rangle|\le C\varepsilon^4r,\qquad
 |\langle\partial_r(a-g)\rangle|\le C_R\varepsilon^4.
\end{gather*}
We differentiate $\rho=f-f\delta_H/2+(a-g)$ and expand the square:
\begin{align*}
 |\nabla\rho|^2-f'^2
 ={}&-ff'\partial_r\delta_H+\frac{f^2}{4}|\nabla \delta_H|^2\\
 &-f'^2\delta_H+2f'\partial_r(a-g)
       +O_R(\varepsilon^4r^{-3}).
\end{align*}
The terms $-f'^2\delta_H$ and $2f'\partial_r(a-g)$ have integrals of
order $O_R(\varepsilon^4)$ after angular averaging. The same holds
for the remainder. Using $ff'-r^{-3}=O(r^{-4})$,
$f^2-1=O(r^{-2})$, and the preceding mean estimates, we obtain
\begin{align*}
 \int_{\mathcal A}(|\nabla\rho|^2-f'^2)\,dy
 &=\int_{\mathcal A}
       \left[-r^{-3}\partial_r\delta_H+\frac14|\nabla \delta_H|^2\right]\,dy
       +O_R(\varepsilon^4)\\
 &=\frac14\int_{\mathcal A}
       \left[|\nabla(|\mathbf{P}_0|^2)|^2-|\nabla(r^{-2})|^2\right]\,dy
       +O_R(\varepsilon^4)\\
 &=\frac1{16}\int_{\partial\mathcal A}
       \partial_\nu(|\mathbf{P}_0|^4-r^{-4})\,ds+O_R(\varepsilon^4).
\end{align*}
The last equality follows from $|\mathbf{P}_0|=|F|$ and
$\Delta|F|^4=4|\nabla(|F|^2)|^2$ for holomorphic $F$.
The outer boundary integral is $O_c(L^{-4})$, since
$F=O(L^{-1})$ and $F'=O(L^{-2})$ there. At the inner boundary,
$F^2-y^{-2}-2d_jA_jy^{-1}$ is holomorphic on the disk of radius
$4cL$ and bounded by $CL^{-2}$. Cauchy's estimates and orthogonality
of the angular powers show that, at $r=R$,
\[
 \partial_r\langle|F|^4-r^{-4}\rangle
 =-8|A_j|^2R^{-3}+O_R(L^{-6})=O_R(L^{-6}).
\]
Here the constant Taylor coefficient contributes a constant to the
angular average, whose radial derivative is zero. Consequently,
$\int_{\mathcal A}(|\nabla\rho|^2-f'^2)\,dy=O_{R,c}(\varepsilon^4)$.

\emph{Step 2. The integral of $\mathcal I(\rho)-\mathcal I(f)$ on $\mathcal A$.}
By \eqref{annumboun}--\eqref{annuaderi} and
$|\nabla\rho|^2-f'^2=2f'a_r+|\nabla a|^2$, we have
\begin{gather*}
 \bigl||\nabla\rho|^2-f'^2\bigr|
       \le C_R(\varepsilon^2r^{-4}+\varepsilon^4r^{-2}),\qquad
 \bigl|\langle|\nabla\rho|^2-f'^2\rangle\bigr|
       \le C_R\varepsilon^4r^{-2},\\
 |\rho^{-2}-f^{-2}|\le C\varepsilon^2,\qquad
 |\langle\rho^{-2}-f^{-2}\rangle|
       \le C_R\varepsilon^4(1+r^2).
\end{gather*}
We use $1+\rho^{-2}=2+(f^{-2}-1)+(\rho^{-2}-f^{-2})$.
The terms with radial factors $f^{-2}-1$ and $f'^2$ are estimated
after angular averaging. The product
$(\rho^{-2}-f^{-2})(|\nabla\rho|^2-f'^2)$ has absolute integral at most
$C_R(\varepsilon^4+\varepsilon^6\log L)=O_R(\varepsilon^4)$.
It follows that
\begin{align*}
 &\int_{\mathcal A}
   \left[2(1+\rho^{-2})|\nabla\rho|^2-2(1+f^{-2})f'^2\right]\,dy\\
 &\qquad=4\int_{\mathcal A}(|\nabla\rho|^2-f'^2)\,dy
       +O_R(\varepsilon^4)=O_{R,c}(\varepsilon^4).
\end{align*}
For the squared Laplacians, $\Delta\rho=\Delta f+\Delta a$ yields
\[
 \left(\frac{\Delta\rho}{\rho}\right)^2
 -\left(\frac{\Delta f}{f}\right)^2
 =\frac{2\Delta f\Delta a+(\Delta a)^2}{f^2}
       +(\rho^{-2}-f^{-2})(\Delta f+\Delta a)^2.
\]
The term linear in $\Delta a$ is estimated using
$|\langle\Delta a\rangle|\le C_R\varepsilon^4$.
For $(\rho^{-2}-f^{-2})(\Delta f)^2$, we use the mean estimate
above. All other products are estimated in absolute value using
$\Delta f=O(r^{-4})$ and $|\Delta a|\le C_R\varepsilon^2r^{-2}$.
Thus
\[
 \left|\int_{\mathcal A}
 \left[\left(\frac{\Delta\rho}{\rho}\right)^2
       -\left(\frac{\Delta f}{f}\right)^2\right]\,dy\right|
 \le C_R\varepsilon^4\int_R^{cL}(r^{-3}+r^{-5}+r^{-7})\,dr
 =O_R(\varepsilon^4).
\]
The two estimates prove
$\int_{\mathcal A}(\mathcal I(\rho)-\mathcal I(f))\,dy
=O_{R,c}(\varepsilon^4)$.

\emph{Step 3. The region away from all centers.}
For $x$ with $r_b(x)\ge cL$, we choose a ball disjoint from the
centers, either Euclidean or centered at the vertex. Its radius is
comparable to $r_b(x)$, it contains $x$ in its concentric ball of
half the radius, and $r_b$ is comparable to $r_b(x)$ throughout it.
The constants in these geometric bounds depend only on $K$.
On this ball, the exact equation is
\begin{align*}
 &(-\Delta_p+\rho(1+\rho))
       \left(1-\rho-\frac12|\mathbf{P}_0|^2\right)\\
 &\qquad=\rho(|\mathbf{P}|^2-|\mathbf{P}_0|^2)
       +\frac12\rho(1-\rho)|\mathbf{P}_0|^2+\frac12\Delta_p|\mathbf{P}_0|^2
       =O\bigl(\varepsilon^2r_b(x)^{-2}+r_b(x)^{-4}\bigr).
\end{align*}
Also, $\bigl|1-\rho-\frac12|\mathbf{P}_0|^2\bigr|\le Cr_b(x)^{-2}$ on this ball,
and $\rho(1+\rho)\ge\frac34$. The maximum principle and Poisson
estimates on disks of radius one, with Lemma~\ref{tippgrad} at the
vertex, imply
\[
 |\nabla_p\rho(x)|\le C\bigl(r_b(x)^{-3}+\varepsilon^2r_b(x)^{-2}\bigr),
 \qquad
 |\Delta_p\rho(x)|\le C\bigl(r_b(x)^{-4}+\varepsilon^2r_b(x)^{-2}\bigr).
\]
For $s>0$, the inclusion
$\{s\le r_b<2s\}\subset\bigcup_{j=1}^NB_{2s}(c_j)$ bounds
the area by $C_KNs^2$. Summing over $s=2^mcL$ yields
\begin{align*}
 \int_{\{r_b\ge cL\}}\mathcal I(\rho)\,dA_p
 &\le C\sum_{m=0}^\infty
       \left[(2^mcL)^{-4}+\varepsilon^4(2^mcL)^{-2}\right]
       \le C_cL^{-4},\\
 \int_{\R^2\setminus\overline{B_{cL}(0)}}\mathcal I(f)\,dy
 &\le C\int_{cL}^\infty(r^{-6}+r^{-8})r\,dr\le C_cL^{-4}.
\end{align*}
Together with the estimates on the $N$ annuli, these bounds prove
\eqref{modudefe2} uniformly in the largest center radius and diameter.
\end{proof}

For an isolated unit vortex $U=f(r)e^{id\vartheta}$ in Euclidean
polar coordinates, with $d=\pm1$, we have
$\nabla\arg U=d e_\vartheta/r$. Hence the integral of $|\nabla\arg U|^4$
outside a disk of radius $R$ equals $\pi/R^2$.
The next proposition shows that the corresponding integral for
$\mathbf{P}_{0,n}$ differs from the sum of these $N$ values by at most
$C_RL_n^{-4}$, uniformly in the configuration diameter. Together with
Lemma~\ref{modudefe}, this estimate prepares the subtraction of the
isolated vortex contributions in Proposition~\ref{poteeesti}.

\begin{proposition}\label{conefpart}
Let $p_n,U_n,c_{j,n},L_n$ satisfy the hypotheses of
Lemma~\ref{modudefe}, and let $\mathbf{P}_{0,n}$ be the explicit phase gradient
\eqref{poinfdefi} for these modulated centers and prescribed signs.
For each fixed $R>0$, with
$\Omega_{R,n}=C_{p_n}\setminus\bigcup_{j=1}^N\overline{B_R(c_{j,n})}$,
\begin{equation*}
 \int_{\Omega_{R,n}}|\mathbf{P}_{0,n}|^4\,dA_{p_n}
 -\frac{N\pi}{R^2}=O_R(L_n^{-4})\qquad(n\to\infty).
\end{equation*}
The constant is independent of the largest center radius and the
configuration diameter.
\end{proposition}
\begin{proof}  On the Euclidean disk of radius $2\ell$ at each center,
where $\ell=4cL$, the identity $F^2=-T_1'+B$ and Cauchy's estimate yield
\begin{gather*}
 T_1=y^{-1}+C_j+\sum_{k\ge1}t_{j,k}y^k,
 \quad |C_j|\le CL^{-1},\quad |t_{j,k}|\le CL^{-1}\ell^{-k},\\
 B=q_jy^{-1}+\sum_{k\ge0}b_{j,k}y^k,
 \quad q_j=2d_jA_j,\quad |b_{j,k}|\le CA_*L^{-1}\ell^{-k}.
\end{gather*}
The double pole of $B$ cancels since $d_j^2=1$. The scalar $|T_1|^2$
and the one form $T_1\overline B\,d\bar z$ are globally defined.
At the vertex $T_1=O_n(r^{p-1})$, $T_1',B=O_n(r^{p-2})$. At infinity
$T_1=O_n(r^{-1})$, $T_1',B=O_n(r^{-2})$.
Their boundary terms are $O_n(r^{2p-2})$ and $O_n(r^{-2})$ and vanish
for fixed $n$. Green's and Stokes' formulas, with clockwise holes, yield
\begin{align*}
 \int_{\Omega_R}|T_1'|^2\,dA_p
 &=\frac14\int_{\partial\Omega_R}\partial_\nu|T_1|^2\,ds_p
 =\frac{N\pi}{R^2}-\pi\sum_{j=1}^N\sum_{k\ge1}k|t_{j,k}|^2R^{2k},\\
 \int_{\Omega_R}(-T_1')\overline B\,dA_p
 &=\frac1{2i}\int_{\partial\Omega_R}T_1\overline B\,d\bar z=\pi\sum_{j=1}^N\left(C_j\overline{q_j}
       +\sum_{k\ge1}t_{j,k}\overline{b_{j,k-1}}R^{2k}\right).
\end{align*}
The first error is $O_R(L^{-4})$ and the mixed integral is
$O(\frac{A_*}{L})+O_R(A_*L^{-3})=O_R(L^{-4})$, since $A_*\le CL^{-3}$.
By Lemma~\ref{pointBdual}, we have
$\int_{\Omega_R}|B|^2\,dA_p=O_R(L^{-6}\log L)=o(L^{-4})$.
Expanding $|F^2|^2=|-T_1'+B|^2$ proves the assertion uniformly in the
configuration diameter.
\end{proof}

Each isolated vortex of degree $1$ or $-1$ has potential energy $2\pi$.
We now estimate the difference between $\mathcal P_{p_n}(U_n)$ and the sum $2\pi N$
for the $N$ vortices. Let the sequence and centers be as in
Proposition~\ref{unitrates2}, and let $L_n$ be their minimum separation.
We fix a sufficiently large $R$ and work on
$\Omega_{R,n}=C_{p_n}\setminus\bigcup_{j=1}^N\overline{B_R(c_{j,n})}$,
where $\rho_n=|U_n|\ge\frac12$ for large $n$.
Locally, $U_n=\rho_ne^{i\theta_n}$, so $\mathbf{P}_n=\nabla_{p_n}\theta_n$ is the
gradient of its angle. The explicit field $\mathbf{P}_{0,n}$ in
\eqref{poinfdefi} is the gradient of
$\sum_{j=1}^N d_j\arg(w-w(c_{j,n}))$, with $w=r^{p_n}e^{i\phi}$.
These gradients are independent of the local choices of angles.
We set
\begin{equation}\label{regucapint}
 E_n=\int_{\Omega_{R,n}}\rho_n^2|\mathbf{P}_n-\mathbf{P}_{0,n}|^2\,dA_{p_n},
 \qquad S_n=L_n^4\sum_{j=1}^N|\beta_{j,n}|^2.
\end{equation}
Thus $E_n$ integrates the squared difference between the two gradients,
with weight $\rho_n^2$. In $S_n$, $\beta_{j,n}$ is the coefficient of $y$ after the singular term
$d_j/y$ in \eqref{poinlexpa}. Thus $S_n$ measures the size of these
linear coefficients at the scale $L_n$. The next proposition relates
the potential energy difference to $E_n$ and bounds $E_n$ from below
in terms of $S_n$. Its estimates hold independently of any positive
lower bound for $S_n$.

\begin{proposition}\label{poteeesti}
Under the preceding assumptions, $S_n=O(1)$ and
\begin{align}
 \mathcal P_{p_n}(U_n)-2\pi N&=4E_n+O_R(L_n^{-4}),
 \label{poteeiden2}\\
 E_n&\ge9\pi S_nL_n^{-4}\log L_n-C_RL_n^{-4}.\notag
\end{align}
Consequently,
\[
 \mathcal P_{p_n}(U_n)-2\pi N
 \ge36\pi S_nL_n^{-4}\log L_n-C_RL_n^{-4}.
\]
The constants are independent of the largest center radius and the
configuration diameter.
\end{proposition}
\begin{proof}
By \eqref{unitrfiel},
$S_n=\frac12\sum_{j=1}^N|M_{j,n}|_{\mathrm F}^2\le C$.
We fix $R$ before taking $n\to\infty$, suppress $n$, and set
$\varepsilon=L^{-1}$ and $v=\mathbf{P}-\mathbf{P}_0$.

\emph{Step 1. Subtracting the identities for isolated vortices.}
We apply \eqref{exacpiden} to $\Omega_R$, first cutting it off at radius
$T$ and deleting a disk of radius $s$ about the vertex. The corresponding
boundary terms are $O_n(T^{-3})$ and $O_n(s^2)$, respectively, and vanish
for each fixed solution as $T\to\infty$ and $s\to0$.
The preceding estimates ensure convergence of the volume integrals.
On each circle $\partial B_R(c_j)$, let $a=\rho-f$ and let brackets
denote the angular mean. Lemma~\ref{fixehforde} implies
\[
 |a|+|a_r|=O_R(\varepsilon^2),\qquad
 |\langle a\rangle|+|\langle a_r\rangle|=O_R(\varepsilon^4).
\]
Since $f$ is radial, the identity
\[
 \frac{1-\rho^2}{\rho}\rho_r-\frac{1-f^2}{f}f'
 =(f^{-1}-f)a_r-(1+f^{-2})f'a+O_R(a^2+|aa_r|)
\]
shows that the integrated boundary difference is $O_R(\varepsilon^4)$.
Both exterior domains have normal $-e_r$ on these circles.
The interior potential integrals differ by $O_R(\varepsilon^4)$ by
\eqref{fixedpote}. Proposition~\ref{conefpart} and
Lemma~\ref{modudefe} imply that the differences involving $|\mathbf{P}_0|^4$ and
$\mathcal I(\rho)$ are also $O_R(\varepsilon^4)$.
We therefore subtract the $N$ identities for isolated vortices and use
$\mathcal P(V_1)=2\pi$ and
$\int_{\R^2\setminus\overline{B_R(0)}}|b|^4\,dy=\pi/R^2$ to obtain
\[
 \mathcal P_p(U)-2\pi N
 =\int_{\Omega_R}(|\mathbf{P}|^4-|\mathbf{P}_0|^4)\,dA_p+O_R(\varepsilon^4).
\]

\emph{Step 2. Relating the fourth powers to $E$.}
All zeros lie in the disks $B_R(c_j)$, with degrees $d_j$.
Thus $v=\nabla_p\psi$ for a single valued real function $\psi$ on $\Omega_R$,
and the integral of $v$ along each boundary circle is zero.
The bounds $|\mathbf{P}_0|\le Cr_b^{-1}$ and $|v|\le C\varepsilon^2r_b^{-1}$ imply
\[
 \int_{\Omega_R}
 \bigl||\mathbf{P}|^4-|\mathbf{P}_0|^4-4|\mathbf{P}_0|^2\mathbf{P}_0\cdot v\bigr|\,dA_p
 \le C\varepsilon^4\int_{\Omega_R}r_b^{-4}\,dA_p
 \le C_R\varepsilon^4.
\]
The equation for $\rho$ reads
\[
 |\mathbf{P}_0|^2-(1-\rho^2)
 =\frac{\Delta_p\rho}{\rho}-2\mathbf{P}_0\cdot v-|v|^2.
\]
After multiplication by $\mathbf{P}_0\cdot v$, the last two terms have integrals
of order $O_R(\varepsilon^4)$. To estimate the first term, we choose
a fixed small $c>0$ so that the annuli
$\mathcal A_j=B_{cL}(c_j)\setminus\overline{B_R(c_j)}$ are disjoint
Euclidean annuli. On each annulus, set $b=d_j e_\vartheta/r$ and
$H=\mathbf{P}_0-b$. Lemma~\ref{annuampl} implies
\[
 \frac{\Delta\rho}{\rho}-\frac{\Delta f}{f}
 =O_R(\varepsilon^2r^{-2}),\qquad
 |\mathbf{P}_0\cdot v|\le C\varepsilon^2r^{-2},\qquad
 |H\cdot v|\le C\varepsilon^4.
\]
The radial term cancels because $\Delta f/f$ depends only on $r$ and
$\int_{\partial B_r(c_j)}v\cdot\tau\,ds_p=0$. Since
$\Delta f=O(r^{-4})$, we obtain
\[
 \int_{\mathcal A_j}\frac{\Delta\rho}{\rho}\mathbf{P}_0\cdot v\,dA_p
 =\int_{\mathcal A_j}\frac{\Delta f}{f}b\cdot v\,dA_p
       +O_R(\varepsilon^4)=O_R(\varepsilon^4).
\]
On $\{r_b\ge cL\}$, the absolute value of the integral is at most
$C\varepsilon^2\int_{r_b\ge cL}r_b^{-4}\,dA_p
\le C_c\varepsilon^4$. Hence
\[
 \int_{\Omega_R}(|\mathbf{P}|^4-|\mathbf{P}_0|^4)\,dA_p
 =4\int_{\Omega_R}(1-\rho^2)\mathbf{P}_0\cdot v\,dA_p+O_R(\varepsilon^4).
\]
We now apply \eqref{unitpiden} with the current centers.
The estimate in the proof of Proposition~\ref{unitrates2} ensures
$v\in L^2(\Omega_R)$ after the centers move. Also,
$(1-\rho^2)\mathbf{P}_0\in L^2(\Omega_R)$, both $j(U)$ and $\mathbf{P}_0$ have zero
normal integral on each circle, and the vertex estimates apply since
$U$ is nonzero there. These are precisely the conditions used to
justify the integration by parts in \eqref{unitpiden}.
The field $Z=\mathbf{P}_0-j(U)=(1-\rho^2)\mathbf{P}_0-\rho^2v$ has zero divergence,
and that identity becomes
\[
 \int_{\Omega_R}(1-\rho^2)\mathbf{P}_0\cdot v\,dA_p-E
 =\sum_{j=1}^N\int_{\partial B_R(c_j)}
       (\psi-\langle\psi\rangle_j)Z\cdot\nu_{\Omega_R}\,ds_p.
\]
Here $\langle\psi\rangle_j$ is the mean on the $j$th circle.
On each circle, $\psi-\langle\psi\rangle_j=O_R(\varepsilon^2)$ and
$Z\cdot\nu_{\Omega_R}=((1-\rho^2)H-\rho^2v)\cdot\nu_{\Omega_R}
=O_R(\varepsilon^2)$, because $b$ is tangent to the circle.
The boundary integral is therefore $O_R(\varepsilon^4)$.
Together with Step 1, this proves \eqref{poteeiden2}.

\emph{Step 3. The logarithmic lower bound.}
The linear term $\varepsilon^2M_jy$ in $H$ determines the estimate on
$\mathcal A_j$. We set $s_j=e_r\cdot M_je_r$ and
$\ell_j=e_\vartheta\cdot M_je_r$. Since $M_j$ is symmetric with zero
trace and $\sigma_j^2=\frac12|M_j|_{\mathrm F}^2=L^4|\beta_j|^2$,
\[
 \partial_\vartheta s_j=2\ell_j,\qquad
 \int_0^{2\pi}s_j^2\,d\vartheta
 =\int_0^{2\pi}\ell_j^2\,d\vartheta=\pi\sigma_j^2.
\]
The fields $j(U)$, $f^2b$, and $H$ have zero divergence. Thus
\[
 \mathcal J_j:=j(U)-f^2b-H=\rho^2v+\mathcal F_j,\qquad
 \mathcal F_j:=(\rho^2-f^2)b+(\rho^2-1)H,\qquad
 \operatorname{div}\mathcal J_j=0.
\]
By \eqref{unitrfiel}, Lemma~\ref{annuampl}, and $a=\rho-f$, we have
\begin{gather*}
 H=\varepsilon^2r(s_je_r+\ell_je_\vartheta)
       +O(\varepsilon^3+\varepsilon^3r^2),\qquad
 f^2-1=-r^{-2}+O(r^{-4}),\\
 a=-fd_j\varepsilon^2\ell_j+
 O_R(\varepsilon^2r^{-2}+\varepsilon^3r+\varepsilon^4r^2+
 \varepsilon^3r^{-1}).
\end{gather*}
In $\mathcal F_j=(2fa+a^2)b+(f^2-1)H+(2fa+a^2)H$, the first
two terms contribute, at leading order, $-2\varepsilon^2\ell_j e_\vartheta/r$ and
$-\varepsilon^2(s_je_r+\ell_je_\vartheta)/r$, respectively. Therefore
\[
 \mathcal F_j=-\frac{\varepsilon^2}{r}
          (s_je_r+3\ell_je_\vartheta)
       +O_R(\varepsilon^2r^{-3}+\varepsilon^3+\varepsilon^4r).
\]
We choose $\eta\in C_c^\infty((R,cL))$, with $0\le\eta\le1$, equal to
one on $[2R,cL/2]$, and $r|\eta'|\le C$ on the two transition annuli.
Let $\varphi_j=\frac32\varepsilon^2\eta(r)s_j(\vartheta)$ on
$\mathcal A_j$, extended by zero. Then
\[
 \nabla_p\varphi_j=\frac32\varepsilon^2\eta's_je_r
       +3\varepsilon^2\frac\eta r\ell_je_\vartheta.
\]
The logarithm comes from the radial integral of $1/r$. More precisely,
for $k=1,2$ the cutoff satisfies
\[
 \int_R^{cL}\eta'\,dr=0,\qquad
 \int_R^{cL}r|\eta'|^2\,dr\le C,\qquad
 \int_R^{cL}\frac{\eta^k}{r}\,dr=\log L+O_{R,c}(1).
\]
The remainder in $\mathcal F_j$, multiplied by $\nabla_p\varphi_j$,
has integral bounded by
$C_R(\varepsilon^4+\varepsilon^5L+\varepsilon^6L^2)=O_R(\varepsilon^4)$.
Testing $\operatorname{div}\mathcal J_j=0$ against $\varphi_j$ and using
the angular identities above, we obtain
\begin{align*}
 \int_{\mathcal A_j}\rho^2v\cdot \nabla_p\varphi_j\,dA_p
 &=-\int_{\mathcal A_j}\mathcal F_j\cdot \nabla_p\varphi_j\,dA_p
 =9\pi\sigma_j^2\varepsilon^4\log L+O_R(\varepsilon^4),\\
 \int_{\mathcal A_j}\rho^2|\nabla_p\varphi_j|^2\,dA_p
 &\le\int_{\mathcal A_j}|\nabla_p\varphi_j|^2\,dA_p
 =9\pi\sigma_j^2\varepsilon^4\log L+O_R(\varepsilon^4).
\end{align*}
The second line uses $\rho^2\le1$. Since the annuli are disjoint,
$\varphi=\sum_{j=1}^N\varphi_j$ and $S=\sum_{j=1}^N\sigma_j^2$ satisfy
\begin{align*}
 E&\ge2\int_{\Omega_R}\rho^2v\cdot \nabla_p\varphi\,dA_p
       -\int_{\Omega_R}\rho^2|\nabla_p\varphi|^2\,dA_p\\
 &\ge9\pi S\varepsilon^4\log L-C_R\varepsilon^4,
\end{align*}
where the first inequality follows by expanding
$\int_{\Omega_R}\rho^2|v-\nabla_p\varphi|^2\,dA_p\ge0$.
The last conclusion follows from \eqref{poteeiden2}.
\end{proof}

The following identity for the third moment shows that, under the stated
algebraic conditions, at least one linear coefficient of the regular parts
is nonzero.

\begin{lemma}\label{cubimiden}
Let $M\ge2$ be an integer, let $b_1,\ldots,b_M\in\C$ be pairwise
distinct, and let $d_1,\ldots,d_M\in\{-1,1\}$. Suppose that
\[
 \sum_{\substack{1\le k\le M\\k\ne j}}\frac{d_k}{b_j-b_k}=0
 \quad(1\le j\le M),\qquad
 Q:=\sum_{j=1}^M d_j,\qquad
 \beta_j:=-\sum_{\substack{1\le k\le M\\k\ne j}}
                 \frac{d_k}{(b_j-b_k)^2}.
\]
Then
\begin{equation*}
 M=Q^2,\qquad 3\sum_jb_j^2\beta_j=Q^3-Q\ne0.
\end{equation*}
\end{lemma}
\begin{proof}
The identity $M=Q^2$ is proved in \cite[\S3, p.~16]{Aref07}.
For $F(z)=\sum_{j=1}^M\frac{d_j}{z-b_j}$, the assumed identities yield
$F(b_j+y)=\frac{d_j}{y}+\beta_jy+O(y^2)$. Thus
\[
 F(z)^3-\sum_{j=1}^M\frac{d_j}{(z-b_j)^3}
          -3\sum_{j=1}^M\frac{\beta_j}{z-b_j}=0,
\]
since the left side is entire and vanishes at infinity.
Comparing coefficients of $z^{-3}$ yields
$Q^3=Q+3\sum_{j=1}^Mb_j^2\beta_j$.
Finally $Q\in\mathbb Z$ and $Q^2=M\ge2$, so $Q^3-Q\ne0$.
\end{proof}

By applying the moment identity at the smallest cluster scale, we obtain a
strictly positive lower bound for the sum of the squared rescaled linear
coefficients. This is the content of the following

\begin{proposition}\label{nonvlcoef2}
Under the assumptions and normalization of Proposition~\ref{unitrates2},
let $S_n$ be defined by \eqref{regucapint}. Then
\[
 \liminf_{n\to\infty}L_n^4\sum_{j=1}^N|\beta_{j,n}|^2
 =\liminf_{n\to\infty}S_n>0.
\]
\end{proposition} 
\begin{proof}
Suppose, to the contrary, that $S_n\to0$ along a subsequence.
We will construct a rational function with at least two distinct poles,
all with residues $\pm1$, whose constant and linear regular coefficients
vanish at every pole. This contradicts Lemma~\ref{cubimiden}.

Relabel the centers so that
$\operatorname{dist}_{p_n}(c_{1,n},c_{2,n})=L_n$.
The separation estimates and Proposition~\ref{unitrates2} imply
\[
 t_n\ge c_*L_n,\quad L_n\to\infty,\quad
 \max_jL_n|A_{j,n}|\to0,\quad
 \max_jL_n^2|\beta_{j,n}|\le S_n^{\frac12}\to0.
\]
Since $p_n\to N$ and $t_n\ge c_*L_n$, we can choose a fixed
$\varrho>0$ such that the Euclidean disks $B_{2\varrho L_n}(c_{j,n})$
are embedded and disjoint. In the corresponding local coordinates,
\eqref{poinlexpa} reads
\[
 L_nF_n(z_{j,n}+L_n\eta)
 =\frac{d_j}{\eta}+L_nA_{j,n}
       +L_n^2\beta_{j,n}\eta+O(|\eta|^2),
 \qquad 0<|\eta|<2\varrho.
\]
The remainder is uniform in $j,n$. Whenever these fields converge
locally uniformly away from their poles, Cauchy's formula on
$|\eta|=\varrho$ therefore shows that the constant and linear regular
coefficients of the limit vanish. This conclusion is independent of
rotations of the local coordinates.

Pass to a further subsequence so that all center radii and pairwise
distances, divided by $L_n$, converge in $[0,\infty]$.
There are two cases, according to the distance from the closest pair
to the vertex after rescaling.

\emph{Case 1: $r(c_{1,n})/L_n$ is bounded.}
Let
\[
 \mathcal S=\left\{j:\lim_{n\to\infty}
                  \frac{r(c_{j,n})}{L_n}<\infty\right\}.
\]
The triangle inequality implies $\{1,2\}\subset\mathcal S$.
Dilate about the vertex by $L_n^{-1}$ and use the polar coordinates
$(r,\phi)$ to compare the cones. After a further subsequence,
\[
 L_n^{-1}c_{j,n}\to\widehat c_j\in C_N\setminus\{o\},\qquad
 \operatorname{dist}_N(\widehat c_j,\widehat c_k)\ge1
 \quad(j,k\in\mathcal S,\ j\ne k).
\]
Here $r(\widehat c_j)\ge c_*$, so every limiting center is distinct
from the vertex. Put $a_j=w_N(\widehat c_j)$. In compatible branches
of the local Euclidean coordinate,
\[
 L_nK_{p_n,j}(L_nz)
 =\frac{p_nz^{p_n-1}}{z^{p_n}-L_n^{-p_n}w(c_{j,n})}
 \to\frac{Nz^{N-1}}{z^N-a_j}\quad(j\in\mathcal S).
\]
For $j\notin\mathcal S$, estimate~\eqref{poinkboun} implies
\[
 |L_nK_{p_n,j}(L_nz)|
 \le\frac{C}{r(c_{j,n})/L_n-M}\to0
 \qquad (|z|\le M).
\]
Thus the rescaled fields converge locally uniformly away from the
centers and the vertex to a field on $C_N$. Since $N$ is an integer,
its lift to $\C$ is the rational function
\[
 \mathcal F(z)=\sum_{j\in\mathcal S}\frac{d_jNz^{N-1}}{z^N-a_j}
 =\sum_{j\in\mathcal S}\sum_{m=0}^{N-1}\frac{d_j}{z-\omega^m\zeta_j},
 \qquad \zeta_j^N=a_j,\quad\omega=e^{\frac{2\pi i}{N}}.
\]
The numbers $a_j$ are distinct and nonzero, so the $N|\mathcal S|$
poles are distinct. By the Cauchy argument above, the constant and
linear regular coefficients vanish at each $\zeta_j$.
The identity \mbox{$\mathcal F(\omega z)=\omega^{-1}\mathcal F(z)$} extends
this vanishing to every rotated pole $\omega^m\zeta_j$.
Lemma~\ref{cubimiden} now yields the desired contradiction.

\emph{Case 2: $r(c_{1,n})/L_n\to\infty$.}
Let
\[
 \mathcal S=\left\{j:\lim_{n\to\infty}
 \frac{\operatorname{dist}_{p_n}(c_{j,n},c_{1,n})}{L_n}
 <\infty\right\}.
\]
Again $\{1,2\}\subset\mathcal S$. For each fixed $M$, the disk
$B_{ML_n}(c_{1,n})$ is Euclidean for all sufficiently large $n$,
because $ML_n/r(c_{1,n})\to0$. Choose a branch of
$z=re^{i\phi/p_n}$ on these disks, write $z_{*,n}=z(c_{1,n})$,
and pass to a subsequence such that
\[
 b_{j,n}=\frac{z_{j,n}-z_{*,n}}{L_n}\to b_j,\quad
 b_1=0,\quad |b_2|=1,\quad |b_j-b_k|\ge1
 \quad(j,k\in\mathcal S,\ j\ne k).
\]
Fix $M>1+\max_{j\in\mathcal S}|b_j|$. Taylor expansion of $z^{p_n}$
in the kernel shows that, for $j\in\mathcal S$,
\[
 L_nK_{p_n,j}(z_{*,n}+L_nz)
 =\frac{1+O_M(L_n/r(c_{1,n}))}{z-b_{j,n}}
 \to\frac1{z-b_j}.
\]
The convergence is uniform on compact subsets of $|z|\le M$ away
from the limiting poles. For $j\notin\mathcal S$, the triangle
inequality and \eqref{poinkboun} imply
\[
 |L_nK_{p_n,j}(z_{*,n}+L_nz)|
 \le\frac{C}{\operatorname{dist}_{p_n}(c_{j,n},c_{1,n})/L_n-M}
 \to0\qquad(|z|\le M).
\]
Consequently the limiting field is the rational function
\[
 \mathcal F(z)=\sum_{j\in\mathcal S}\frac{d_j}{z-b_j}.
\]
It has at least two distinct poles, and the Cauchy argument shows that
its constant and linear regular coefficients vanish at each pole.
This again contradicts Lemma~\ref{cubimiden}. Both cases are therefore
impossible, which proves $\liminf_nS_n>0$.
\end{proof}

\begin{proposition}\label{uniteside}
Let $p_n,U_n,c_{j,n},L_n$ satisfy the assumptions of
Proposition~\ref{nonvlcoef2}, and define $S_n$ by
\eqref{regucapint}. For every sufficiently large fixed $R$,
there is a constant $C_R$ such that
\begin{equation*}
 p_n-N\ge18S_nL_n^{-4}\log L_n-C_RL_n^{-4}>0
 \qquad\text{for all sufficiently large }n.
\end{equation*}
\end{proposition}
\begin{proof}
Potential energy quantization and Proposition~\ref{poteeesti} yield
\[
 2\pi(p_n-N)=\mathcal P_{p_n}(U_n)-2\pi N
 \ge36\pi S_nL_n^{-4}\log L_n-C_RL_n^{-4}.
\]
Since $\liminf_nS_n>0$ by Proposition~\ref{nonvlcoef2},
division by $2\pi$ yields the stated inequality after changing $C_R$.
Its right side is positive for all sufficiently large $n$ because
$L_n\to\infty$.
\end{proof}

We can now conclude the argument for the positive reflection class
with $p\in[8,9]$. If vortices escaped along a sequence, then
\eqref{posiearit} would require the limiting parameter to be an odd
integer. The only possibility in $[8,9]$ is $9$, whereas
Proposition~\ref{uniteside} would imply $p_n>9$ for all sufficiently
large $n$, contradicting $p_n\le9$. Thus all zeros lie in a fixed
ball about the vertex, uniformly over this class.

\section{Weighted compactness and the formulation for degree theory}
\label{weighted} 
The preceding section excludes vortex escape for $p\in[8,9]$ in the
positive reflection class. To use degree theory, we need convergence
in a norm that also controls the behavior at the vertex and at infinity.
We first prove uniform estimates for the solutions and their derivatives.
Outside a fixed disk, write
$U=(f_p+\mathfrak a)e^{i(\phi+\psi)}$, where $\mathfrak a$ is the
difference between the modulus and $f_p$, and $\psi$ is the phase
correction. Proposition~\ref{posiwcomp} shows that
$\mathfrak a=O(r^{-4})$ and $\psi=O(r^{-2})$, uniformly in $p$.
The weighted norms measure these decay rates together with the
corresponding derivative bounds. Combined with the estimates at the
vertex and local regularity, these bounds yield convergence of a
subsequence in weaker weighted norms. Thus the convergence controls
the whole cone, including its behavior at infinity.

The second task is to put the equations for all $p\in[8,9]$ in one
fixed real Banach space. We obtain a uniformly bounded inverse by
adding a nonnegative term supported in a fixed disk to the linearized
operator. Combining this inverse with modulus and phase coordinates
at infinity yields an equivalent equation
$G(p,\xi)=\xi+\mathcal C(p,\xi)=0$, with the radial vortex represented
by $\xi=0$. The map $\mathcal C$ is compact on bounded subsets that
stay a positive distance from the relative boundary of the coordinate domain,
and the positive solutions form a compact set in these coordinates.

\subsection{Uniform estimates and compactness}
We consider a parameter interval of the form
\begin{equation}\label{onesstri}
 K=[a,b]\Subset(2,\infty),\qquad b\text{ odd},\qquad
 [a,b)\cap(2\mathbb Z+1)=\varnothing.
\end{equation}
The application is $K=[8,9]$. In the coordinate $y=re^{i\phi}$, set
\[
 A=\{y\in\R^2:1<|y|<2\},\qquad
 E_R=\{y\in\R^2:|y|>R\},\qquad
 \mathcal V=B_{2r_0}(0)\setminus\{0\},\quad 0<r_0<\frac1{16}.
\]
Fix $0<\alpha<1$. For $R\ge2$, we define the exterior weighted norms by
\begin{align*}
 \|h\|_{C_\mu^{k,\alpha}(E_R)}
 &=\sup_{s\ge R}s^\mu\|h(s\cdot)\|_{C^{k,\alpha}(A)},
 \\
 \|h\|_{U_\mu^{k,\alpha}(E_R)}
 &=\sup_{y_0\in E_R}|y_0|^\mu
                  \|h\|_{C^{k,\alpha}(B_1(y_0)\cap E_R)}.
\end{align*}
For $U=c_0+c_2r^2+\eta$ on $\mathcal V$, we use the vertex norm
\begin{equation}\label{tipwnorm}
 |c_0|+|c_2|+\sup_{0<s\le r_0}s^{-2-\sigma}
                      \|\eta(s\cdot)\|_{C^{2,\alpha}(A)},
 \qquad 0<\sigma<\min(a-2,2).
\end{equation} 

The coefficients of an actual solution are real and satisfy
$c_2=-\frac14c_0(1-c_0^2)$.
For $(p,U)\in\mathscr P_K$ with $|U|\ge\frac12$ on $E_R$, set
$\mathfrak a=|U|-f_p$. Positivity, reflection and normalization determine
a unique continuous real phase $\psi$ such that
\[
 U=(f_p+\mathfrak a)e^{i(\phi+\psi)},\quad
 \psi(r,-\phi)=-\psi(r,\phi),\quad \psi(r,0)=\psi(r,\pi)=0.
\]
The weighted topology below consists of convergence
in $p$, in the vertex norm~\eqref{tipwnorm}, and in
$C^{2,\alpha}_{\rm loc}(\R^2\setminus\{0\})$, together with convergence
of $(\psi,\mathfrak a)$ in
$C_1^{2,\alpha}(E_R)\times U_3^{2,\alpha}(E_R)$ for a fixed $R$.

\begin{proposition}\label{posiwcomp}
For $K,\sigma,\alpha$ as above, there are $R_K\ge2$ and $C_{K,\alpha}>0$
such that every $(p,U)\in\mathscr P_K$ satisfies $|U|\ge\frac12$ on
$E_{R_K}$ and
\begin{equation}\label{unifpinfty}
 \|\psi\|_{C_2^{2,\alpha}(E_{R_K})}+
 \|\mathfrak a\|_{U_4^{2,\alpha}(E_{R_K})}\le C_{K,\alpha}.
\end{equation}
Moreover, $\mathscr P_K$ is compact in the preceding weighted topology
with $R=R_K$.
\end{proposition}
\begin{proof}
\emph{Step 1: Exclusion of escape.}
By \eqref{posiearit} and \eqref{onesstri}, an escaping sequence has
a subsequence with $p_n\to b$. Proposition~\ref{uniteside} would then
imply $p_n>b$, a contradiction. Thus every decomposition in
Proposition~\ref{poteddeco} has $J=0$ and $c=1$.
Set $\rho=|U|$ and $h=1-\rho$. Strong convergence of the potential
energy integrands implies $h\to0$ uniformly at infinity: by the
uniform gradient bound, a fixed positive value of $h$ would imply
a fixed positive potential energy on a disk. Thus $\rho\ge\frac12$
outside a uniform radius and the normalized phase $\psi$ is defined
there. With $\mathbf{P}=\nabla_p(\phi+\psi)$ and $\mathbf{P}_0=\nabla_p\phi$,
the energy convergence and \eqref{phaseinfty} imply
\[
 \sup_{(p,U)\in\mathscr P_K}
 \int_{E_R}(h^2+\rho^2|\nabla_p\psi|^2)\,dA_p\to0
 \quad\text{as }R\to\infty.
\]

\emph{Step 2: Preliminary decay.}
The cutoff test in Proposition~\ref{quantize}, applied to
\eqref{modueequa} between $R$ and $2R$, and \eqref{bernboun} yield
\[
 \int_{E_{2R}}|\nabla_p\rho|^2\,dA_p\le C\int_{E_R}h^2\,dA_p.
\]
Since $|\mathbf{P}_0|\le C/r$, the energy on $B_{2\kappa r}(x)$ is at most
$C_K\kappa^2+o(1)$, uniformly as $r=r(x)\to\infty$.
Choosing $\kappa$ small and then $R_K$ large, Lemma~\ref{unifhein} yields
\begin{equation}\label{stripfirs}
 |\nabla_pU|+|\mathbf{P}|\le Cr^{-1}.
\end{equation}
We use radial barriers $Cr^{-m}$ for equations with a uniformly
positive zeroth order coefficient. Equation~\eqref{modueequa} with
$m=2$ and unit disk estimates for $\nabla_pU$ and $h$ imply
$\|h\|_{C^{2,\alpha}(B_1(x))}\le C_\alpha r^{-2}$.
For a local phase $\Theta$ of $U$, rescaling
$\Delta_p\Theta=-2\nabla_p\log\rho\cdot \nabla_p\Theta$ on disks of radius
$\kappa r$ now yields $\nabla_p\mathbf{P}=O(r^{-2})$: the rescaled drift and its
first derivative are bounded, and $r|\mathbf{P}|\le C$.
Differentiating \eqref{modueequa}, we obtain
\[
 (-\Delta_p+3\rho^2-1+|\mathbf{P}|^2)\nabla_ph=2\rho(\nabla_p\mathbf{P})^T\mathbf{P}.
\]
The coefficient is at least $1$ for large $r$ and the right hand side
is $O(r^{-3})$. Kato's inequality and the $m=3$ barrier therefore imply
$\nabla_p\rho=O(r^{-3})$.

\emph{Step 3: Weighted estimates.}
The phase equation and Step 2 imply $|\Delta_p\psi|\le Mr^{-4}$.
Moreover, \eqref{stripfirs} bounds $\partial_\phi\psi$ uniformly.
Since $\psi$ vanishes on both boundary rays,
$|\psi(r,\phi)|\le M_0\sin\phi$ for $0\le\phi\le\pi$.
The positive function
$w_p(\phi)=\frac14[\cos((2\phi-\pi)/p)/\cos(\pi/p)-1]$
vanishes at $0,\pi$ and satisfies $p^2w_p''+4w_p=-1$.
On $R<r<T$, the maximum principle bounds $\pm\psi$ by
\[
 Mr^{-2}w_p(\phi)+M_0\left(\frac Rr\right)^p\sin\phi
                  +M_0\left(\frac rT\right)^p\sin\phi.
\]
Letting $T\to\infty$ and using $p>2$, we obtain $\psi=O(r^{-2})$.
Rescaled $W^{2,q}$ estimates, with $q>2/(1-\alpha)$, then imply
$\|\nabla_p\psi\|_{C^{0,\alpha}(B_1(x))}\le C_\alpha r^{-3}$.
For the modulus difference, subtracting the radial equation yields
\[
 \bigl(-\Delta_p+|\mathbf{P}_0|^2-1+\rho^2+\rho f_p+f_p^2\bigr)\mathfrak a
 =-\rho(|\mathbf{P}|^2-|\mathbf{P}_0|^2).
\]
Its coefficient is at least $1$ for large $r$, and
$|\mathbf{P}|^2-|\mathbf{P}_0|^2=2\mathbf{P}_0\cdot \nabla_p\psi+|\nabla_p\psi|^2$ is $O(r^{-4})$
in the unit disk $C^{0,\alpha}$ norm. The $m=4$ barrier and Schauder
estimates imply $\|\mathfrak a\|_{C^{2,\alpha}(B_1(x))}\le Cr^{-4}$.
Since $\rho=f_p+\mathfrak a$, Proposition~\ref{radial} implies
$\nabla_p^2\rho=O(r^{-4})$. Differentiating
$\Delta_p\psi=-2\rho^{-1}\nabla_p\rho\cdot \mathbf{P}$ then yields
$\nabla_p\Delta_p\psi=O(r^{-5})$.
Rescaled Schauder estimates for $\psi$ prove~\eqref{unifpinfty} and
\begin{equation}\label{striuainft}
 \|U(r,\cdot)-e^{i\cdot}\|_{C^1(S^1)}\le Cr^{-2}.
\end{equation}

\Needspace{4\baselineskip}
\emph{Step 4: Compactness.}
By compactness of $K$ and Lemma~\ref{locacomp}, every sequence in
$\mathscr P_K$ has a subsequence with $p_n\to p_*\in K$ and local
uniform convergence to $U_*$, smooth away from the vertex.
Choose $\sigma<\sigma'<\min(a-2,2)$. Lemma~\ref{tipesti} bounds the
vertex remainders at scales below $\epsilon$ by
$C\epsilon^{\sigma'-\sigma}$ in~\eqref{tipwnorm}, while
\eqref{unifpinfty} implies
\[
 \|\psi_n-\psi_*\|_{C_1^{2,\alpha}(E_T)}
 +\|\mathfrak a_n-\mathfrak a_*\|_{U_3^{2,\alpha}(E_T)}\le\frac CT.
\]
The vertex coefficients converge, and local smooth convergence together
with continuity of $p\mapsto f_p$ controls the fixed annuli.
Letting $n\to\infty$, then $\epsilon\to0$ and $T\to\infty$, proves
convergence in the stated topology. The limit has finite potential
energy and the same reflection symmetry. Estimate~\eqref{striuainft}
ensures its normalization and degree one at infinity. Its imaginary
part is nonnegative, nonzero, and hence strictly positive by the strong
maximum principle. Thus $(p_*,U_*)\in\mathscr P_K$.
\end{proof}

\subsection{A coercive inverse}
Fix $I=[8,9]$, $\sigma=\delta=\frac12$ and $\alpha=\frac14$.
All spaces in this section are real spaces of fields on the fixed
plane $\R^2$. For $\mathcal V=B_{2r_0}(0)\setminus\{0\}$, we set
\[
 \|h\|_{C^{k,\alpha}(\mathcal V;r^\nu)}
 =\sup_{0<s\le r_0}s^{-\nu}\|h(s\cdot)\|_{C^{k,\alpha}(A)}.
\]
Choose $R>R_I$ and set $\mathcal A=B_{3R}(0)\setminus
\overline{B_{r_0}(0)}$. The space $X$ consists of fields satisfying
$\xi(r,-\phi)=\overline{\xi(r,\phi)}$, $\xi\in C^{2,\alpha}(\mathcal A)$,
and
\begin{equation}\label{fixedXtip}
 \xi=c_0+c_2r^2+\eta,\quad c_0,c_2\in\R,\quad
 \eta\in C^{2,\alpha}(\mathcal V;r^{2+\sigma})
\end{equation}
on $\mathcal V$, with exterior representation
\begin{equation}\label{commXinfty}
 e^{-i\phi}\xi=a+ib,\qquad
 a\in U_{\delta+2}^{2,\alpha}(E_R),\quad b\in C_\delta^{2,\alpha}(E_R),
 \quad a\text{ even},\quad b\text{ odd}
\end{equation}
on $E_R$, with norm
\begin{align*}
 \|\xi\|_X={}&|c_0|+|c_2|
 +\|\eta\|_{C^{2,\alpha}(\mathcal V;r^{2+\sigma})}
 +\|\xi\|_{C^{2,\alpha}(\mathcal A)}\\
 &+\|a\|_{U_{\delta+2}^{2,\alpha}(E_R)}
 +\|b\|_{C_\delta^{2,\alpha}(E_R)}.
\end{align*}
The space $Y$ consists of fields $\zeta$ with the same reflection,
with $\zeta-d_0\in C^{0,\alpha}(\mathcal V;r^\sigma)$ for some
$d_0\in\R$, regularity $C^{0,\alpha}(\mathcal A)$, and exterior components
$\operatorname{Re}(e^{-i\phi}\zeta)\in U_{\delta+2}^{0,\alpha}(E_R)$,
$\operatorname{Im}(e^{-i\phi}\zeta)\in C_{\delta+2}^{0,\alpha}(E_R)$.
Its norm is
\begin{align*}
 \|\zeta\|_Y={}&|d_0|+\|\zeta-d_0\|_{C^{0,\alpha}(\mathcal V;r^\sigma)}
 +\|\zeta\|_{C^{0,\alpha}(\mathcal A)}\\
 &+\|\operatorname{Re}(e^{-i\phi}\zeta)\|_{U_{\delta+2}^{0,\alpha}(E_R)}
 +\|\operatorname{Im}(e^{-i\phi}\zeta)\|_{C_{\delta+2}^{0,\alpha}(E_R)}.
\end{align*}
The vertex coefficients are unique. Agreement on overlaps and reflection
are closed conditions in the products of these weighted H\"older spaces
and their coefficient spaces of finite dimension. Thus $X$ and $Y$ are
real Banach spaces independent of $p$. The coefficients $c_0,c_2$ vary independently in $\R$. Enlarging $R$ to another fixed radius yields equivalent
norms on the same spaces.
The linearization at $V_p$ is
\begin{equation*}
 J_p\xi=-\Delta_p\xi-(1-f_p^2)\xi
               +2\operatorname{Re}(\overline{V_p}\xi)V_p.
\end{equation*}
Its relative exterior expression is
\begin{equation}\label{commjinfty}
 \begin{pmatrix}-\Delta_p+2+v_a&2p^2r^{-2}\partial_\phi\\
 -2p^2r^{-2}\partial_\phi&-\Delta_p+v_b\end{pmatrix},\quad
 v_a=p^2r^{-2}+3f_p^2-3=O(r^{-2}),\quad
 v_b=p^2r^{-2}+f_p^2-1=O(r^{-4}).
\end{equation}
For a radial cutoff $\chi_0\in C_c^\infty(\R^2)$ with
$0\le\chi_0\le1$ and $\chi_0=1$ near the vertex,
let $\mathcal H$ be the completion of compactly supported smooth
reflection fields, allowed real constants near the vertex, in the norm
\begin{equation}\label{commHnorm}
 \|\xi\|_{\mathcal H}^2=\int_{\R^2}
 \left[|\nabla_y\xi|^2+\chi_0|\xi|^2+
 (1-\chi_0)|\operatorname{Re}(e^{-i\phi}\xi)|^2\right]\,dy.
\end{equation}
The angular Poincar\'e inequality for odd functions shows that its exterior norm is equivalent to
\begin{equation}\label{commHinfty}
 \int_{E_R}(|\nabla_ya|^2+|\nabla_yb|^2+a^2+r^{-2}b^2)\,dy.
\end{equation}
This inequality also shows that $\mathcal H$ is independent of the cutoff, with equivalent norms. The forms below use $\,dy=r\,dr\,d\phi=p\,dA_p$.

For the radial solution $f_p$ of Proposition~\ref{radial},
the quadratic form of $J_p$ on $\mathcal H$ is
\begin{equation}\label{jacofform}
 Q_p[\xi]=\int_0^\infty\!\int_0^{2\pi}
 \left[|\xi_r|^2+\frac{p^2}{r^2}|\xi_\phi|^2-(1-f_p^2)|\xi|^2
 +2\bigl(\operatorname{Re}(\overline{V_p}\xi)\bigr)^2\right]r\,d\phi\,dr.
\end{equation}
 We add a fixed compactly supported zeroth order term to $J_p$.
This yields a uniformly coercive quadratic form and an inverse
that varies continuously with $p$.
\begin{proposition}\label{coerline}
There exist a radial $\chi_0\in C_c^\infty(\R^2)$, with
$0\le\chi_0\le1$ and $\chi_0=1$ near $0$, and constants $C_0,c>0$,
independent of $p\in I$, such that
\begin{equation*}
 B_p=J_p+C_0\chi_0:X\to Y
\end{equation*}
is a bounded isomorphism for every $p\in I$. The families $J_p,B_p$ and $B_p^{-1}$
are continuous in the corresponding operator norms.
Moreover, $J_p:X\to Y$ is Fredholm of index zero and
$\ker_XJ_p=\ker_{\mathcal H}J_p$. For every $p\in I$ and $\xi\in\mathcal H$,
\begin{equation}\label{commcoer}
 Q_p[\xi]+C_0\int_{\R^2}\chi_0|\xi|^2\,dy
       \ge c\|\xi\|_{\mathcal H}^2.
\end{equation}
\end{proposition}

\begin{proof}
To solve $B_p\xi=f$ for $f\in Y$, we first obtain weighted estimates
for $J_p$ outside a large disk. We then choose $\chi_0$ and $C_0$
so that coercivity yields a unique weak solution in $\mathcal H$.
Regularity at the vertex and the exterior estimates show that this
solution belongs to $X$.

\emph{Exterior inverse.}
For $S\ge R$, all exterior inverses have zero trace on
$\partial B_S(0)$. We use the even and odd components of
\eqref{commXinfty}, with the corresponding component norms of $Y$.
The barrier $r^{-\delta-2}$ and unit disk Schauder estimates yield
the first of the scalar inverses
\begin{align*}
 (-\Delta_p+2)^{-1}&:U_{\delta+2}^{0,\alpha}(E_S)
                    \longrightarrow U_{\delta+2}^{2,\alpha}(E_S),\\
 (-\Delta_p)^{-1}&:C_{\delta+2}^{0,\alpha}(E_S)
                    \longrightarrow C_\delta^{2,\alpha}(E_S).
\end{align*}
For the second, acting on odd functions, use $r^{-\delta}q_p(\phi)$,
where $-p^2q_p''-\delta^2q_p=1$ and $q_p(0)=q_p(\pi)=0$.
Since $\delta<p$, the function $q_p$ is uniformly comparable to
$\sin\phi$. The maximum principle, odd reflection and rescaled
boundary Schauder estimates prove the required bound.
Solving first on bounded annuli and taking limits constructs both
inverses, with norms uniform in $p\in I$ and large $S$.

We first solve the upper triangular part of~\eqref{commjinfty}:
\[
 -\Delta_p b=\operatorname{Im}(e^{-i\phi}f),\qquad
 (-\Delta_p+2)a=\operatorname{Re}(e^{-i\phi}f)-2p^2r^{-2}b_\phi.
\]
The map $b\mapsto r^{-2}b_\phi$ is uniformly bounded between the
required spaces. For the other coupling, $a_\phi=re_\phi\cdot\nabla_ya$
and rescaling the H\"older norm imply
\[
 \|r^{-2}a_\phi\|_{C_{\delta+2}^{0,\alpha}(E_S)}
 \le CS^{-1+\alpha}\|a\|_{U_{\delta+2}^{2,\alpha}(E_S)}.
\]
Multiplication by $v_a$ and $v_b$ has norms $O(S^{-2})$ and
$O(S^{-2+\alpha})$, respectively. Thus the remaining terms composed
with the triangular inverse have norm at most
$C(S^{-2}+S^{-1+\alpha})$. After increasing $R$, this is less than
$1/2$ for every $p\in I$ and $S\ge R$. A Neumann series yields a unique
solution $z=e^{i\phi}(a+ib)$ of $J_pz=f$ with zero trace and
\begin{align*}
 \|a\|_{U_{\delta+2}^{2,\alpha}(E_S)}
 +\|b\|_{C_\delta^{2,\alpha}(E_S)}
 \le C\left(
 \|\operatorname{Re}(e^{-i\phi}f)\|_{U_{\delta+2}^{0,\alpha}(E_S)}
 +\|\operatorname{Im}(e^{-i\phi}f)\|_{C_{\delta+2}^{0,\alpha}(E_S)}\right).
\end{align*}

\emph{Coercivity and the global inverse.}
The weights in~\eqref{commHinfty} imply the continuous embeddings
$X\hookrightarrow\mathcal H$ and $Y\hookrightarrow\mathcal H^*$,
with $X$ dense in $\mathcal H$. In particular, the imaginary component
of $f\in Y$ pairs with $b$ through
\[
 \left|\int_{E_R}\operatorname{Im}(e^{-i\phi}f)b\,dy\right|
 \le\|r\operatorname{Im}(e^{-i\phi}f)\|_{L^2(E_R)}
       \|b/r\|_{L^2(E_R)},
\]
where the first norm is finite because $\delta>0$.
The terms in the exterior quadratic form satisfy
\[
 \frac{4p^2}{r^2}|ab_\phi|\le\varepsilon a^2+
 \frac{C_{I,\varepsilon}}{r^4}|b_\phi|^2,\qquad
 \int_{E_R}|v_b|b^2\,dy\le\frac C{R^2}
                \int_{E_R}\frac{|b_\phi|^2}{r^2}\,dy.
\]
Choose $\varepsilon$ small and then enlarge $R$.
We fix a radial cutoff $\chi_0=1$ on $B_{2R}(0)$ and
$\chi_0=0$ outside $B_{3R}(0)$. For sufficiently large $C_0$,
the added term absorbs the bounded negative potential in the fixed
disk and proves~\eqref{commcoer}, uniformly in $p\in I$.
The theorem of Lax and Milgram therefore constructs
$B_p^{-1}:\mathcal H^*\to\mathcal H$.

At the vertex, the conformal representative $g$ of an $\mathcal H$
solution belongs locally to $H^1$ and satisfies
$-\Delta_wg=a_p(w)(M(w)g+h(w))$, with $M$ a bounded real linear map
and $h$ bounded.
For a fixed $1<s<9/8$, the weights $a_p$ are uniformly bounded in $L^s$.
The embedding $H^1\subset L^t$ for every finite $t$, H\"older's
inequality and interior $W^{2,q}$ estimates for some $q>1$ imply
boundedness of $g$.
The argument of Lemma~\ref{tipesti} then yields $\xi=c_0+O(r^2)$
for a bounded right hand side. For $B_p\xi=f\in Y$, substitution
into the equation and rescaled Schauder estimates yield
\begin{gather*}
 \xi=c_0+c_2r^2+\eta,\qquad
 c_2=\frac{(C_0-1)c_0-f(0)}4,\\
 |c_0|+|c_2|+\|\eta\|_{C^{2,\alpha}(\mathcal V;r^{2+\sigma})}
 \le C(\|f\|_Y+\|\xi\|_{\mathcal H}).
\end{gather*}
The iteration is uniform since $2+\sigma<8\le p$, and reflection
makes the coefficients real. For $\xi\in X$, direct substitution
also shows that $J_p\xi=-4c_2-c_0+O_{C^{0,\alpha}}(r^\sigma)$.
Together with~\eqref{commjinfty} and interior estimates, this proves
uniform boundedness of $J_p:X\to Y$.

Let $f\in Y$ and $\xi=B_p^{-1}f\in\mathcal H$.
Choose a fixed $T>3R$ and extend the trace of $\xi$ on $\partial B_T(0)$
to a $C^{2,\alpha}$ reflection field $e$ on $\overline{E_T}$,
supported in $T\le r\le2T$.
The exterior inverse applied to $f-J_pe$, followed by addition of $e$,
produces a weighted solution $z$ of $J_pz=f$ on $E_T$ with the same
boundary values as $\xi$. Its weighted norm is bounded by
$C_T(\|f\|_Y+\|\xi\|_{\mathcal H})$.
Since $\delta>0$, both solutions have finite $\mathcal H$ norm there.
The zero extension of $\xi-z$ belongs to $\mathcal H$, and radial
cutoffs justify testing its equation with $\xi-z$. Since $\chi_0=0$
on $E_T$, coercivity implies
\[
 0=Q_p[\xi-z]\ge c\|\xi-z\|_{\mathcal H}^2.
\]
Thus $\xi=z$ on $E_T$. The vertex and interior estimates now imply
\[
 \|B_p^{-1}f\|_X
 \le C(\|f\|_Y+\|B_p^{-1}f\|_{\mathcal H})\le C\|f\|_Y.
\]
Hence the variational inverse restricts to a bounded inverse from
$Y$ to $X$.

\emph{Compactness and parameter dependence.}
Multiplication by $\chi_0$ is compact from $X$ to $Y$, by the finite
dimension of the vertex coefficients, the factor $r^{2-\sigma}$
near zero and compact embedding on annuli. It is compact from
$\mathcal H$ to $\mathcal H^*$ by Rellich compactness.
Consequently, $J_p=B_p(\operatorname{Id}-C_0B_p^{-1}\chi_0)$ is
Fredholm of index zero. If $J_p\xi=0$ in $\mathcal H$, the vertex
argument above and interior regularity imply $\chi_0\xi\in Y$.
Thus $\xi=C_0B_p^{-1}(\chi_0\xi)\in X$, proving equality of the kernels.

The coefficient estimates at zero, on fixed annuli and at infinity yield
\[
 \|J_p-J_q\|_{\mathcal L(X,Y)}
 \le C|p^2-q^2|
 +C_{\epsilon,T}\|f_p-f_q\|_{C^{2,\alpha}([\epsilon,T])}
 +C(\epsilon^{2-\sigma}+T^{-2+\alpha}).
\]
We first choose $\epsilon$ small and $T$ large, and then let $q\to p$.
Continuity of $p\mapsto f_p$ proves continuity of $J_p$ in operator
norm. Since $B_p-J_p=C_0\chi_0$ is independent of $p$, the same holds
for $B_p$, and the resolvent identity proves continuity of $B_p^{-1}$.
The same coefficient bounds and the angular Poincar\'e inequality
prove continuity of the quadratic forms $Q_p$ on $\mathcal H$.
\end{proof}

For $p\in I$, let $B_p,\chi_0,C_0$ be as in
Proposition~\ref{coerline}, and set

\begin{align*}
 A_p=B_p^{-1}J_p=\operatorname{Id}-C_0B_p^{-1}\chi_0,\quad\mbox{and}\quad
 \mathfrak b_p(\xi,\zeta)=Q_p(\xi,\zeta)+
 C_0\int_{\R^2}\chi_0\operatorname{Re}(\xi\bar\zeta)\,dy.
\end{align*}
 
\begin{lemma}\label{negaelift}
For $p\in I$, the operator $A_p$ extends to a compact perturbation of
the identity on $\mathcal H$ and is selfadjoint with respect to
$\mathfrak b_p$. Its eigenspaces in $X$ and $\mathcal H$ coincide for
every $\lambda\le0$. For each negative eigenvalue, its algebraic
multiplicity in either space equals the dimension of this eigenspace.
Moreover,
\[
 \ind(Q_p;\mathcal H)=\sum_{\lambda<0}\dim\ker_X(A_p-\lambda).
\]
\end{lemma}
\begin{proof}
By the continuity of $Q_p$ and~\eqref{commcoer}, the inner product
$\mathfrak b_p$ defines a norm
equivalent to $\|\cdot\|_{\mathcal H}$. The proof of
Proposition~\ref{coerline} constructs $B_p^{-1}:\mathcal H^*\to\mathcal H$
and proves that multiplication by $\chi_0$ is compact from
$\mathcal H$ to $\mathcal H^*$. Hence $C_0B_p^{-1}\chi_0$ is compact
on $\mathcal H$. The variational definition of $B_p^{-1}$ implies,
for $\xi,\zeta\in\mathcal H$,
\[
 \mathfrak b_p(A_p\xi,\zeta)=Q_p(\xi,\zeta)
 =\mathfrak b_p(\xi,A_p\zeta),
\]
which proves selfadjointness.

The weights defining $X$ imply $X\subset\mathcal H$.
Conversely, if $\xi\in\mathcal H$ satisfies $A_p\xi=\lambda\xi$
with $\lambda\le0$, then
\[
 B_p\xi=\frac{C_0}{1-\lambda}\chi_0\xi.
\]
The vertex estimates and interior regularity used in the proof of
Proposition~\ref{coerline} imply $\chi_0\xi\in Y$.
Since the variational inverse agrees with the inverse from $Y$ to $X$,
we obtain
\[
 \xi=\frac{C_0}{1-\lambda}B_p^{-1}(\chi_0\xi)\in X.
\]
Thus the eigenspaces coincide.

We now apply the compact selfadjoint spectral theorem
\cite[Theorems~6.8 and~6.11]{BrezisFA} to $C_0B_p^{-1}\chi_0$ in
$(\mathcal H,\mathfrak b_p)$. The negative eigenspaces of $A_p$ have
finite total dimension, and every generalized eigenvector for a
negative eigenvalue belongs to its eigenspace. Since $X\subset\mathcal H$
and the eigenspaces coincide, the same assertion holds in $X$.
This proves the equality of algebraic and geometric multiplicities
in both spaces. Finally, the identity
$Q_p[\xi]=\mathfrak b_p(A_p\xi,\xi)$ shows that $Q_p$ is negative
definite on the sum of the negative eigenspaces and nonnegative on
its $\mathfrak b_p$ orthogonal complement. Their total dimension
is therefore $\ind(Q_p;\mathcal H)$.
\end{proof}
 
\subsection{A compact perturbation of the identity}
We use $X$ from \eqref{fixedXtip}--\eqref{commXinfty} and $Y,B_p$
from Proposition~\ref{coerline}. To express the equation in these spaces, we use additive
coordinates near the vertex and modulus and phase coordinates at infinity.
Choose a smooth radial cutoff $0\le\chi\le1$ with $\chi=0$ for $r\le R$
and $\chi=1$ for $r\ge2R$. For $e^{-i\phi}\xi=a+ib$, define
\begin{align*}
 T_p(\xi)&=V_p+\xi+
 \chi\left[e^{i\phi}(f_p+a)e^{ib/f_p}-V_p-\xi\right],\\
 \mathcal D&=I\times\{\xi\in X:\sup_{r\ge R}(|a|+|b|)<\eta\},
 \qquad \mathcal D_p=\{\xi:(p,\xi)\in\mathcal D\}.
\end{align*}
Thus $T_p(0)=V_p$, and $\mathcal D$ is relatively open in $I\times X$
and contains $I\times\{0\}$.

\begin{proposition}\label{commpequa}
For suitable $R$ and $\eta>0$, each
$T_p:\mathcal D_p\to C^0(\R^2)$ is injective, and there exists
$G(p,\xi)=\xi+\mathcal C(p,\xi)$ from $\mathcal D$ to $X$ such that
\[
 G(p,\xi)=0\quad\Longleftrightarrow\quad(p,T_p(\xi))\in\mathscr R_I,
 \qquad D_\xi G(p,0)=B_p^{-1}J_p.
\]
The map $G$ is $C^1$ in $\xi$, with $G$ and $D_\xi G$ jointly continuous.
The map $\mathcal C$ is jointly compact on bounded subsets at positive
distance from the relative boundary of $\mathcal D$.
Every member of $\mathscr P_I$ has a unique preimage, and these
preimages form a compact subset of $\mathcal D$.
Convergence of zeros in $I\times X$ implies convergence
in~\eqref{posictopo}.
\end{proposition}
\begin{proof}
\emph{Coordinates.}
For small $\eta$, the real pointwise derivative $M_p(\xi)=D_\xi T_p(\xi)$
satisfies $\|M_p(\xi)-\operatorname{Id}\|<\frac12$.
The pointwise domain is convex, so integration along segments yields
\[
 |T_p(\xi)(y)-T_p(\zeta)(y)|\ge\frac12|\xi(y)-\zeta(y)|.
\]
The contraction theorem also constructs an inverse near $V_p$.
By Proposition~\ref{posiwcomp}, all positive solutions lie in this
image once $R$ is large, and their preimages stay a positive distance
from the complement of $\mathcal D$. The inverse is additive for
$r\le R$, smooth on the transition annulus, and satisfies
$a=|U|-f_p=\mathfrak a$, $b=f_p\psi$ for $r\ge2R$.
The compactness in Proposition~\ref{posiwcomp}, the inequalities
$3>\delta+2$, $1>\delta$, and continuity of $p\mapsto f_p$
therefore imply compactness of these preimages in $I\times X$.

\emph{The equation in $X$.}
We define
\[
 \mathcal F(p,\xi)=M_p(\xi)^{-1}
 [-\Delta_pT_p(\xi)-T_p(\xi)(1-|T_p(\xi)|^2)],\qquad
 \mathcal N(p,\xi)=\mathcal F(p,\xi)-J_p\xi.
\]
The factor $M_p(\xi)^{-1}$ reduces the principal part to $-\Delta_p\xi$,
so $\mathcal N$ contains derivatives of order at most one.
Since $T_p(0)=V_p$ and $M_p(0)=\operatorname{Id}$,
$\mathcal F(p,0)=0$ and $D_\xi\mathcal F(p,0)=J_p$.
On $E_{2R}$, put $\rho=f_p+a$ and $\psi=b/f_p$.
Using $-\Delta_pb+v_bb=-f_p\Delta_p\psi-2\nabla_pf_p\cdot\nabla_p\psi$
to subtract the linear part, we obtain $e^{-i\phi}\mathcal N=N_a+iN_b$, where
\begin{align*}
 N_a&=3f_pa^2+a^3+2a\nabla_p\phi\cdot\nabla_p\psi
                         +(f_p+a)|\nabla_p\psi|^2,\\
 N_b&=\frac{2a}{f_p+a}\nabla_pa\cdot\nabla_p\phi
       +\frac{2(a\nabla_pf_p-f_p\nabla_pa)}{f_p+a}\cdot\nabla_p\psi.
\end{align*}

We work on a bounded subset of $\mathcal D$ at positive distance
from its relative boundary. Here $\rho\ge c>0$,
$\|a\|_{C^{2,\alpha}(B_1(x))}\le Cr^{-\delta-2}$ and
$\|\psi(s\cdot)\|_{C^{2,\alpha}(A)}\le Cs^{-\delta}$.
The two terms of $N_b$ have respective unit disk H\"older bounds
$O(r^{-2\delta-5})$ and $O(r^{-2\delta-3})$.
Product estimates and the factor $s^\alpha$ from rescaling imply,
for $S\ge2R$,
\[
 \|N_a\|_{U_{\delta+2}^{0,\alpha}(E_S)}\le CS^{-\delta},\qquad
 \|N_b\|_{C_{\delta+2}^{0,\alpha}(E_S)}
 \le CS^{-\delta-1+\alpha}=CS^{-5/4}.
\]
Near the vertex, $T_p(\xi)=V_p+\xi$ implies
$\mathcal N=c_0^3+O_{C^{0,\alpha}}(r^2)$, whose nonconstant
contribution to the $Y$ norm on $B_\epsilon(0)$ is
$O(\epsilon^{2-\sigma})$.
For any sequence in this bounded subset, we first take subsequences
in $p$ and $c_0$ and use compact embedding from $C^{2,\alpha}$ into
$C^{1,\alpha}$ on fixed annuli. The uniform bounds above then imply
convergence of $\mathcal N$ in $Y$, proving joint compactness.

The formal derivative $D_\xi\mathcal N(p,\xi)[\zeta]$ satisfies the same
exterior bounds uniformly for $\|\zeta\|_X\le1$. At the vertex its
constant term is $3c_0^2\zeta(0)$, with the same
$O(\epsilon^{2-\sigma})$ bound for the remaining contribution.
Product and chain rules on fixed annuli, followed by
$\epsilon\to0$ and $S\to\infty$, prove joint continuity in
$\mathcal L(X,Y)$ of this derivative. For small $h\in X$, the integral
formula then yields
\begin{align*}
 &\|\mathcal N(p,\xi+h)-\mathcal N(p,\xi)-D_\xi\mathcal N(p,\xi)[h]\|_Y\\
 &\quad\le\sup_{0\le t\le1}
 \|D_\xi\mathcal N(p,\xi+th)-D_\xi\mathcal N(p,\xi)\|_{\mathcal L(X,Y)}
 \|h\|_X=o(\|h\|_X).
\end{align*}
Thus $\mathcal N$ is $C^1$ in $\xi$.
By Proposition~\ref{coerline}, the maps
\begin{align*}
 \mathcal C(p,\xi)&=B_p^{-1}[\mathcal N(p,\xi)-C_0\chi_0\xi],\\
 G(p,\xi)&=B_p^{-1}\mathcal F(p,\xi)=\xi+\mathcal C(p,\xi)
\end{align*}
have the stated compactness, continuity and derivative properties.

\emph{Zeros and solutions.}
Since $B_p$ and $M_p(\xi)$ are invertible, $G(p,\xi)=0$ is equivalent
to the equation for $U=T_p(\xi)$ away from the vertex.
The vertex expansion implies $U\in H^1(B_{r_0}(o))$ and, for each
smooth compactly supported test function $\varphi$,
\[
 \left|\int_{\partial B_\epsilon(o)}\nabla_pU\cdot\nu\,\varphi\,ds_p\right|
 \le C_\varphi\epsilon(\epsilon+\epsilon^{1+\sigma}+\epsilon^{p-1})\to0.
\]
Thus the weak equation extends across the vertex.
At infinity, $\rho=1+O(r^{-2})$ and $\psi=O(r^{-\delta})$, so
\begin{gather*}
 \int_{E_{2R}}(1-|U|^2)^2\,dA_p\le C\int_{2R}^\infty r^{-3}\,dr<\infty,\\
 \deg(U,\infty)=\frac1{2\pi}\int_0^{2\pi}(1+\psi_\phi)\,d\phi=1.
\end{gather*}
The formula for $T_p$ also implies reflection and normalization.
Hence every zero represents a member of $\mathscr R_I$.
Finally, convergence in $X$ implies local convergence and controls
$b,b_\phi$ and $|a_\phi|\le C\|\xi\|_Xr^{-\delta-1}$ at infinity.
Together with continuity of $p\mapsto f_p$, this proves convergence
in~\eqref{posictopo}.
\end{proof}

\section{Morse indices and the proof of the main theorem}\label{indices}

In this final section, we will first compute the Morse indices of the radial vortices at $p=8,9$
in the real reflection space $\mathcal H$ of \eqref{commHnorm}. (We emphasize that this is different from the full Morse index.)  Fourier mode decomposition then reduces the index
and kernel calculations to a single pair of angular modes, described
by two coupled radial functions. In this respect, we follow the Fourier reduction and positive supersolution argument
of \cite[\S2--3]{DPLMWY}, developed there for degrees two and three.
Positive supersolutions bound the index from above, while explicit
negative subspaces yield the lower bound. With the radial
inequalities proved below and in Appendix~\ref{radialineq}, we obtain
indices $3$ and $4$ at $p=8$ and $p=9$, respectively. 

These index values enter the existence argument through their parity.
In the formulation $G(p,\xi)=0$ from Section~\ref{weighted}, the radial
vortex corresponds to $\xi=0$. The index calculation therefore yields
local degrees $-1$ and $1$ at $p=8$ and $p=9$, respectively.
The compactness of the positive solution class allows us to define
a total degree that stays constant as $p$ varies over $[8,9]$.
If the radial vortex were the only solution at both parameters,
these two local degrees would equal the same total degree.
Their opposite signs therefore imply an additional positive solution
for at least one $q\in\{8,9\}$.   

We now decompose the quadratic form~\eqref{jacofform} into angular Fourier modes.
Reflection pairs the Fourier modes in the orthogonal subspaces
\[
 \mathcal H_0=\{e^{i\phi}h(r):h\text{ real}\}\cap\mathcal H,\qquad
 \mathcal H_k=\{Ae^{i(1-k)\phi}+Be^{i(1+k)\phi}:A,B\text{ real}\}
 \cap\mathcal H\quad\text{for }k\ge1.
\]
The radial equation and integration by parts yield
\begin{align*}
 \frac{Q_p[e^{i\phi}h]}{2\pi}
 &=\int_0^\infty\left[rf_p^2\left|\left(\frac h{f_p}\right)'\right|^2
                  +2rf_p^2h^2\right]\,dr,\\
 \frac{Q_p[Ae^{i(1-k)\phi}+Be^{i(1+k)\phi}]}{2\pi}
 &=\int_0^\infty\biggl[rf_p^2\left(
 \left|\left(\frac A{f_p}\right)'\right|^2+
 \left|\left(\frac B{f_p}\right)'\right|^2\right)\\
 &\hspace{8mm}+\frac{p^2}{r}\bigl(k(k-2)A^2+k(k+2)B^2\bigr)
                 +rf_p^2(A+B)^2\biggr]\,dr.
\end{align*}
For $k\ge2$, the integral of $r|u'|^2+r^{-1}u^2$ is finite for
$u=A,B$, and the same holds for $u=h$. Together with
$0<f_p'/f_p<p/r$, this justifies radial cutoff approximation in these
identities. Their equality cases imply $h=0$ and $B=A+B=0$,
respectively. Hence $Q_p$ is strictly positive on $\mathcal H_0$ and
$\mathcal H_k$ for $k\ge2$. Parseval's identity shows that projection
onto $\mathcal H_1$ is injective on every negative subspace. Testing a
weak kernel element against its Fourier components also yields
\begin{equation}\label{morseone}
 \ind(Q_p;\mathcal H)=\ind(Q_p;\mathcal H_1),\qquad
 \ker_{\mathcal H}J_p\subset\mathcal H_1.
\end{equation}

In $\mathcal H_1$, let $R=r/p$, $F_p(R)=f_p(pR)$ and
$\xi=(u(R)+v(R)e^{2i\phi})/\sqrt2$. By integration in the angular variable, we obtain
$Q_p[\xi]=\pi p^2q_p[u,v]$, where
\begin{equation}\label{rescaform}
 q_p[u,v]=\int_0^\infty R\left[
 \frac{u'^2+v'^2}{p^2}+(2F_p^2-1)(u^2+v^2)
 +2F_p^2uv+\frac{4v^2}{R^2}\right]\,dR.
\end{equation}
Its domain $\mathcal K$ consists of real pairs $z=(u,v)$ in
$H^1_{\rm loc}((0,\infty))^2$ satisfying
\begin{equation}\label{rescfdoma}
 \int_0^\infty\left(R|z'|^2+\frac{v^2}{R}\right)\,dR
 +\int_0^1R|z|^2\,dR+\int_1^\infty R(u+v)^2\,dR<\infty.
\end{equation}
This condition also implies $\int_1^\infty R^{-1}|z|^2\,dR<\infty$,
while $u$ may have a nonzero value at the origin.

\begin{proposition}\label{endpindi}
For the operators and space defined above,
\[
 \ind(Q_8;\mathcal H)=3,\qquad\ind(Q_9;\mathcal H)=4,\qquad
 \ker_{\mathcal H}J_8=\ker_{\mathcal H}J_9=\{0\}.
\]
\end{proposition}
\begin{proof}
Let $K_8=3$ and $K_9=4$. We apply the scalar reduction and evaluation
argument of \cite[Lemma~3.1 and Proposition~3.2]{DPLMWY}, checking its
domain and the required strict inequalities for $p=8,9$.

\emph{Upper bound.}
Let $t=\log R$, $x=e^{2t}$, with dots denoting derivatives in $t$.
For the scalar reduction, define
\[
 \mathfrak p=x(1-F_p(\sqrt x)^2),\qquad
 V=-\mathfrak p+\frac{(4-\mathfrak p)(x-\mathfrak p)}{x+4-2\mathfrak p},\qquad
 \mathscr S_p=-p^{-2}\partial_t^2+V.
\]
Lemma~\ref{algerboun} implies
$x+4-2\mathfrak p>4-x>0$ for $x<4$ and $x+4-2\mathfrak p>x+7/4$ for $x\ge4$.
Completing the square in \eqref{rescaform}, we obtain
\begin{equation}\label{scalscomp}
\begin{split}
 q_p[u,v]={}&\int_\R\left[\frac{\dot v^2}{p^2}
 +(x+4-2\mathfrak p)\left(v+\frac{x-\mathfrak p}{x+4-2\mathfrak p}u\right)^2\right]\,dt
 +\ell_p[u],\\
 \ell_p[u]={}&\int_\R\left(\frac{\dot u^2}{p^2}+Vu^2\right)\,dt.
\end{split}
\end{equation}
In the domain \eqref{rescfdoma}, we have
$\dot u,\dot v,v\in L^2(\R)$, $e^tu\in L^2(({-\infty},0))$ and
$e^t(u+v),u\in L^2((0,\infty))$. Since
$V=O(\min\{e^{2t},1\})$, every integral in
\eqref{scalscomp} is finite.

Let $j_m$ be the $m$th positive zero of $J_0$. Choose
$a_{p,j}=j_j^2/p^2$ for $1\le j\le3$ and $a_{9,4}=5/3$.
The Bessel bounds in \cite[(10.18.18), the remainder estimate following
(10.18.21), and (10.21.12)]{DLMF}
imply $173/20<j_3<9$ and $j_4>93/8$, and consequently
$7/6<a_{8,3}<4$ and $a_{9,3}<5/3<j_4^2/81$.
Let $t_j=\frac12\log a_{p,j}$, $t_0=-\infty$,
$t_{K_p+1}=\infty$ and $I_j=(t_{j-1},t_j)$. On the first three intervals
take $\psi_j=(-1)^{j-1}J_0(pe^t)$, and for $p=9$ take
$\psi_4=-J_0(9e^t)$ on $I_4$. Since
$(-p^{-2}\partial_t^2-x)\psi_j=0$ and
\[
 V+x=(x-\mathfrak p)\left(1+\frac{4-\mathfrak p}{x+4-2\mathfrak p}\right)>0
 \quad\text{for }x<4,
\]
these functions satisfy $\psi_j>0$ and $\mathscr S_p\psi_j>0$.
On the last interval use, respectively,
\[
 \psi_4=(x-7/6)x^{-17/3}e^{-6/x}\quad\text{for }p=8,\qquad
 \psi_5=(x-5/3)x^{-2}\quad\text{for }p=9.
\]
Lemma~\ref{eightcomp} proves the first strict inequality.
For $p=9$, the bound $\mathfrak p<9/8$ in Lemma~\ref{algerboun}
and direct differentiation imply $\mathscr S_9\psi_5>0$ on $I_5$. 

Suppose $u(t_j)=0$ for $1\le j\le K_p$. On each $I_j$, integration by
parts against the positive function $\psi_j$ yields
\[
 \ell_p[u]\ge\sum_{j=1}^{K_p+1}
 \int_{I_j}\frac{\mathscr S_p\psi_j}{\psi_j}u^2\,dt.
\]
The identity extends from compactly supported functions by cutoff
approximation. At finite boundaries, the zero trace and local
$H^1$ regularity suffice. At $-\infty$, $\dot u\in L^2$ implies
$u(t)=o(\sqrt{|t|})$, so $T^{-2}\int_{-2T}^{-T}u^2\,dt=o(1)$.
At $+\infty$, $u\in H^1((0,\infty))$.
Thus \eqref{scalscomp} is strictly positive for
every nonzero pair whose first component vanishes at all the chosen
points. Evaluation at those $K_p$ points is therefore injective on
every nonpositive subspace $S\subset\mathcal K$, and
\begin{equation}\label{analnboun}
 q_p|_S\le0\quad\Longrightarrow\quad\dim S\le K_p.
\end{equation}

\emph{Negative directions.}
The radial equation for $F_p$ is
\begin{equation}\label{scalrequa}
 \mathcal L_pF_p=-p^{-2}(F_p''+R^{-1}F_p')
                    +(R^{-2}-1+F_p^2)F_p=0.
\end{equation}
For $U_p=(x/(x+1/p))^{p/2}$, we have
\[
 \frac{\mathcal L_pU_p}{U_p}
 =\frac{p}{1+px}+\frac{p+2}{(1+px)^2}-1+
   \left(\frac{px}{1+px}\right)^p>0.
\]
The quotient identity \cite[(2.5)]{DPLMWY}, with
$RU_p^2(F_p/U_p)'\to0$ at zero and $F_p/U_p\to1$ at infinity,
implies $F_p<U_p$ and $F_p^2\le p^px^p$.
For $\beta_p=29/(5p^2)$, take
\[
 z_{p,1}(R)=\bigl((1-R^2/\beta_p)(1-R^2/(3\beta_p)),0\bigr)
 \quad\text{for }0\le R^2\le\beta_p,
\]
extended by zero. Integrating this polynomial and using the preceding
bound, we obtain $q_p[z_{p,1}]<0$.

For the next tests let $r=pR$ and $u(R)=r^{-1/2}w(r)$, where
$w\in H_0^1((a,b))$ is extended by zero. Integrating by parts, we obtain
\begin{equation}\label{radirtran}
 p^2q_p[u,0]=\int_a^b
 \left[w'^2+\left(2f_p(r)^2-1-\frac1{4r^2}\right)w^2\right]\,dr.
\end{equation}
For $0<\alpha<p$, evaluation of the radial equation at a maximum of
$r^{-\alpha}f_p(r)$ and maximization of $s(1-s)^\alpha$ imply
\begin{equation}\label{poweruppe}
 f_p(r)^2\le\frac{\alpha^\alpha}{(\alpha+1)^{\alpha+1}}
                  \left(\frac{r^2}{p^2-\alpha^2}\right)^\alpha.
\end{equation}
For $p=8$, choose
$w=\sin(\pi(r-a)/(b-a))$ on $[241/100,561/100]$ and on $[6,11]$.
Using \eqref{poweruppe} with $\alpha=4$ and $\alpha=1$,
respectively, direct integration in \eqref{radirtran} yields
$q_8[u,0]<0$ for both intervals.
For $p=9$, the same sine function on $[241/100,111/20]$ and
\eqref{poweruppe} with $\alpha=5$ yield
$q_9[u,0]<0$. On $[a,b]=[111/20,437/50]$, choose instead
$w=\sigma-2\sigma^3+\sigma^4$, where $\sigma=(r-a)/(b-a)$.
Using $F_9(R)^2<R^{10}/9$ from
Lemma~\ref{algerboun}, polynomial integration in
\eqref{radirtran} yields $q_9[u,0]<0$.

For the last direction at $p=9$, set $a=17/18$, $b=5/4$ and
$g(x)=(4x-3)/10$. Lemma~\ref{algerboun} implies
$F_9(\sqrt a)^2<1/13<g(a)$ and $F_9(\sqrt b)^2<1/5=g(b)=1-1/b$.
Direct differentiation shows that $\sqrt{g(R^2)}$ on
$[\sqrt a,\sqrt b]$ and $\sqrt{1-R^{-2}}$ on $[\sqrt b,\infty)$
are strict supersolutions for $\mathcal L_9$. The quotient identity,
together with $F_9(R)/\sqrt{1-R^{-2}}\to1$, thus yields
\[
 F_9(\sqrt x)^2\le
 \begin{cases}g(x),&a\le x\le b,\\1-1/x,&x\ge b.\end{cases}
\]
Let $u(x)=x(x-a)(3-x)^3$ on $[a,3]$, extended by zero, and take
$z_{9,4}(R)=(u(R^2),-u(R^2)/8)$. Substitution into
\eqref{rescaform}, using these upper bounds reduces the sign to
polynomial integrals. Direct evaluation yields $q_9[z_{9,4}]<0$.

These $K_p$ functions belong to $\mathcal K$ and have disjoint support
interiors.  Their span $E_p$ is therefore negative definite
and has dimension $K_p$. In view of \eqref{analnboun},
$\ind(q_p;\mathcal K)=K_p$. If $z$ were a nonzero weak kernel element,
$E_p\oplus\R z$ would be a nonpositive space of dimension $K_p+1$.
This contradicts the same bound, and \eqref{morseone} proves the proposition.
\end{proof}
 
We now combine the index calculation with the compactness from
Section~\ref{weighted} to complete the proof of the main theorem.

\begin{proof}[Proof of Theorem~\ref{main}]
For the maps of Proposition~\ref{commpequa}, the set
\[
 \mathcal Z=\{(p,\xi)\in\mathcal D:(p,T_p(\xi))\in\mathscr P_I\}
\]
is compact and contains $I\times\{0\}$, where $I=[8,9]$.
We use the degree of Leray and Schauder with the identity orientation.
Compactness, relative openness of $\mathscr P_I$ and convergence of
zeros in \eqref{posictopo} imply
$\operatorname{dist}(\mathcal Z,G^{-1}(0)\setminus\mathcal Z)>0$.
Hence, for sufficiently small $\varepsilon>0$, the bounded open set
\[
\mathcal O=\{(p,\xi)\in I\times X:
       \operatorname{dist}((p,\xi),\mathcal Z)<\varepsilon\}
\]
has closure inside $\mathcal D$ with positive distance from its
complement and satisfies
$G^{-1}(0)\cap\overline{\mathcal O}=\mathcal Z\subset\mathcal O$.
All closures and boundaries here are relative to $I\times X$.
If $G(p_n,\xi_n)\to0$ on $\partial\mathcal O$, joint compactness of
$\mathcal C$ and $\xi_n=G(p_n,\xi_n)-\mathcal C(p_n,\xi_n)$ yield a
convergent subsequence whose limit is a boundary zero, a contradiction.
The homotopy theorem for variable domains
\cite[Theorem~2.1(iii), p.~198]{Mawhin99} therefore implies, for
$\Omega_p=\{\xi:(p,\xi)\in\mathcal O\}$,
\begin{equation}\label{degreequa}
 \deg(G(8,\cdot),\Omega_8,0)=\deg(G(9,\cdot),\Omega_9,0).
\end{equation}

Proposition~\ref{endpindi} and
Proposition~\ref{coerline} imply that $J_p:X\to Y$
is an isomorphism for $p=8,9$. Thus $0$ is an isolated zero of
$G(p,\cdot)$, whose derivative is $A_p=B_p^{-1}J_p$.
Lemma~\ref{negaelift} equates the total algebraic
multiplicity of the real eigenvalues $\mu>1$ of $\operatorname{Id}-A_p$ with
$\ind(Q_p;\mathcal H)$, using $\lambda=1-\mu$.
The local degree formula \cite[\S2, p.~197]{Mawhin99} yields
\[
 \deg_{\rm loc}(G(p,\cdot),0)=(-1)^{\ind(Q_p;\mathcal H)},\qquad
 \deg_{\rm loc}(G(8,\cdot),0)=-1,\quad
 \deg_{\rm loc}(G(9,\cdot),0)=1.
\]
If both integer parameters admitted only the standard solution in
$\mathscr P_I$, this would contradict
\eqref{degreequa}. Consequently, there exist
$q\in\{8,9\}$ and $(q,U)\in\mathscr P_I$ with $U\ne V_q$.

At this integer parameter, Lemma~\ref{intelift} provides
the smooth lift $u(re^{i\theta})=U(r,q\theta)$ with the required
equation, degree, normalization, symmetries and sector positivity.
Its potential energy is
\[
 \int_{\R^2}(1-|u|^2)^2\,dx=q\mathcal P_q(U)=2\pi q^2.
\]
If $u=\zeta V_q(\cdot-a)$ with $|\zeta|=1$, its unique zero and
rotational symmetry imply $a=e^{2\pi i/q}a$, hence $a=0$. 
The normalization at infinity implies $\zeta=1$, so $U=V_q$, a contradiction.
\end{proof}

\appendix
\section{Radial inequalities}\label{radialineq}

We prove the explicit radial bounds used in the Morse index argument.
The functions below satisfy strict differential inequalities, so the
quotient identity \cite[(2.5)]{DPLMWY} applies with the stated boundary data.

\begin{lemma}\label{algerboun}
For $p\in\{8,9\}$, let $F_p(R)=f_p(pR)$ and $x=R^2$. Then
\[
\begin{gathered}
 x(1-F_8(R)^2)<43/40\quad\text{for }x\ge7/6,\qquad
 x(1-F_9(R)^2)<9/8\quad\text{for }x\ge5/3,\\
 F_9(R)^2<x^5/9\quad\text{for }x>0,\qquad
 F_9(\sqrt{17/18})^2<1/13,\qquad F_9(\sqrt{5/4})^2<1/5.
\end{gathered}
\]
\end{lemma}
\begin{proof}
For $\mathcal L_p$ in \eqref{scalrequa}, choose
$\beta_8=1/4$, $\beta_9=1/3$ and
\[
 \Pi_8(s)=s^4+\frac{17}{16}s^3+\frac12s^2+\frac78s+2,\qquad
 \Pi_9(s)=s^4+\frac98s^3+\frac12s^2+\frac12s+2.
\]
Let $G_p(R)=(R^2-\beta_p)^2/\sqrt{\Pi_p(R^2-\beta_p)}$ for
$R>\sqrt{\beta_p}$, and $G_p=0$ otherwise. Direct differentiation
shows that $\mathcal L_pG_p<0$ on $R>\sqrt{\beta_p}$.
The extension is $C^1$, and the radial asymptotics imply
\[
 G_p(\sqrt{\beta_p})=0<F_p(\sqrt{\beta_p}),\qquad
 F_8-G_8=\frac1{32x}+O(x^{-2}),\qquad
 F_9-G_9=\frac1{16x}+O(x^{-2}).
\]
At an interior maximum $z=G_p/F_p\ge1$, the quotient identity
\[
 \frac1{p^2}(RF_p^2z')'
 =RF_p^2z\left[F_p^2(z^2-1)-\frac{\mathcal L_pG_p}{G_p}\right]
\]
would have a nonpositive left side and a positive right side.
Thus $G_p<F_p$. Direct algebra yields
$x(1-G_8^2)<43/40$ for $x\ge7/6$ and
$x(1-G_9^2)<9/8$ for $x\ge5/3$, proving the first two bounds.

For the upper bounds, take
$H(R)=R^6/\sqrt{R^{12}+R^{10}+R^8+4R^4+4}$.
Differentiation shows that $\mathcal L_9H>0$ on $(0,\infty)$.
Since $F_9/H=O(R^3)$ at zero and $F_9/H\to1$ at infinity,
the same quotient argument implies $F_9<H$.
For $x>0$, a direct calculation shows that
\[
 x^6+x^5+x^4+4x^2+4>9x.
\]
Evaluation of this polynomial divided by $x^6$ at $17/18$ and
$5/4$ yields values greater than $13$ and $5$, respectively.
These inequalities prove all three upper bounds.
\end{proof}

\Needspace{5\baselineskip}
The final lemma treats the last interval for $p=8$ in
Proposition~\ref{endpindi}.
\begin{lemma}\label{eightcomp}
For the scalar operator $\mathscr S_8$ associated with
\eqref{scalscomp}, let
$\eta(x)=(x-7/6)x^{-17/3}e^{-6/x}$ on $x>7/6$.
Then $\mathscr S_8[\eta(e^{2t})]>0$ for $t>\frac12\log(7/6)$.
\end{lemma}
\begin{proof}
The preceding lemma implies $\mathfrak p<43/40$. Since $\partial_{\mathfrak p}V<0$,
\[
 V\ge-\frac{13}{12}
       +\frac{(35/12)(x-13/12)}{x+11/6}.
\]
Substituting this lower bound and differentiating with
$\partial_t=2x\partial_x$, we find
$\mathscr S_8[\eta(e^{2t})]/\eta(e^{2t})>0$ for $x>7/6$.
Indeed, after multiplication by the positive denominator
$x^2(6x-7)(6x+11)$, the resulting polynomial is a sum of squares
with a positive constant term.
\end{proof}

\section*{Acknowledgements}
Hong-Ge Chen was supported by the National Natural Science Foundation of China (Grant Nos.~12671239 and~12571249) and the Postdoctoral Project of Hubei Province (Grant No.~2024HBBHXF095). Yong Liu was supported by National Key R\&D Program program of China 2022YFA1005400, the National Natural Science Foundation of China (Grant No.~12471204). Juncheng Wei was supported by the China National Key R\&D Program of China (Grant No.~2022YFA1005602) and the Hong Kong General Research Fund (projects ``New frontier in singular limits of nonlinear partial differential equations" and "On Fujita equations in the critical or supercritical regime"). Wen Yang was supported by the National Key Research and Development Program of China (Grant No.~2022YFA1006800), the National Natural Science Foundation of China
(Grant No.~12531010), the Science and Technology Development Fund of the Macao SAR (FDCT, Grant No.~0070/2024/\allowbreak RIA1), the Multi-Year Research Grants of the University
of Macau (Grant Nos.~MYRG-GRG2025-00051-FST and~MYRG-GRG2026-00034-FSC), and the University of Macau Development Foundation (Grant No.~TISF/2025/006/FST). The authors used AI-assisted tools, principally ChatGPT, for brainstorming, literature searches and manuscript editing. AI output was not treated as mathematical authority. The authors independently worked and verified all theorem statements,
proofs, calculations, constants, citations, and final prose, and take full responsibility for the contents of the paper.

\end{document}